\documentclass[Afour,sageh,times, doublespace]{sagej}

\usepackage{moreverb,url}

\usepackage{epstopdf}
\usepackage{graphicx}
\usepackage{subcaption}
\usepackage[title]{appendix}

\usepackage[colorlinks,bookmarksopen,bookmarksnumbered,citecolor=red,urlcolor=red]{hyperref}

\newcommand*\rfrac[2]{{}^{#1}\!/_{#2}}
\newcommand{\real}{{\mathbb{R}}}

\newtheorem{mydef}{\bf Definition}
\newtheorem{prop}{\bf Proposition}

\newtheorem{remark}{\bf Remark}

\newtheorem{theorem}{\bf Theorem}
\newtheorem{corollary}{\bf Corollary}

\newcommand\BibTeX{{\rmfamily B\kern-.05em \textsc{i\kern-.025em b}\kern-.08em
T\kern-.1667em\lower.7ex\hbox{E}\kern-.125emX}}

\def\volumeyear{2017}

\begin{document}

\runninghead{Da and Grizzle}

\title{Combining Trajectory Optimization, Supervised Machine Learning, and Model Structure for Mitigating the Curse of Dimensionality in the Control of Bipedal Robots}

\author{Xingye Da\affilnum{1} and Jessy Grizzle\affilnum{2}}

\affiliation{\affilnum{1} Mechanical Engineering Department, University of Michigan, Ann Arbor, MI, US \\
\affilnum{2} Electrical and Computer Engineering and Robotics Institute, University of Michigan, Ann Arbor, MI, US}

\corrauth{Xingye Da, Mechanical Engineering Department, University of Michigan, Ann Arbor, MI 48109, US}

\email{xda@umich.edu}

\begin{abstract}
	To overcome the obstructions imposed by high-dimensional bipedal models, we embed a stable walking motion in an attractive low-dimensional surface of the system's state space. The process begins with trajectory optimization to design an open-loop periodic walking motion of the high-dimensional model and then adding to this solution, a carefully selected set of additional open-loop trajectories of the model that steer toward the nominal motion. A drawback of trajectories is that they provide little information on how to respond to a disturbance. To address this shortcoming, Supervised Machine Learning is used to extract a low-dimensional state-variable realization of the open-loop trajectories. The periodic orbit is now an attractor of the low-dimensional state-variable model but is not attractive in the full-order system. We then use the special structure of mechanical models associated with bipedal robots to embed the low-dimensional model in the original model in such a manner that the desired walking motions are locally exponentially stable. The design procedure is first developed for ordinary differential equations and illustrated on a simple model. The methods are subsequently extended to a class of hybrid models and then realized experimentally on an Atrias-series 3D bipedal robot.
\end{abstract}


\maketitle

\section{Introduction and Problem Statement}

%
%
%
%
%
%
%
%

\begin{figure}
	\centering
	\includegraphics[width=1\columnwidth]{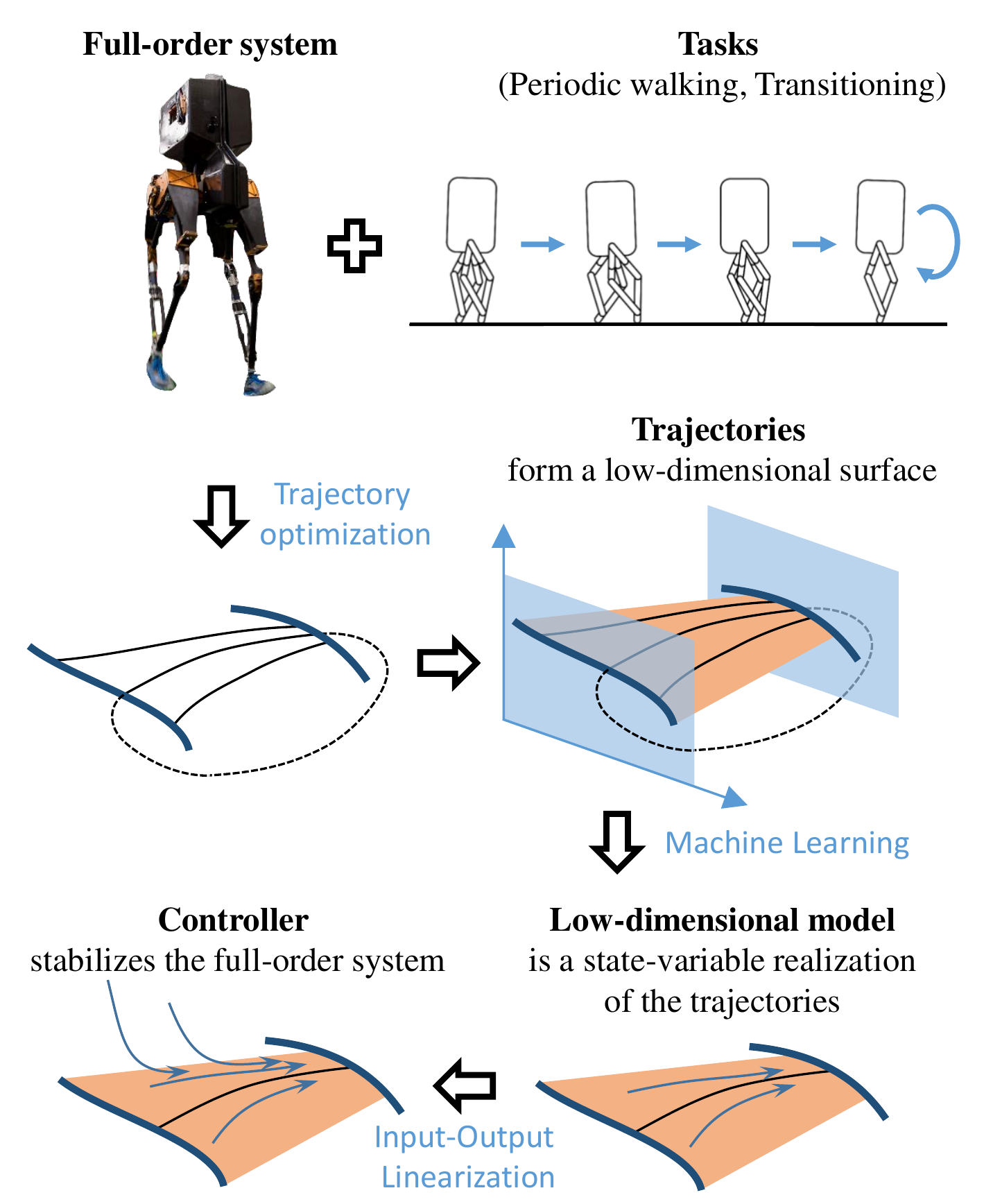}
	\caption{Our overall approach. The full-order model, the desired objectives, and physical constraints are combined into a trajectory optimization problem for designing a periodic gait. Using model structure or physical insight, a low-dimensional surface of initial conditions is selected for trajectory building, with the trajectories designed to approach the periodic orbit. If the trajectories form a low-dimensional surface, Supervised Machine Learning is used to a extract a vector field from the data that realizes the trajectories. System structure is then used again to render the low-dimensional model (surface and vector field) invariant and attractive.}
	\label{fig:OurApproach}
\end{figure}

We seek to design controllers for high degree-of-freedom (DoF) bipedal robots with several degrees of underactuation (DoU), or, if the robot is ``fully actuated'', we wish to take into account the limited ability of ankle torques to affect the overall evolution of the robot. The paper will focus on the tasks of walking stably forward, backward, or in place, and transitioning among such motions. We want the gaits to be dynamic in the sense that they can use the full capability of the robot regarding speed, terrain type, and other forms of agility. Moreover, of course, we need to embed the controller on the robot for real-time implementation. Our unique approach to this well-studied problem is illustrated in Figure~\ref{fig:OurApproach}.

\subsection{Proposed Approach to Controller Design}

We begin with a full-order dynamical model of the robot, expressing its many degrees of freedom and actuation capability, and a simplified representation of its contact with the environment so that the overall model can be expressed as a system with impulse effects, a special class of hybrid models. To overcome the obstructions imposed by high-dimensional models, we first seek to embed a stable walking motion in a low-dimensional surface of the system's state space. Subsequently, we seek to stabilize the gait in the full-order model by rendering the surface attractive.

\begin{figure}
	\centering
	\includegraphics[width=0.8\columnwidth]{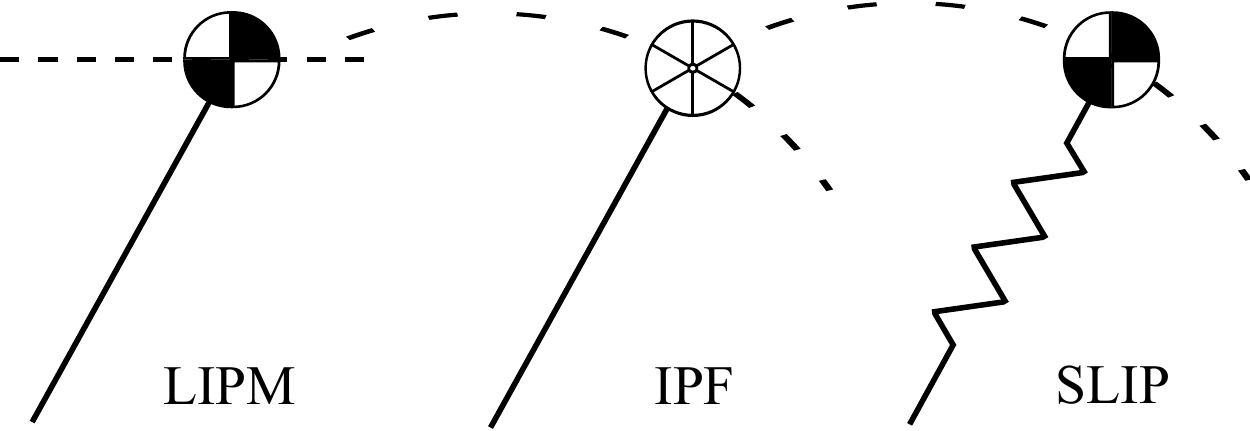}
	\caption{Adapted from \cite[Figure~5]{grizzle2014models}. Three low-dimensional models that are frequently used as approximate representations of walking robots. From left to right: the Linear Inverted Pendulum (LIP) lumps the mass of the robot at a point moving at a constant height and assumes massless legs;
		the Inverted Pendulum with Flywheel (IPF) relaxes the assumption on constant height and adds a flywheel to account for internal angular momentum; and
		the Spring-Loaded Inverted Pendulum (SLIP) adds a spring to model a robot's legs as a massless pogo stick. There is no obvious way to embed these low-dimensional models into the full model of a robot.}
	\label{fig:biped-models}
\end{figure}

The most common approach in the literature to getting around the high dimension of the model is to represent the walking task through the dynamics of a low-dimensional inverted pendulum (e.g., LIP, SLIP or others in Figure~\ref{fig:biped-models}), which equipped with a foot-placement strategy for stability \cite{KAYAKO92, RA86a,RA86b, PRTE06}. The robot is then controlled in such a way that its center of mass approximately follows the target dynamics. The many challenges associated with this more common approach include: achieving stable solutions in the full model; exploiting the full capability of the machine, especially in light of physical constraints of the hardware or environment; deciding how to associate the states of the low-dimensional pendulum with the full-order system; and finally, even deciding upon the appropriate pendulum model for a given task is not evident: what is the correct model for turning while stepping off a platform?


For these reasons, we do not rely on a pre-specified pendulum model to encode the walking motion. In this introduction, we ask that the reader allow us to sketch the main ideas of our approach without worrying too much about technicalities. Later developments will be more formal.

For the sake of simplicity, let's pretend the model of the robot can be captured by an ordinary differential equation,
\begin{equation}
\label{eq:PretendModel}
\dot{x}=F(x,u),
\end{equation}
with state variables $x \in {\cal X}$ and control inputs $u\in {\cal U}$. The design process begins with construction of a periodic solution meeting relevant constraints. Denote the period by $T_p>0$ and the initial condition by $\xi^*$. The next step is to make an initial selection of a low-dimensional set, let's call it $Z_0 \subset {\cal X}$, such that $\xi^*\in Z_0$. Now begins the real work; we seek to design open-loop trajectories of the full-order model
\begin{equation}
\label{eq:PeriodicBehavior}
(x_\xi(t), u_\xi(t)),~ 0\le t \le T_p,~\xi \in Z_0
\end{equation}
that over the interval $[0, T_p]$, ``approach'' the periodic solution. Specifically, for some $0 \le c < 1$, and for each $\xi \in Z_0$, we have
\begin{equation}
\label{eq:BoundaryCondition}
\varphi_\xi(T_p) \in Z_0 ~~\text{and}~~
||\varphi_\xi(T_p) - \xi^*||  \le c ||\xi - \xi^*||.
\end{equation}
After designing the trajectories, we seek to construct  a low-dimensional subsystem that realizes them, namely,
\begin{equation}
\label{eq:PeriodicRealization}
\begin{aligned}
\dot{z}(t)&=G(t,z)\\
x(t)&=H(t,z),
\end{aligned}
\end{equation}
with $z \in Z \subset {\cal X}$, such that (a) for $\xi \in Z_0$,
$$z(0)=\xi \Rightarrow x_\xi(t) = H(t, z(t)),$$ and (b) the periodic motion is a ``locally exponentially stable output'' of the model.  If this can be done, we would argue that \eqref{eq:PeriodicRealization} is a more desirable target model than a typical pendulum because the target has been constructed directly from the full-order model and its ``specification'', that is, the constraints imposed when designing the trajectories.

\begin{figure}
	\centering
	\includegraphics[width=0.6\columnwidth]{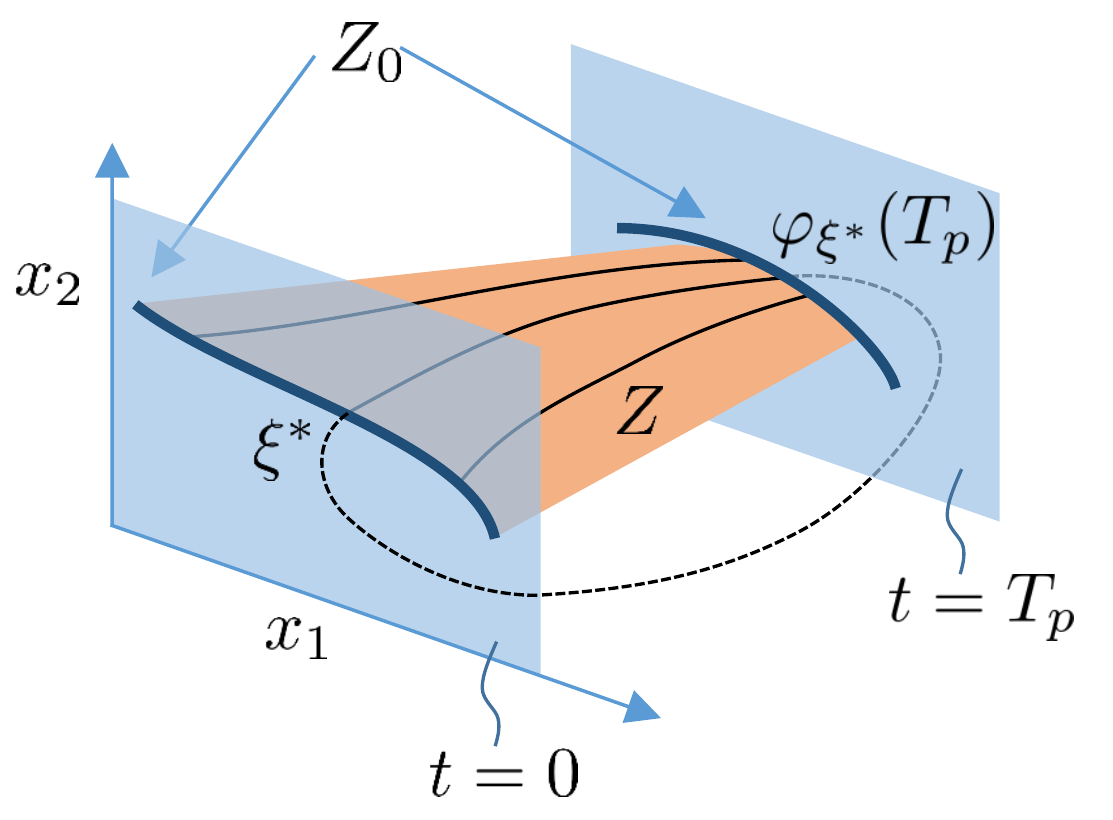}
	\caption{Based on model structure, the system's state is decomposed into $x=(x_1, x_2)$, where the dimension of $x_1$ is much smaller than the dimension of $x_2$. The surface $Z_0$, shown here as a line, is the set of initial conditions used to build a set of trajectories that will fill out the surface $Z$. By construction, this model, if it exists, will be easy to embed in the robot's full state space.}
	\label{fig:ControlConcept}
\end{figure}

To turn this into a viable feedback design process for the original system, we have to address the following issues:
\begin{enumerate}
	\item[(i)] How to compute the low-dimensional model \eqref{eq:PeriodicRealization} for realistic bipedal robots and the surface\endnote{$Z_0$ only contains the initial conditions for the trajectories. The evolution of the trajectories will determine $Z$. The vector field and output map in \eqref{eq:PeriodicRealization} will be computed through a combination of Supervised Machine Learning and model structure.} $Z \supset Z_0$ on which it is defined?
	
	\item[(ii)] As with any low-dimensional target model, how to embed it in the full-order model with stability? In other words, how to design a feedback for the original system \eqref{eq:PretendModel} that does two things: (a) creates an invariant surface $Z$ in its state space with restriction dynamics given by \eqref{eq:PeriodicRealization}, thereby encoding the desired stable walking motion; and (b) solutions of the closed-loop system starting near the surface asymptotically converge to the surface, thereby realizing the walking motion in a stable manner in the overall system.
\end{enumerate}

The remainder of the paper is dedicated to addressing these challenges for a class of hybrid models and tasks of interest to bipedal locomotion. Section~\ref{sec:MainIdeas} develops the basic ideas in the simpler setting of ordinary differential equations. The results are of independent interest for tasks such as rising from a sitting position or standing in place. Trajectory optimization is used to generate the low-dimensional set of open-loop trajectories \eqref{eq:PeriodicBehavior} that includes a metric for attractivity to a periodic solution, a family of periodic solutions, or transitions among such solutions. Model structure and Supervised Machine Learning are proposed as a means to extract functions from the open-loop trajectories to build the low-dimensional model \eqref{eq:PeriodicRealization}. Finally, the appeal is made once again to model structure to embed the low-dimensional model in the full-order model while guaranteeing local exponential attractivity.

Section~\ref{sec:PendulumOnCart} illustrates the design process on the well-known inverted pendulum on a cart. This will allow the reader to explore the method on a simple model. Section~\ref{sec:HybridSystemControl} develops the results for hybrid models, preparing the ground for the simulations and experiments reported in Section~\ref{sec:biped} for the bipedal robot of Figure~\ref{fig:OurApproach}.  The discussion and conclusions are given in Section~\ref{sec:Conclusions} and all proofs are given in Appendix~\ref{sec:Proofs}.

\section{Related Literature}

In the following, related work in the legged robotics literature is summarized and contrasted to the work in the present paper. Appendix~\ref{app:Relations} provides a more technical comparison of several nonlinear control methods, namely, backstepping, zero dynamics, and immersion and invariance.

\subsection{Online Optimization} One of the earliest uses of online optimization for bipedal walking was done on a \textit{simulation model} of the planar robot RABBIT \cite{AzPoEs02, AzPoEs04}. The controller could only be done in simulation because its computation time was approximately forty times slower than the duration of a step. More recently, Model Predictive Control (MPC) was applied in the \textit{DARPA Virtual Robotics Challenge} \cite{ErLoTaKuKoTo13}. In that work, a real ``simulation-time'' implementation on a full-order model of Atlas was achieved through the use of a novel physics engine and a relaxed contact map. MPC was applied on \cite{KuDeFa16} to a simplified model of Atlas that captured the kinematics and centroidal dynamics; this resulted in walking at 0.4~m/s. On a planar biped, higher walking speeds from 0.43~m/s to 0.97~m/s were achieved in \cite{Hereid2016Online} using online Partial Hybrid Zero Dynamics (PHZD) gait generation. Average computational time was 0.5~s. To stabilize a robot that has to make and break multiple contacts with the environment, a piecewise affine approximation of a nonlinear hybrid system was controlled using a mix of online approach and explicit MPC \cite{MarcucciTedrake2017}, or a Piecewise-Affine Quadratic Regulator \cite{HanTedrake2017}. These optimizations are more effective than the purely online MPC, but the model so far has not taken the centroidal angular momentum into account.


\subsection{Pendulum Models}

The pendulum models illustrated in Figure~\ref{fig:biped-models} are ubiquitous means in the bipedal robotics literature to reduce the online computational burden. The LIP model is especially prevalent for the design of flat-footed walking gaits based on the Zero Moment Point criterion \cite{Kajita1991,Yamaguchi99,Takanishi06_WABIAN2,SAWAAO02,robots_HRP4,robots_KHR3,PFLOGI02}. The optimization or closed-form compute the CoM trajectories and swing foot positions on the reduced-order models. A low-level controller and inverse kinematics then realize these on the full-order model or robot. Recent experimental uses of this approach can be found in \cite {PRKODBRECOJONE12, KrEnWiOt12, FaPoAtIj14, ReHuJoPeVaJoAbHu15}.

The bottom line, however, is that when a pendulum model is pre-specified as a template \cite{FUKO99} the full order model needs to compromise its achievable motions to follow the template. Moreover, for each different task of the robot, such as walking or running, the designer is faced with the selection of the ``best'' target model. In our approach, a low-dimensional model is generated from the full-order model and the task. It is dynamically feasible and uses the full capability of the robot to accomplish the task.


\subsection{Pre-computed Gaits}

A means to get around the limitation of online computation is to pre-compute a set of controllers and design a control policy to ``stitch'' them together. A policy that switches the task (target walking speed, running vs. walking, stairs vs. flat ground) is employed in \cite{SrPaPoGr13,Martin2014Design,PoHeAm13}. Finite-state machines or motion primitives are used in \cite{PaRaGr13, SaBy15,MaUm14} for rough terrain and in \cite{ApLeBu15} for reducing settling time. Interpolation among gaits has been used to create a continuous family of gaits in \cite{EmViGr16, DaHaHaGrGr16, NguyenRSS2017, NguyenWAFR2017}. Transient trajectories that approach the nominal periodic orbit were added in \cite{LiAtSu13, DaHartleyGrizzle2017} to enlarge the basin of attraction. The current paper provides a formal mathematical framework for the work in \cite{DaHartleyGrizzle2017} and increases its applicability.


\subsection{Hybrid Zero Dynamics}

The work in the present paper is related to the method of virtual constraints and hybrid zero dynamics (HZD) proposed in \cite{WEGRKO02}. In that method, a monotonically increasing function of a robot's generalized coordinates $q$ is first identified, often denoted by $\theta(q)$, and then a family of virtual constraints of the form
\begin{equation}
\label{eq:VCs}
y = h_0(q) - h_d(\theta(q), \alpha)
\end{equation}
are posted, where $\dim(y)=\dim(u)$, the number of inputs, $h_0(q)$ represents quantities to be regulated, such as knee angles and hip angles, and $h_d(\theta(q), \alpha)$ is a vector of splines representing the two be determined desired evolution of $h_0(q)$. A parameter optimization problem is posted to select the values of $\alpha$ (if they exist) so that $y\equiv0$ along a periodic solution of the model, torque bounds are met, as are ground contact forces. If the outputs $y$ have vector relative degree two \cite[pp.~119]{WGCCM07}, the robot model with outputs \eqref{eq:VCs} is input-output linearizable, and hence a feedback controller can be designed that drives the virtual constraints asymptotically to zero. If the surface defined by the outputs being zero is invariant in the hybrid model of the robot, then \cite[Ch.~5]{WGCCM07} provides strong stability theorems for the closed-loop system.

While this theory has been successfully implemented on many robots \cite{ Ames2012Dynamically, Ames2014Human, Buss2014Preliminary, Chevallereau2003RABBIT, Hereid2014Dynamic, Martin2014Design, KOPAPOG210, Sreenath2012Design, Reher2016Algorithmic}, lower-limb prostheses \cite{GrLeHaSe14, aghasadeghi2013learning, Zhao2015hybrid} and exoskeleton \cite{Agrawal2017First}, it has important limitations that the current paper overcomes:
\begin{itemize}
\item In basic HZD, only one optimization is done, namely the determination of the periodic orbit. Hence, only that solution is guaranteed to be feasible. Here, we build a low-dimensional surface of feasible solutions.
    	\item The stability theorems in \cite{WGCCM07} only apply to robots with one degree of underactuation. Here, multiple degrees of underactuation can be handled. Moreover, unlike the method proposed in \cite{BuHaGrGr16, Kaveh2016BMI} which works with the linearization of the Poincar\'e map to achieve local exponential stability for robots with more than one degree of underactuation, here, at least on the low-dimensional surface, large perturbations away from the periodic orbit can be included in the controller design, providing better robustness.
    \item For robots with one degree of underactuation, the stability mechanism in \cite[Ch.~5.4]{WGCCM07} is through energy loss at impacts, similar to the stability proofs for passive robots walking down a slope. Here, more general stability mechanisms in bipedal robots \cite{PRTE06}, such as posture adjustments through foot placement and knee bend, automatically arise.
        \item Only gaits for which a monotonic variable can be identified can be treated with the current HZD method. Here, a much richer set of locomotion primitives can be realized, such as stepping in place or transitioning from walking forward to walking backward. As in \cite{WANGChevallereau2011}, \cite{ReCoHeHuAm16}, \cite{DaHaHaGrGr16}, the feedback is allowed to be time-dependent, enriching the set of possible closed-loop solutions.
            \item Similar in spirit to \cite{GriffinIJRR2016}, which added conjugate momenta of the underactuated coordinates into the virtual constraints, this paper can include velocity as well as positions in the low-dimensional surface that eventually defines a generalized hybrid zero dynamics manifold; see also \cite{HartleyGrizzleCCTA2017}.
\end{itemize}

\section{Presentation of Main Ideas}
\label{sec:MainIdeas}

For the class of robot problems of interest to us, optimal gaits can now be computed in minutes \cite{Jo2014, HeCoHuAm16}, but not in tens of milliseconds, which is what would be required for online use. In the simpler setting of a non-hybrid system, this section develops our main ideas for mitigating the curse of dimensionality in optimization-based controller design. Section~\ref{sec:MPCstyleController} provides the first example of conditions for a family of open-loop trajectories of a model from which a realization can be extracted and its equilibrium will be locally exponentially stable; see also \cite{co94,schurmann2017convex}. The process of building the realization from the trajectories is based on regression, namely Supervised Machine Learning. The size of the model that can be treated with these initial results is limited by both the number of optimizations it takes to create the family of open-loop trajectories and the number of features that can be included in the Supervised Machine Learning. Section~\ref{sec:StabilizingReducedOrderModel} extends the design process to building reduced-order target model for a high-dimensional model. Importantly, the design is based on a far smaller set of open-loop trajectories.  Sections~\ref{sec:InnerOuterLoop} and \ref{sec:ExtendedModels} then provide conditions for embedding the target model in the full-order model such that the origin of the full-order model is locally exponentially stable. The proofs of the results developed in the section are given in Appendix~\ref{sec:Proofs}. The relationship to other controller design methods is addressed in Appendix~\ref{app:Relations}.

In presenting the main ideas, we will deliberately organize them as a design philosophy. We choose to let the user rely on his or her wits to meet our conditions, rather than muddying the waters with a set of highly technical sufficient conditions that no one would ever check. We know as well as the readers know that optimization problems are very tricky: it is easy to paint oneself into a corner that only yields non-smooth solutions. On the other hand, many problems, such as the examples worked out in the paper, seem to have very nicely behaved solutions. We are confident that we did not cherry-pick the only nice problems and that the readers will find a host of further interesting examples.

\subsection{Model Assumptions}

To keep the connection to stabilizing periodic gaits in bipeds, we consider a periodically time-varying nonlinear system with equilibrium point at the origin,
\begin{equation}
\label{eq:PeriodicallyTimeVaryingSystem}
\dot{x}=f(t,x,u).
\end{equation}
The simple coordinate transformation required for shifting a periodic solution of a nonlinear model to the origin is provided in \cite[pp.~147]{KHA02}. An equilibrium point is treated as a special case of a periodic orbit where the period can be any number $T_p>0$; in particular, when discussing periodicity, $T_p$ is \textit{not required to be a fundamental period}.

The ODE \eqref{eq:PeriodicallyTimeVaryingSystem} is assumed to satisfy the following conditions.
\begin{enumerate}
\setcounter{enumi}{0}
\renewcommand\labelenumi{{\bf A-\theenumi}}
\item \label{itm:A1} $f:[0, \infty) \times \real^n \times \real^m \to \real^n$ is locally Lipschitz continuous in $x$ and $u$, piecewise continuous in $t$, and there exists $T_p>0$ such that $\forall~(t,x,u) \in [0, \infty) \times \real^n \times \real^m$, $f(t+T_p,x,u)=f(t,x,u)$.
\item \label{itm:A2}  $f(t,0,0)=0$ for all $t\ge 0$.
 \item \label{itm:A3}  The user has selected an open ball about the origin, $B\subset \real^n$, a positive-definite, locally Lipschitz-continuous function $V:B \to \real$, and constants $0<\alpha_1 \le \alpha_2$ such that, $\forall~x\in B$
 $$ \alpha_1 x^\top x \le V(x) \le \alpha_2 x^\top x.$$
\mbox{ } \hfill \(\blacksquare\)
\end{enumerate}

As alluded to above, the reader is encouraged to view the assumptions made throughout this section as requirements to impose on an open-loop trajectory optimizer. We have found them straightforward to meet when using the optimizers of \cite{Jo2014, HeCoHuAm16}. In many cases, the positive-definite function indicated in A-\ref{itm:A3} comes ``for free'' from the optimization problem used to compute trajectories; this is standard in model predictive control \cite{mayne2000constrained}. Because the Lyapunov condition in A-\ref{itm:A4} below can also be included as a constraint in the trajectory generation process, the user has much freedom in its selection, even something as simple as the $2$-norm squared could be used.

\subsection{Extracting a Feedback from Open-loop Trajectories}
\label{sec:MPCstyleController}

Two feedback controllers can be constructed from the following solutions of the model \eqref{eq:PeriodicallyTimeVaryingSystem}. The obvious relation to MPC is discussed in the remarks following Prop.~\ref{prop:XDaMPCStability}.

\begin{enumerate}
\setlength{\itemsep}{.15cm}
\setcounter{enumi}{3}
\renewcommand\labelenumi{{\bf A-\theenumi}}

  \item \label{itm:A4} There is a constant $0 \le c < 1$,  such that, for each initial condition $\xi \in B$, there exists a continuous input $u_\xi:[0, T_p] \to \real^m$ and a corresponding solution of the ODE, $\varphi_\xi: [0, T_p]  \to \real^n$ satisfying $\varphi_\xi(T_p) \in B$, and
\begin{equation}
\label{eq:LyapunovConstraintA4}
V(\varphi_\xi(T_p)) \le c V(\xi);
\end{equation}
moreover, for $\xi=0$, $u_\xi(t) \equiv 0$.
For clarity, solutions are taken in the sense of equation (C.2) of \cite[pp.~657]{KHA02}, namely
\begin{equation}
\label{eq:DefSolutionODE}
\varphi_\xi(t) = \xi + \int_{0}^t f(\tau,\varphi_\xi(\tau),u_\xi(\tau)) d \tau.
\end{equation}
\mbox{ } \hfill \(\blacksquare\)
\end{enumerate}


A $T_p$-periodic Continuous-Hold (CH) feedback is defined by periodic extension of $u_\xi(t)$, namely,
\begin{equation}
\label{eq:pwcontfeedback}
u^{ch}(t,\xi)=u_\xi(\hat{t}),~ \hat{t}=t~{\rm mod}~ T_p.
\end{equation}
Jumps are allowed at multiples of the period, with continuity taken from the right. Due to the reset or hold-nature of the above feedback, the stability of solutions of
\eqref{eq:PeriodicallyTimeVaryingSystem} in closed-loop with \eqref{eq:pwcontfeedback}
should be studied as a sampled-data system, that is, the solutions should be evaluated at times $t_k=kT_p$. We will not analyze its stability, however, because it will clearly perform poorly in the face of perturbations occurring between samples, where the system is open loop.

We proceed directly instead to a feedback controller that allows continual updates in the state variables, and yet, under certain conditions, can be built from the open-loop trajectories given in Assumption A-\ref{itm:A4}. To understand the feedback controller, a \textit{thought experiment} is helpful: Suppose at time $t_0=0$ the system's initial state value is $\xi^0$, and the continuous-hold feedback \eqref{eq:pwcontfeedback} is being applied. Then the system is evolving along the trajectory $\varphi_{\xi^0}(t)$. Suppose subsequently at time $0 < t_d < T_p$, an ``impulsive disturbance'' affects the system, displacing the system's state to a value $x(t_d)\neq \varphi_{\xi^0}(t_d)$.  What input might be applied, given the information in A-\ref{itm:A4}? If there exists a $\xi^d \in B$ such that $x(t_d)=\varphi_{\xi^d}(t_d)$, then applying the input $u_{\xi^d}(t)$ for $t_d \le t <T_p$ will move the system toward the equilibrium in the sense that $V( \varphi_{\xi^d}(T_p))\le cV(\xi^d)$. The next result builds on this idea; see also \cite{co94,schurmann2017convex}.

%

\begin{prop}
\label{prop:XDaMPCStability}
Assume the open-loop system \eqref{eq:PeriodicallyTimeVaryingSystem} satisfies Assumptions A-\ref{itm:A1} to A-\ref{itm:A4}.
Assume in addition there exists an open set $B^e \supset B$ and a feedback $$\mu:[0, \infty) \times B^e \to \real^m$$ that is piecewise continuous in $t$, $T_p$-periodic, locally Lipschitz continuous in $x$, and, such that, for $0\le t < T_p$ and $\xi \in B$,
\begin{equation}
\label{eq:LearningConditionA}
\mu(t,\varphi_{\xi}(t))=u_\xi(t).
\end{equation}
Then the origin of the closed-loop system,
\begin{equation}
\label{eq:ClosedLoopCUFullStateFeedback}
\dot{x}=f^{cl}(t,x):=f(t,x,\mu(t,x)),
\end{equation}
is locally exponentially stable, uniformly in $t$, and the trajectories in A-\ref{itm:A4} are solutions of \eqref{eq:ClosedLoopCUFullStateFeedback}. Said another way, \eqref{eq:ClosedLoopCUFullStateFeedback} is a realization of the trajectories in A-\ref{itm:A4}.
\mbox{ } \hfill \(\blacksquare\)
\end{prop}


\begin{table}
	\centering
		\begin{tabular}{c|c||c}
			\multicolumn{2}{c||}{\bf Features } &   { \bf Labels}  \\
			\hline
			$t_j$& $x^{j, i}=\varphi_{\xi^i}(t_j) $& $\mu^{j, i}=u_{\xi^i}(t_j)$\\
			\hline
        & & \\
			$t_0=0$ & $x^{0,1} $  & $\mu^{0,1}$ \\
 & & \\
			$t_0=0$ & $x^{0,2} $  & $\mu^{0,2}$ \\
			$\vdots$ & $\vdots$& $\vdots$\\
			$t_0=0$ & $x^{0,M} $  & $\mu^{0,M}$ \\
 & & \\
			$t_1$ & $x^{1,1} $  & $\mu^{1,1}$ \\
			$\vdots$ & $\vdots$& $\vdots$\\
			$t_1$ & $x^{1,M} $  & $\mu^{1,M}$ \\
			$\vdots$ & $\vdots$& $\vdots$  \\
			$t_N=T_p$ & $x^{N,1} $  & $\mu^{N,1}$ \\
			$\vdots$ & $\vdots$& $\vdots$  \\
			$t_N=T_p$ & $x^{N,M} $ & $\mu^{N,M}$ \\
			\end{tabular}
		\caption[]{Conceptual arrangement of the data in A-\ref{itm:A4} from which a controller satisfying \eqref{eq:LearningConditionA} may be determined by ``regression''. In practice, not only time must be discretized, but also the initial conditions $\xi^i \in B$. This is where the ``curse of dimensionality'' rears its ugly head. Placing ten points per dimension leads to $10^n$ optimizations to compute, which quickly becomes impractical.}
		\label{tab:How2DoRegression}
\end{table}

Returning to \eqref{eq:PeriodicRealization} in the Introduction, the key point is that when a function can be found that is compatible with the open-loop trajectories in A-\ref{itm:A4} in the sense that the ``learning condition '' \eqref{eq:LearningConditionA} is satisfied, then \eqref{eq:ClosedLoopCUFullStateFeedback} provides an ``exponentially stable realization'' of the trajectories; this idea will be extended to a lower-dimensional system in the next subsection. Secondly, if at time $0\le t \le T_p$ there exists $\xi$ such that $x(t)=\varphi_{\xi}(t)$, the value of the ``learned function'' is being set to $$\mu(t,x(t)):=u_\xi(t).$$
In practice, the trajectories will only be computed for a finite grid of initial conditions $\xi^j$, $j\in J$. Hence, \eqref{eq:LearningConditionA} is an \textit{interpolation of the data} from the trajectory optimizations. For $x(t)\in B^e$ for which there does not exist $\xi$ such that $x(t)=\varphi_{\xi}(t)$, the function in \eqref{eq:LearningConditionA} is an \textit{extrapolation of the data}. In \cite{co94}, the fitting of a function to trajectory data was done \textit{in principle} in closed form and by hand; here, the fitting is being done with Supervised Machine Learning. Moreover, the learning algorithms available today provide easy tools for checking the quality of a fit and hence for checking how closely a function was found that meets the learning condition \eqref{eq:LearningConditionA}.

\begin{remark}
Because the solutions in A-\ref{itm:A4} will be computed via a trajectory optimization algorithm, it is useful to understand how the assumptions on $\mu$ relate to requirements on the trajectories.
\begin{enumerate}


\item[(i)] Suppose $T_h > T_p$ and that for $\xi \in B$, ${u^o_\xi:[0, T_h] \to \real^m}$ minimizes a cost function of the form
\begin{equation}
\label{eq:CostFunctionOne}
J(\xi)= \min_{u} \int_{0}^{T_h} L(  \varphi_\xi(t), u_\xi(t)) dt + N(\varphi_\xi(T_h))
\end{equation}
where, as before, $\varphi_\xi(t)$ is the solution of  \eqref{eq:PeriodicallyTimeVaryingSystem} with initial condition $\xi$ at $t_0=0$, and suppose furthermore that $u_\xi:[0, T_p] \to \real^m$ is the restriction of $u^o_\xi$ to $[0, T_p]$, that is,
$$\left. u_\xi =  u^o_\xi \right|_{[0, T_p]}.$$
By the principle of optimality, for $0 \le t_0 < T_p$,
$$\left. u^o_\xi \right|_{[t_0, T_h]} $$
is a minimizer of
\begin{equation}
\label{eq:CostFunctionTwo}
J(x_0)= \min \int_{t_0}^{T_h} L( \varphi_\xi(t), u_\xi(t)) dt + N(\varphi_\xi(T_h)),
\end{equation}
where, $\varphi_\xi(t)$ is the solution of \eqref{eq:PeriodicallyTimeVaryingSystem} with initial condition $\varphi_\xi(t_0)$ at $t_0$.
Hence, the condition \eqref{eq:LearningConditionA} can be interpreted as arising from an MPC-style controller with a shrinking horizon, $[t_0, T_h]$, for $0 \le t_0 \le T_p$ and fixed final-time $T_h$. This control strategy is visualized in Figure~\ref{fig:optimialPrincipal}. The condition \eqref{eq:LearningConditionA} comes from the shrinking of the optimization horizon; it will be essential in allowing a judiciously chosen set of open-loop trajectories to be realized with a low-dimensional state-variable model.

\begin{figure}
	\centering
	\includegraphics[scale = 0.5]{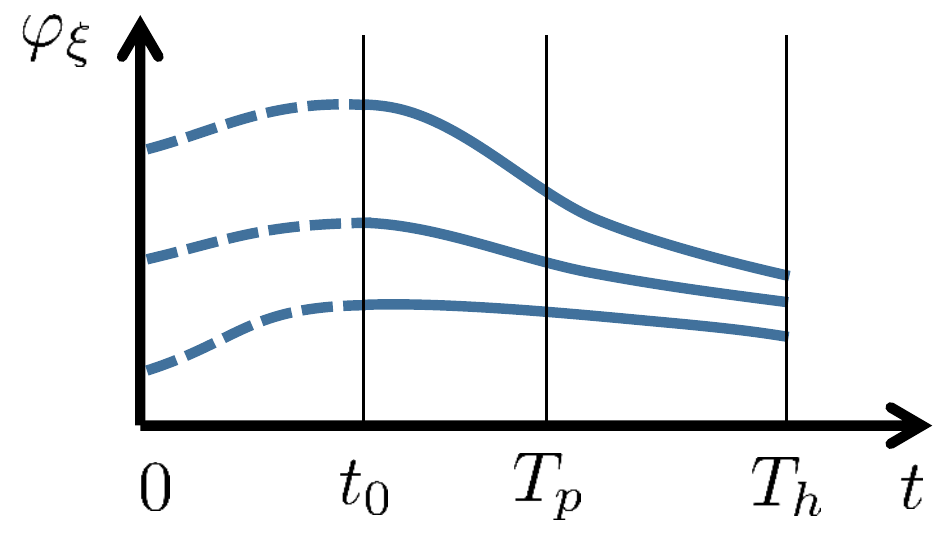}
	\caption{Principle of Optimality. If the system is initialized at $\varphi_\xi(t_0)$ and the cost function is modified from \eqref{eq:CostFunctionOne} to \eqref{eq:CostFunctionTwo}, then $\varphi_\xi:[t_0, T_h] \to \real^n$ is optimal.}
	\label{fig:optimialPrincipal}
\end{figure}

\begin{figure}
	\centering
	\includegraphics[scale = 0.5]{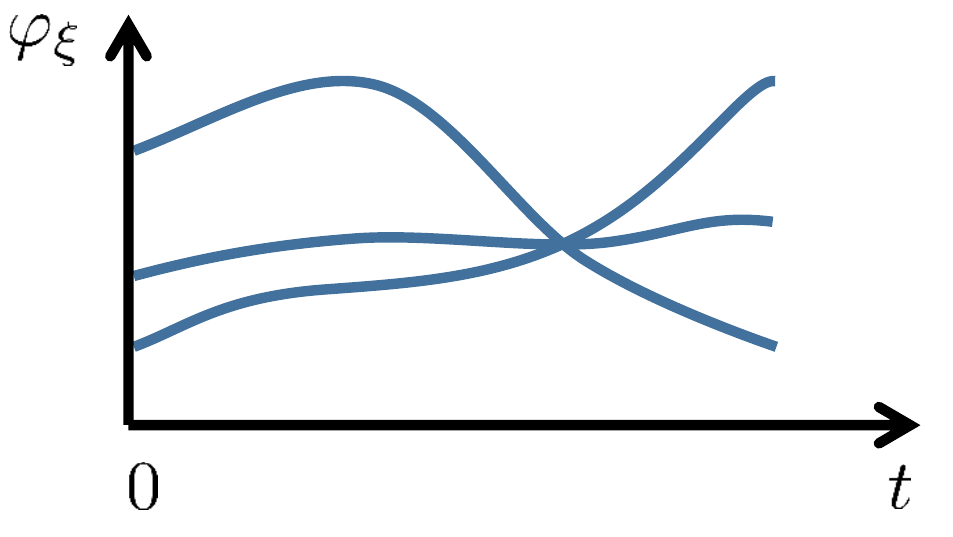}
	\caption{This shows the trajectories crossing one another, which means that the mapping $\Psi_t:B \to \real^n$ in \eqref {eq:PSIt} is not injective at certain moments of time. In this case, a feedback function cannot be extracted from the data.}
	\label{fig:tangleDemostration}
\end{figure}

\item[(ii)] Supervised Machine Learning will be used to extract the function $ \mu(t, x) $ in Prop.~\ref{prop:XDaMPCStability} from the trajectories and control inputs given in A-\ref{itm:A4} The method is sketched in Table~\ref{tab:How2DoRegression}. An example is given in Sect.~\ref{sec:PendulumOnCart}.

\item[(iii)] The local Lipschitz continuity of $\mu$ imposes conditions on the solutions given in A-\ref{itm:A4}. Indeed, for each $t \in [0, T_p]$, the mapping
\begin{equation}
\label{eq:PSIt}
\Psi_t:B \to \real^n,~\text{by} ~\Psi_t(\xi):=\varphi_\xi(t)
\end{equation}
 must be injective. This follows by the Gronwall-Bellman inequality \cite[pp.~651]{KHA02}; see also \cite[Exercise~3.17]{KHA02}. Hence, the optimization problem must be set up so as to avoid the existence of trajectories that cross one another, which can easily occur as shown in Figure~\ref{fig:tangleDemostration}. For example, if the user selected $T_h=T_p$ and imposed that the origin be attained at $T_h$, that is, dead-beat control, then a locally Lipschitz continuous $\mu$ would not exist.


 \item[(iv)] Conditions are known under which the value function \eqref{eq:CostFunctionOne} meets the Lyapunov conditions in A-\ref{itm:A3} and A-\ref{itm:A4}. Roughly speaking, they require that the terminal weight $N$ either be replaced with finite-time convergence to the origin or, the terminal weight be selected as $ \beta N(x(T_h))$, where $N(x)$ is positive definite and $\beta>0$ is sufficiently large. Hence, in practice, the Lyapunov constraint can be replaced by careful formulation of the trajectory optimization problem. In our limited experience, it is never an active constraint.

 \item[(v)] There is a long history of work in the nonlinear control literature that relates asymptotic controllability to an equilibrium point and the existence of stabilizing feedback controllers. The reader is referred to \cite{co94, coron1995feedback,clarke1997asymptotic,shim2003asymptotic} and references therein. The methods employed are not nearly as constructive as the present paper.
\end{enumerate}

\end{remark}

\subsection{Building a Reduced-Order Target Model}
\label{sec:StabilizingReducedOrderModel}

To begin the construction of a reduced-order target model as in \eqref{eq:PeriodicRealization}, we now assume that the system \eqref{eq:PeriodicallyTimeVaryingSystem} is decomposed in the form
\begin{equation}
\label{eq:TwoBlockSystem}
\begin{aligned}
\dot{x}_1&= f_1(t,x_1,x_2)\\
\dot{x}_2&= f_2(t,x_1,x_2,u),
\end{aligned}
\end{equation}
where $x_1\in \real^{n_1}$ and $x_2\in \real^{n_2}$.   For clarity of exposition, the input is assumed not to appear in $f_1$; the changes required to include inputs in $f_1$ are given in Sect.~\ref{sec:ExtendedModels}. Assumptions A-\ref{itm:A1} through A-\ref{itm:A3} are assumed to hold for \eqref{eq:TwoBlockSystem}.

Because of how the decomposition will arise in the case of bipeds, we think of the $x_1$-states as the ``weakly actuated part'' of the system and the $x_2$-states as the ``strongly actuated part'' of the system. With the model expressed in this form, it is clear that the $x_2$-states are virtual controls for the $x_1$-states. We will continue to build open-loop trajectories by the full-order model, except now the trajectories will be computed for a \textit{reduced set of initial conditions defined by the $x_1$-subsystem}.


\begin{mydef}
\label{def:InsertionFunction}
An \textit{insertion map}, ${\gamma:\real^{n_1} \to \real^{n_2}}$, is a  function that preserves the equilibrium point, namely $\gamma(0)=0$.
\mbox{ } \hfill \(\blacksquare\)
\end{mydef}

The insertion map specifies initial conditions for $x_2$ as a function of $x_1$; in other words, it specifies the surface $Z_0$ in the Introduction, just before \eqref{eq:PeriodicBehavior}.

\begin{enumerate}
\setlength{\itemsep}{.15cm}
\setcounter{enumi}{4}
\renewcommand\labelenumi{{\bf A-\theenumi}}
\item \label{itm:A5}  There is a constant $0 \le c < 1$,  such that, for each initial condition $\xi=(\xi_1, \gamma(\xi_1)) \in B$, there exists a piecewise continuous input $u_{\xi_1}:[0, T_p] \to \real^m$ and a corresponding solution of the ODE, $\varphi_{\xi_1}: [0, T_p]  \to \real^n$  satisfying $\varphi_{\xi_1}(T_p) \in B$, and
\begin{equation}
\label{eq:LyapunovConstraintReduced}
V((\varphi_{1\xi_1}(T_p), \gamma(\varphi_{1\xi_1}(T_p)) ) \le c V(\xi_1, \gamma(\xi_1)),
\end{equation}
where  the solution of the $(n_1 + n_2)$-dimensional model \eqref{eq:TwoBlockSystem} has been decomposed as $$\varphi_{\xi_1}(t)=:(\varphi_{1\xi_1}(t), \varphi_{2\xi_1}(t)).$$
\mbox{ } \hfill \(\blacksquare\)
\end{enumerate}

\begin{prop}
\label{prop:ReducedSystemStability} Assume the open-loop system \eqref{eq:TwoBlockSystem} satisfies Assumptions A-\ref{itm:A1} to A-\ref{itm:A3} and A-\ref{itm:A5}, and define $B_1:=\{ \xi_1 \in \real^{n_1}~|~ (\xi_1,\gamma(\xi_1))\in B \}$. Assume in addition there exists an open set $B^e_{1} \supset B_1$ and a function $$\nu:[0, \infty) \times B^e_{1} \to \real^{n_2}$$ that is piecewise continuous in $t$, $T_p$-periodic, locally Lipschitz continuous in $x_1$, and, such that, for $0\le t < T_p$ and $\xi_1 \in B_1$,
\begin{equation}
\label{eq:LearningCondition}
\nu(t,\varphi_{1\xi_1}(t))=\varphi_{2\xi_1}(t).
\end{equation}
Then the origin of the reduced-order system
\begin{equation}
\label{eq:ClosedLoopReducedFeedback}
\dot{x}_1=f_{\rm red}^{cl}(t,x_1):=f_1(t,x_1,\nu(t,x_1)),
\end{equation}
is locally uniformly exponentially stable, and the trajectories in A-\ref{itm:A5} are solutions of \eqref{eq:ClosedLoopReducedFeedback}.
\mbox{ } \hfill \(\blacksquare\)
\end{prop}

\begin{table}
	\centering
		\begin{tabular}{c|c||c}
			\multicolumn{2}{c||}{\bf Features } &   { \bf Labels}  \\
			\hline
			$t_j$& $x_{1}^{j, i}=\varphi_{1\xi_{1}^i}(t_j) $ & $\nu^{j, i}=\varphi_{2\xi_{1}^i}(t_j)$ \\
			\hline
            & & \\
			$t_0=0$ & $x_{1}^{0,1} $  & $\nu^{0,1}$ \\
            & & \\
			$t_0=0$ & $x_{1}^{0,2} $  & $\nu^{0,2}$ \\
			$\vdots$ & $\vdots$& $\vdots$\\
			$t_0=0$ & $x_{1}^{0,M} $  & $\nu^{0,M}$ \\
            & & \\
			$t_1$ & $x_{1}^{1,1} $  & $\nu^{1,1}$ \\
			$\vdots$ & $\vdots$& $\vdots$\\
			$t_1$ & $x_{1}^{1,M} $  & $\nu^{1,M}$ \\
			$\vdots$ & $\vdots$& $\vdots$  \\
			$t_N=T_p$ & $x_{1}^{N,1} $  & $\nu^{N,1}$ \\
			$\vdots$ & $\vdots$& $\vdots$  \\
			$t_N=T_p$ & $x_{1}^{N,M} $ & $\nu^{N,M}$ \\
		\end{tabular}
	\caption[]{Conceptual arrangement of the data in A-\ref{itm:A5} from which a controller satisfying \eqref{eq:LearningCondition} may be determined by ``regression''. Since only the $x_1$-component of the state is sampled, far fewer optimizations are required. The number of time samples remains the same.}
	\label{tab:How2DoRegressionReduced}
\end{table}

\begin{remark}
\label{rmk:gamma}
\
\begin{enumerate}

    \item[(i)] Assumption A-\ref{itm:A5} and Prop.~\ref{prop:ReducedSystemStability} represent our first result to  mitigate the curse of dimensionality. Assumption A-\ref{itm:A5}  leads to a greatly reduced training set for building a realization than A-\ref{itm:A4} because, in many practical examples, $n_1 \ll (n_1 + n_2)$. Proposition~\ref{prop:ReducedSystemStability} says that this reduced training set can encode a stabilization goal that is a feasible action of the full-order model. The next section embeds the target model \eqref{eq:ClosedLoopReducedFeedback} in the full-order model, completing our basic plan for mitigating the curse of dimensionality.

    \item[(ii)] The numerical burden of developing the ``training sets'' for $\nu(t,x_1)$ is exponential in the dimension of $x_1$, at least if a uniform grid is used to sample $B_1$.

\item[(iii)] Table \ref{tab:How2DoRegressionReduced} shows how to extract the function $\nu(t,x_1)$ from the optimization data.

\item[(iv)] The local Lipschitz continuity of $\nu(t,x_1)$ imposes stronger conditions on the solutions given in A-\ref{itm:A5} than those encountered in A-\ref{itm:A4}. This is because the mapping defined by, for each $t \in [0, T_p]$,
\begin{equation}
\label{eq:PsiReduced}
\Psi_{1t}:B_1 \to \real^{n_1},~ \Psi_{1t}(\xi_1):=\varphi_{1\xi_1}(t)
\end{equation}
being injective is stronger than the mapping
\begin{equation}
\label{eq:PsiFull}
\Psi_{t}:B_1 \to \real^{n},~ \Psi_{t}(\xi_1):=\left[ \begin{array}{c} \varphi_{1\xi_1}(t) \\ \varphi_{2\xi_1}(t) \end{array} \right]
\end{equation}
      being injective.  If $\Psi_{t}$ is continuously differentiable and full rank, then there does exist a new choice of $x_1$-coordinates for which the corresponding mapping $\Psi_{1t}$ is full rank and hence is locally injective. This will be illustrated on the cart-pendulum model in Sect.~\ref{sec:PendulumOnCart}.

\item[(v)] Under the assumptions of Prop.~\ref{prop:ReducedSystemStability}, for each $t \in [0, T_p)$, $\Psi_{t}:B_1 \to \real^{n}$ is a homeomorphism onto its image. It follows that $\Psi_{e}: [0, T_p) \times B_1 \to \real^{n+1}$, by
    $$\Psi_{e}(t,\xi_1):= \left[ \begin{array}{c} t \\ \Psi_{t}(\xi_1) \end{array} \right],$$
    is also a homeomorphism onto its image. After augmenting the state with time in the usual manner, the low-dimensional model \eqref{eq:PeriodicRealization} discussed in the Introduction can be seen as evolving on the surface
\begin{equation}
\label{eq:ZSurface}
Z_{T_p}:=\Psi_{e}( [0, T_p) \times B_1 ),
\end{equation}
    with the dynamics and output given by
    \begin{equation}
\label{eq:PreviewOfTargetModel}
\begin{aligned}
\dot{\tau}&=1\\
\dot{x}_1&=f_1(\tau,x_1,\nu(\tau,x_1))  \\
 x&=\left[\begin{array}{c} x_1 \\ \nu(\tau,x_1) \end{array}  \right].
\end{aligned}
\end{equation}
At this point, the direct relation with trajectories of the original model is only true for $0 \le t <T_p$.


 \item[(vi)] In Sect.~\ref{sec:OrbitLibrary}, we will provide a concrete way to select the insertion map. For now, we propose an insertion map inspired by backstepping
 \begin{equation}
 	\label{eq:GammaExampleFromV}
 	\gamma(x_1) := K x_1,
 \end{equation}
that the equilibrium of the reduced-order model,
 \begin{equation}
\label{eq:linearizedSystemf1}
\dot{x}_1 = f_1(t,x_1,\gamma(x_1))
\end{equation}
is stable. Other relations to backstepping are noted in Appendix~\ref{app:Relations}.


 \item[(vii)] Suppose the system \eqref{eq:TwoBlockSystem} is time invariant, so that one is stabilizing a trivial periodic orbit (i.e.,  an equilibrium). Then $T_p>0$ is a free parameter available to the designer. How to choose it? If the insertion map actually stabilizes the equilibrium of the reduced-order model \eqref{eq:linearizedSystemf1},
     then in principle, $T_p$ can be taken to be arbitrarily small, subject to choosing $c>0$ and the positive function in \eqref{eq:LyapunovConstraintReduced} properly. Otherwise, if the system is locally asymptotically controllable to the origin \cite{clarke1997asymptotic}, a larger $T_p$ makes it easier to meet the Lyapunov contraction condition.

    \end{enumerate}
\end{remark}

\subsection{Embedding the Target Dynamics in the Original System}
\label{sec:InnerOuterLoop}

Consider the system \eqref{eq:TwoBlockSystem} with the assumptions and notation of Prop.~\ref{prop:ReducedSystemStability}. Assume there exists a feedback $u(t,x_1,x_2)$ such that in the coordinates
\begin{equation}
\label{eq:NewCoordinates}
y:=x_2-\nu(t,x_1),
\end{equation}
the closed-loop system has the form
\begin{equation}
\label{eq:TwoBlockSystemNewCoordinates}
\begin{aligned}
\dot{x}_1&= f_1(t,x_1,\nu(t,x_1)+y)\\
\dot{y}&= Ay,
\end{aligned}
\end{equation}
with $A$ Hurwitz. Then the surface $y\equiv0$ is invariant and the restriction dynamics is given by \eqref{eq:ClosedLoopReducedFeedback}. While $\nu(t,x_1)$ is $T_p$-periodic, its limits from the left and the right are not necessarily equal at $T_p$, and $y$ in \eqref{eq:NewCoordinates} inherits this property. Hence, without further assumptions, it cannot be a solution of the ODE \eqref{eq:TwoBlockSystemNewCoordinates}, in the usual sense \cite[pp.~657]{KHA02}, namely
$$y(t) = y(t_0) + \int_{t_0}^{t} A y(\tau) d \tau.$$

In the following, we impose continuity in $t$ on $\nu(t,x_1)$, and then after the theorem, analyze what this means in terms of the trajectories coming out of the optimizer.

\begin{theorem}
\label{them:OverallSystemStability} Assume the open-loop system \eqref{eq:TwoBlockSystem} satisfies Assumptions A-\ref{itm:A1} to A-\ref{itm:A3}, and  A-\ref{itm:A5}, and define $B_1:=\{ \xi_1 \in \real^{n_1}~|~ (\xi_1,\gamma(\xi_1))\in B \}$. Assume in addition there exists an open set $B^e_{1} \supset B_1$ and a feedback
\begin{equation}
\label{eq:TheGoldenNu}
\nu:[0, \infty) \times B^e_{1} \to \real^{n_2}
\end{equation}
that is continuous in $t$, $T_p$-periodic, locally Lipschitz continuous in $x_1$, and, such that, for $0\le t < T_p$ and $\xi_1 \in B_1$,
\begin{equation}
\label{eq:LearningConditionC}
\nu(t,\varphi_{1\xi_1}(t))=\varphi_{2\xi_1}(t).
\end{equation}
Then any feedback $u(t, x_1,x_2)$, piecewise  continuous in $t$ and locally Lipschitz continuous in $x_1$ and $x_2$ that transforms the system to \eqref{eq:TwoBlockSystemNewCoordinates}, with $A$ Hurwtiz, renders the origin of \eqref{eq:TwoBlockSystemNewCoordinates} locally uniformly exponentially stable. Moreover, the surface defined by
\begin{equation}
\label{eq:ZeroDynamicsManifold}
Z_e:=\{(t,x_1,x_2)~|~ (t~{\rm mod}~ T_p, x_1,x_2)\in Z_{T_p} \},
\end{equation}
is invariant with restriction dynamics given by \eqref{eq:PreviewOfTargetModel}.
\mbox{ } \hfill \(\blacksquare\)
\end{theorem}

\begin{remark}
\
\begin{enumerate}
\item[(i)] In fact, \eqref{eq:ZeroDynamicsManifold} is the Isidori-Byrnes  \cite{ISI95} Zero Dynamics Manifold and the Zero Dynamics is given by
\begin{align*}
\dot{\tau}&=1\\
\dot{x}_1&=f_1(\tau,x_1,\nu(\tau,x_1)).
\end{align*}
The output that is being ``zeroed'' is
$$ y = x_2 - \nu(t, x_1) $$
as long as the domain is properly specified.
\item[(ii)] If $\nu$ in \eqref{eq:TheGoldenNu} is continuous in $t$, then for all $x_1\in B^e_{1}$, $\nu(T_p,x_1)=\nu(0,x_1)$. How does this relate to the trajectories in the training set used to generate $\nu$? Because $V$ in \eqref{eq:LyapunovConstraintReduced} is positive definite, and $V$ decreases after $T_p$ seconds, there exists an open ball $B_2$ contained in $B_1$ such that $\xi_1\in B_2 \Rightarrow \varphi_{1\xi_1}(T_p) \in B_1$. Because $\hat{\xi}_1:=\varphi_{1\xi_1}(T_p) \in B_1$,
\begin{equation}
\varphi_{2\hat{\xi}_1}(0)=\gamma(\hat{\xi}_1):=\gamma(\varphi_{1\xi_1}(T_p)).
 \end{equation}
From \eqref{eq:LearningCondition}, because $\hat{\xi}_1 \in B_1$,
\begin{equation}
\label{eq:LearningCondition22}
\nu(0,\varphi_{1 \hat{\xi}_1}(0))=\varphi_{2\hat{\xi}_1}(0),
\end{equation}
and because $\varphi_{1 \hat{\xi}_1}(0)=\hat{\xi}_1$, we have
\begin{equation}
\label{eq:LearningCondition23}
\nu(0,\varphi_{1 \hat{\xi}_1}(0))=\nu(0,\hat{\xi}_1).
\end{equation}
From the continuity of $\nu$ and using the definition of $\hat{\xi}_1$,
\begin{equation}
\label{eq:LearningCondition24}
\nu(0,\varphi_{1\xi_1}(T_p))=\nu(T_p,\varphi_{1\xi_1}(T_p)).
\end{equation}
From \eqref{eq:LearningCondition} again,
\begin{equation}
\label{eq:LearningCondition25}
\nu(T_p,\varphi_{1\xi_1}(T_p))=\varphi_{2 \xi_1}(T_p).
\end{equation}
Putting these together, the corresponding condition on the trajectories used in the training set for $\nu$ is given in A-\ref{itm:A6}.
\end{enumerate}
\end{remark}

\begin{enumerate}
\setlength{\itemsep}{.15cm}
\setcounter{enumi}{5}
\renewcommand\labelenumi{{\bf A-\theenumi}}
\item \label{itm:A6} The solutions in A-\ref{itm:A5} also satisfy
\begin{equation}
\label{eq:boundaryCondition}
\gamma(\varphi_{1\xi_1}(T_p))=\varphi_{2 \xi_1}(T_p).
\end{equation}
\mbox{ } \hfill \(\blacksquare\)
\end{enumerate}

Section \ref{sec:PendulumOnCart} will illustrate these ideas on a simple low-dimensional example to make it easy for the interested reader to reproduce the results. The true benefits of the approach will not be clear until Sect.~\ref{sec:biped}, where it will be applied to a high-dimensional hybrid model of a bipedal robot, and subsequently implemented in hardware on the robot of Figure~\ref{fig:MARLOv2}.

\subsection{Extended class of models}
\label{sec:ExtendedModels}

We discuss the case with the input appearing in both blocks. To keep the presentation brief and simple, it is supposed that the model has the form
\begin{equation}
\label{eq:TwoBlockSystemTV}
\begin{aligned}
\dot{x}_1&= f_1(t,x_1,x_2,u)\\
\dot{x}_2&= f_2(t,x_1,x_2,u),
\end{aligned}
\end{equation}
with
$$x_2=\left[ \begin{array}{c} x_{2a} \\ x_{2b}  \end{array} \right]~\text{and}~f_2= \left[ \begin{array}{c} x_{2a} \\ u \end{array} \right].$$

\begin{corollary}
\label{cor:OverallSystemStabilityExtraU} Assume the open-loop system \eqref{eq:TwoBlockSystemTV} satisfies Assumptions A-\ref{itm:A1} to A-\ref{itm:A3}, A-\ref{itm:A5}, and A-\ref{itm:A6}, and define $B_1:=\{ \xi_1 \in \real^{n_1}~|~ (\xi_1,\gamma(\xi_1))\in B \}$. Assume in addition there exists an open set $B^e_{1} \supset B_1$, a function
\begin{equation}
\label{eq:TheGoldenNu2}
\nu:[0, \infty) \times B^e_{1} \to \real^{n_2}
\end{equation}
satisfying the conditions of Theorem~\ref{them:OverallSystemStability}, and a second function
\begin{equation}
\label{eq:TheGoldenMu}
\mu:[0, \infty) \times B^e_{1} \to \real^{m}
\end{equation}
that is piecewise continuous in $t$, $T_p$-periodic, locally Lipschitz continuous in $x_1$, and such that, for $0\le t < T_p$ and $\xi_1 \in B_1$,
\begin{equation}
\label{eq:LearningConditionCorollary}
\mu(t,\varphi_{1\xi_1}(t))=u_{\xi1}(t).
\end{equation}
Then for all $\rfrac{n_2}{2} \times \rfrac{n_2}{2} $ positive definite matrices $K_p$ and $K_d$, the origin of
\begin{equation}
\label{eq:ModelWithFeedback}
\begin{aligned}
\dot{x}_1&= f_1(t,x_1,x_2,u)\\
\dot{x}_2&= f_2(t,x_1,x_2,u) \\
u&= \mu(t,x_1) - [K_p~~K_d]\big(x_2 - \nu(t, x_1) \big)
\end{aligned}
\end{equation}
is locally exponentially stable, uniformly in $t_0$. Moreover, the surface defined by \eqref{eq:ZeroDynamicsManifold}
is invariant with restriction dynamics given by
    \begin{equation}
\label{eq:ReducedDynamicsExtraUterm}
\begin{aligned}
\dot{\tau}&=1\\
\dot{x}_1&=f_1(t,x_1,\nu(t,x_1), \mu(t,x_1))  \\
 x&=\left[\begin{array}{c} x_1 \\ \nu(\tau,x_1) \end{array}  \right].
\end{aligned}
\end{equation}
\mbox{ } \hfill \(\blacksquare\)
\end{corollary}


\begin{remark}
\label{rmk:delta}
\
\begin{enumerate}
\item[(i)] There is no extra boundary condition, such as A-\ref{itm:A6}, associated with \eqref{eq:LearningConditionC} because the term $\mu$ arises from the inputs instead of the states of the ODE, as in the case of $\nu$. In particular, $\mu$ can be piecewise continuous in $t$. The learning of $\mu$ is done the same as for $\nu$ in Table~\ref{tab:How2DoRegressionReduced}.

\item[(ii)] \label{re:pre-feedback} Most systems will require a pre-feedback to arrive at the form \eqref{eq:TwoBlockSystemTV}; this must be taken into account when implementing the feedback indicated in \eqref{eq:ModelWithFeedback}.

    \end{enumerate}
\end{remark}

\subsection{Orbit Library and Design of the Insertion Map}
\label{sec:OrbitLibrary}

The objective of this section is to provide a systematic means for designing the insertion map in a way that takes into account the ``physics'' of a model.
\begin{mydef}
\label{def:OrbitLibrary}
An orbit library ${\cal L}$ is a set of periodic trajectories of the model \eqref{eq:TwoBlockSystemTV} that are parameterized by the $x_1$-states. We denote the library consisting of the periodic solutions by
\begin{equation}
\label{eq:OrbitLibraryB1}
{\cal L}:=\{ \varphi_{\xi_1}:[0, T_p]\to \real^n~|~\xi_1 \in B_1 \},
\end{equation}
with $B_1$ as A-\ref{itm:A5}.
	\mbox{ } \hfill \(\blacksquare\)
\end{mydef}

\begin{mydef}
\label{def:libraryGamma}
An insertion map associated to an orbit library \eqref{eq:OrbitLibraryB1} is a function $\gamma_{\cal L}:B_1 \to \real^{n_2}$ such that
\begin{equation}
\label{eq:GammaFromLibraryB1}
\gamma_{\cal L}(\xi_1):=\varphi_{2\xi_1}(0).
\end{equation}
\mbox{ } \hfill \(\blacksquare\)
\end{mydef}

One should think of the above insertion map as taking the states of the $x_1$-coordinates, associating them to periodic orbits of the full model, and then defining the initial condition of the $x_2$-coordinates (in the trajectory optimization) to be its value at a point on the associated periodic orbit. Hence, the overall model is being initialized in a physically meaningful manner. Moreover, the trajectories in A-\ref{itm:A6} can now be interpreted as affecting a transition \textit{from} a family of periodic solutions \textit{to} a  desired periodic solution \textit{in a way that} leads to stabilization of the desired periodic solution, via Theorem~\ref{them:OverallSystemStability} or Corollary~\ref{cor:OverallSystemStabilityExtraU}. The authors have found this to be very useful on bipedal robots.

\section{Inverted Pendulum on a Cart}
\label{sec:PendulumOnCart}

This section will illustrate the controller designs of Section~\ref{sec:MainIdeas} on the well-known inverted pendulum on a cart model. The MATLAB code for the calculations is available for download in Extension~\hyperlink{ex:code}{1}. The optimization setup and the learning method used here are nearly identical to what will be implemented on a bipedal robot in Section~\ref{sec:biped}; the only significant change involves the hybrid aspect of a biped model.

\subsection{System Model}

\begin{figure}[t!]
	\centering
	\includegraphics[scale=0.5]{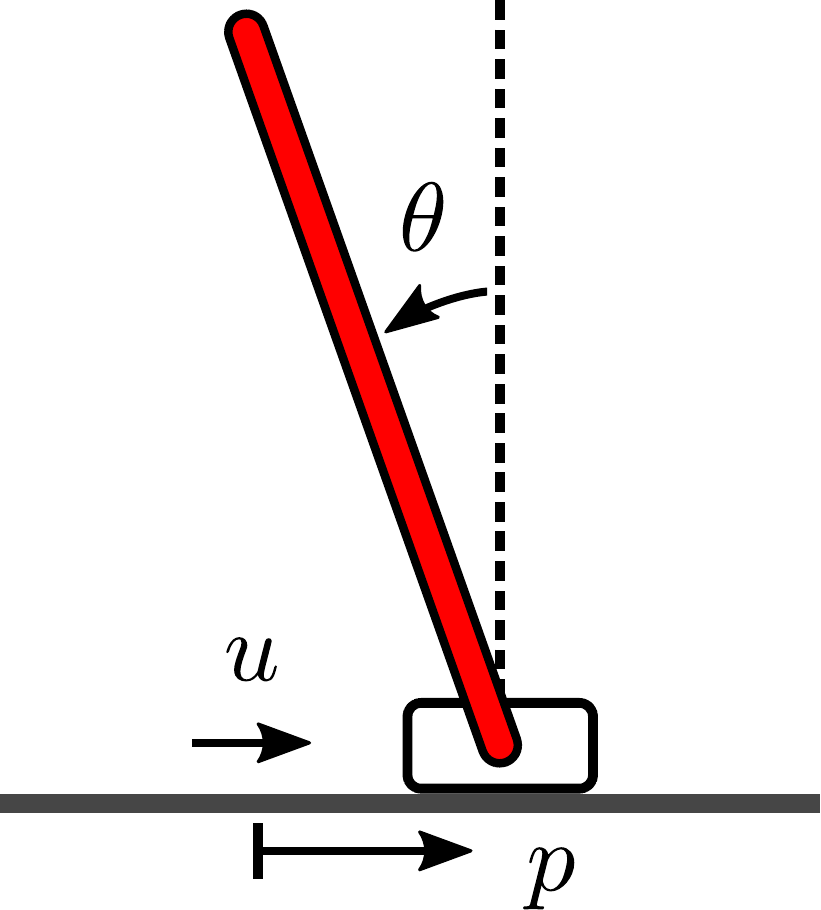}
	\caption{An inverted pendulum on a cart model is used to illustrate the controller designs of Sect.~\ref{sec:MainIdeas}. The objective is to stabilize a continuum of periodic motions, including a trivial periodic orbit corresponding to the pendulum upright and the cart at the origin. In part of the analysis, a barrier is imposed.}
	\label{fig:Pendulum}
\end{figure}

The system consists of a unit length, uniformly distributed unit-mass pendulum attached via a revolute joint to a planar unit-mass cart, shown in Figure~\ref{fig:Pendulum}. A driving force is applied on the cart, and there is no torque acting on the revolute joint of the pendulum. The motion of the cart and the pendulum are free of friction forces. The configuration variable $ q := (p, \theta) $ consists of the cart position and the pendulum angle. The system is written in state variable form as
\begin{equation}
\label{eq:Decomposedx1x2}
\begin{aligned}
	\dot{x}_1 &= \begin{bmatrix}
\dot{p}
\\
\frac{2\sin(\theta)\dot{\theta}^2 - 3g\cos(\theta)\sin(\theta) - 4u}{3\cos(\theta)^2 - 8}
\end{bmatrix}
\\
\dot{x}_2 &= \begin{bmatrix}
\dot{\theta}
\\
\frac{3\cos(\theta)\sin(\theta)\dot{\theta}^2 - 12g\sin(\theta)- 6\cos(\theta)u}{3\cos(\theta)^2 - 8}
\end{bmatrix},
\end{aligned}
\end{equation}
where $u$ is the force acting on the cart and the system state $ x $ is decomposed into $ x_1 =(p, \dot{p})$ and $ x_2 =( \theta, \dot{\theta})$.The equilibrium point of the upright pendulum is $ x^* = 0 $ and $ u^* = 0 $.  Assumptions A-\ref{itm:A1} and A-\ref{itm:A2} are then trivially satisfied.


The overall control objective will be to locally exponentially stabilize a continuum of periodic motions with a common period $ T_p = 2 $ seconds. We first illustrate the control design method on a trivial periodic orbit corresponding to the pendulum upright and the cart at the origin.

\subsection{Stabilizing the Upright Equilibrium While Respecting a Barrier}
\label{sec:StablizeEquilibrium}

The presentation follows the basic steps of the design, from learning a full-state feedback as in Prop.~\ref{prop:XDaMPCStability} to embedding a target model as in Corollary~\ref{cor:OverallSystemStabilityExtraU}.

\subsubsection{Trajectory Generation and Learning for the Full-Order Model}

The set $B$ and positive definite function of  A-\ref{itm:A3} are discussed shortly. For an initial state $ \xi \in B $, the direct collocation algorithms of \cite{Jo2014, HeCoHuAm16} are used to generate a trajectory $ \varphi_\xi(t) $ and corresponding input $ u_\xi(t) $ over an interval  $ [0, T_p] $ to meet the conditions of A-\ref{itm:A4}. To emphasize the ability to handle interesting constraints in the control design, the cart position is heavily penalized if it moves out of a ``safe region'' $ [-p_b, p_b] $, with $ p_b = 2 $.

The cost function for determining the trajectories is a standard quadratic form with an additional penalty for the safety region:
\begin{equation}
\label{eq:CostFunctionPendulum}
\begin{aligned}
J(\xi) &= \min_{u} \int_{0}^{T_h} \big( ||x||^2_Q + ||u||^2_R  + L(p, p_b) \big)dt \\
L(p, p_b) &= w p^2(e^{p - p_b} + e^{-p - p_b}).
\end{aligned}
\end{equation}
The weights $ Q $ and $ R $ are taken as identity matrices and the penalty weight is $ w = 10 $.
The optimization is subjected to the system dynamics constraints\eqref{eq:Decomposedx1x2} and the terminal constraint $ x(T_h) = 0 $,
with $ T_h = 3T_p = 6 $ seconds. One could also use a terminal cost $ N(x(T_h)) $ in place of the terminal constraint. Even though a terminal constraint may make the optimization problem infeasible for some initial conditions, we have found it to be quite practical in bipedal robots. As discussed earlier, the cost function in \eqref{eq:CostFunctionPendulum} can often used as a Lyapunov function meeting the conditions in A-\ref{itm:A4}. Here, we do not add this as a constraint to the optimization and will illustrate the satisfaction of the Lyapunov condition.

The function $ \mu(t, x) $ in Prop.~\ref{prop:XDaMPCStability} is learned for the ball of initial conditions
\begin{equation}
\label{eq:sampleFullState}
\begin{aligned}
B= &\{-1 \le p \le1, -\frac{\pi}{6} \le \theta \le \frac{\pi}{6}, \\
   &-2 \le \dot{p} \le 2, -2  \le \dot{\theta} \le 2 \},
\end{aligned}
\end{equation}
with samples $ \xi^i \in B$ selected from a uniform grid of the state. Five points are used in each dimension, for a total of 625 input sequence $ u_{\xi^i}(t) $ and solutions $ \varphi_{\xi^i}(t) $. At each time
\begin{equation}
\label{eq:sampleTime}
t_j \in \{t~|~ j\dfrac{T_p}{40}, j = 0, 1,\ldots, 40 \},
\end{equation}
the time-state pair $ (t_j, x^{j,i}) $ is a feature and the input $ u^{j, i} $ is a label. The complete list of features and labels is shown in Table.~\ref{tab:How2DoRegression} in Section~\ref{sec:MPCstyleController}. We use the MATLAB Neural Network Fitting Toolbox to approximate $\mu$. The fitting setup is shown in Table~\ref{tab:NNparameters}. The mean squared error of the validation set is around $10^{-4}$.

\begin{figure}[t!]
	\centering
	\includegraphics[width=0.9\columnwidth]{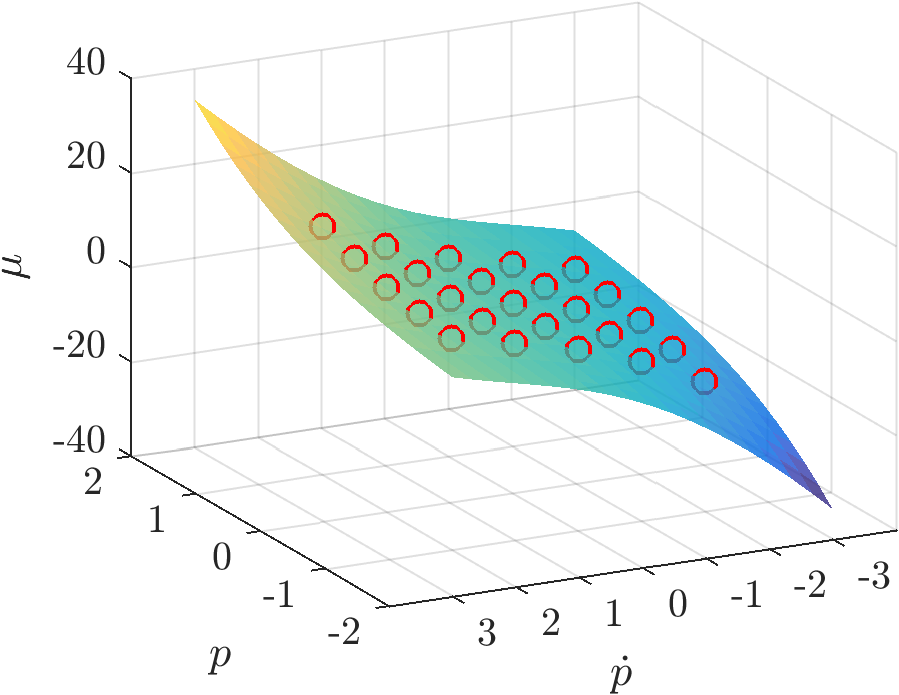}
	\caption{A slice of the function $\mu(t,p, \dot{p}, \theta,  \dot{\theta})$, with $t=0$, $\theta=0$ and $\dot{\theta}=0$. The presence of the (soft-penalty) barrier is most evident near $p=-2$ and $\dot{p}=-3$. The circles are training and validation data. Both interpolation and extrapolation can be seen in the surface.}
	\label{fig:fitting}
\end{figure}

\begin{table}
	\caption[]{MATLAB Neural Network Fitting Parameters}
	\label{tab:NNparameters}
	{
	\renewcommand{\arraystretch}{1.5}%
	\begin{center}
		\begin{tabular}{c|c}
			\hline
			Hidden neurons & 50 \\
			Training Ratio & 80\% \\
			Validation Ratio & 20\% \\
			Training Algorithm & Bayesian Regularization \\
			Max Iteration & 4000 \\
		\end{tabular}
	\end{center}
	}
\end{table}

\begin{figure}[t!]
	\centering
	\includegraphics[width=0.9\columnwidth]{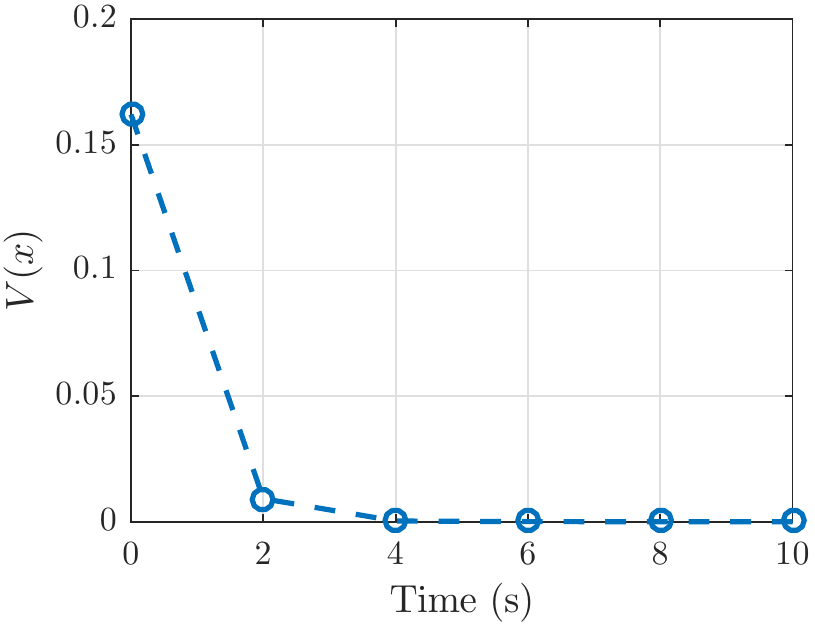}
	\caption{The plot shows the typical evolution of the optimization-cost function, confirming that it serves as a Lyapunov function. It is to be noted that $V(t)$ is only required to monotonically decrease from sample to sample, that is, from $kT_p$ to $(k+1)T_p$, with $T_p=2$.}
	\label{fig:LyapunovFunctionFullState}
\end{figure}

We show there exists a Lyapunov function
\begin{equation}
V(x) := x^\top P x
\end{equation}
as required in A-\ref{itm:A4} that is built from the cost function in \eqref{eq:CostFunctionPendulum}. The matrix
\begin{equation}
P = \left[ \begin{array}{rrrr}
0.04 & -0.11 &  0.03 & -0.03 \\
-0.11 & 0.94 & -0.12 & 0.18 \\
0.03 & -0.12 &  0.03 & -0.03 \\
-0.03 & 0.18 & -0.03 & 0.04
\end{array} \right]
\end{equation}
is from a regression of $ J(x) $. The matrix is positive definite. We next find the constant $ c $ in A-\ref{itm:A4} using the data set to train $ \mu(t, x) $ that $ c $ satisfies
\begin{equation}
c \ge \max_{i} \frac{V(\varphi_{\xi^i}(T_p))}{V(\varphi_{\xi^i}(0)}.
\end{equation}
The maximum ratio over 625 points of $ \xi^i $ is 0.22, then we set $ c = 0.25 $. Notice when $ \xi^i $ is the equilibrium point, $ V(\varphi_{\xi^i}(T_p) =  V(\varphi_{\xi^i}(0)) = 0$, which has to be ignored when finding $ c $. Figure~\ref{fig:LyapunovFunctionFullState} shows that $V(\varphi_\xi(kT_p))$ in the simulation is exponentially decreased.

\begin{figure*}
	\centering
	\begin{subfigure}{\columnwidth}
		\centering
		\includegraphics[width = \columnwidth]{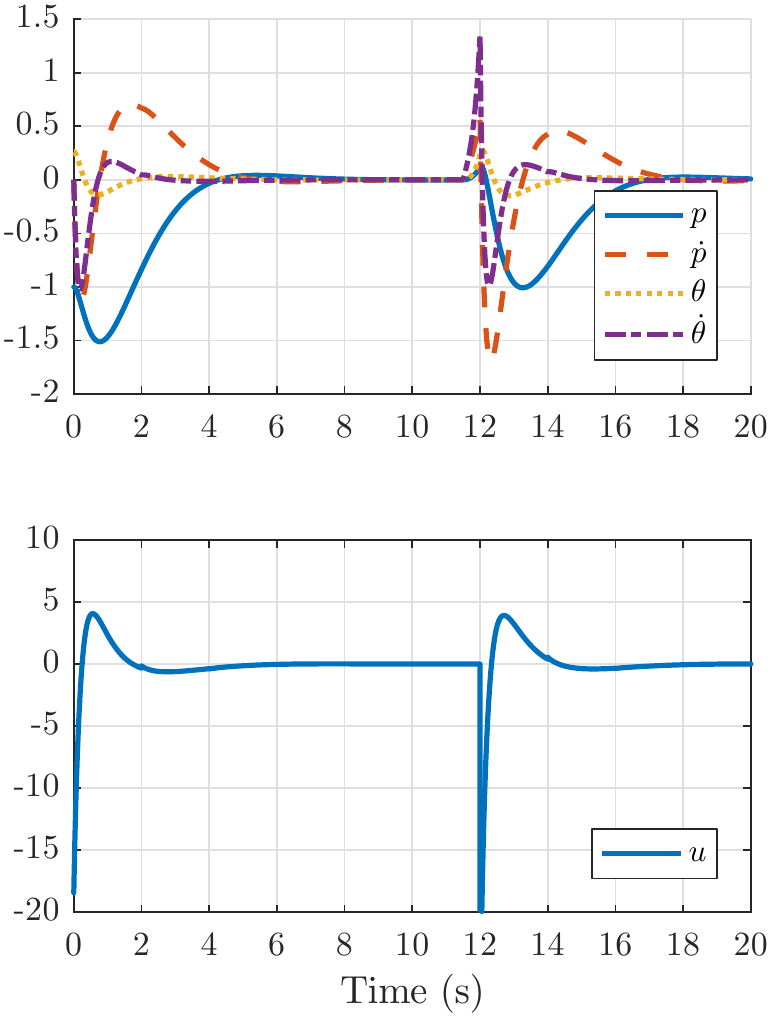}
		\caption{Continuous-hold controller $ u^{ch} $}
	\end{subfigure}
	~
	\begin{subfigure}{\columnwidth}
		\centering
		\includegraphics[width = \columnwidth]{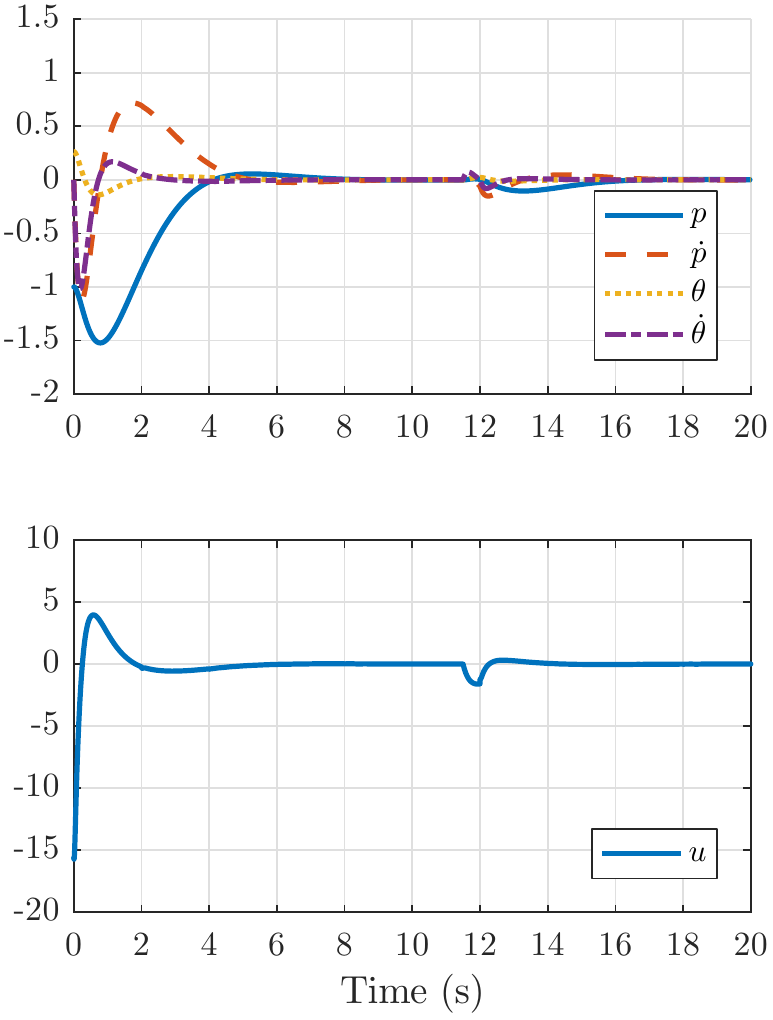}
		\caption{Learned feedback $ \mu $}
	\end{subfigure}
\caption{Comparison of the continuous-hold controller in (a), versus the learned controller in (b). In both cases, the initial condition is $(p, \dot{p}, \theta,  \dot{\theta})=(-1, 0, \rfrac{\pi}{12}, 0)$. A disturbance is applied for $t \in [11.5, 12]$ seconds. The classical MPC plot is in an Appendix-\ref{app:Standard_MPC} }
\label{fig:CH_CUcontrollerDiscussion}
\end{figure*}

\subsubsection{Comparison of Continuous Hold vs Learned Feedback}

Figure \ref{fig:CH_CUcontrollerDiscussion} compares $u^{ch}$ in \eqref{eq:pwcontfeedback} and $\mu$ in \eqref{eq:LearningConditionA}. The continuous-hold feedback $u^{ch}$ updates the state every $T_p = 2$ seconds. This type of MPC-style controller is guaranteed to perform poorly in the face of disturbances occurring within the sample period. Figure \ref{fig:CH_CUcontrollerDiscussion}-(a) shows the continuous-hold controller stabilizing the system from the initial condition
\begin{equation}
\label{eq:pendulumInitialCondition}
(p, \dot{p}, \theta,  \dot{\theta})=(-1, 0, \rfrac{\pi}{12}, 0).
\end{equation}
Figure \ref{fig:CH_CUcontrollerDiscussion}-(b) shows that the learned feedback $\mu$ performs identically to $u^{ch}$ given the same initial condition and a perfect model. The difference is obvious when it comes to disturbance rejection. A constant external force $ d = 1N $ is applied to the cart for $ t \in [11.5, 12] $. The continuous-hold feedback $u^{ch}$ has to wait till $ t = 12 $ to update the state and respond to the disturbance; on the other hand, the learned feedback $\mu$ updates the state continuously and hence reacts to the disturbance immediately.

\subsubsection{Building a Reduced-Order Target Model}

The previous subsection illustrated how feedback is extracted from data via Supervised Machine Learning and will provide a benchmark for later designs. Because the cart-pendulum system has only four states, the number of trajectories required for training in the previous subsection was quite manageable, and with random sampling techniques \cite{parisini1995receding}, it could be further reduced. Eventually, however, the number of required optimizations will become untenable. Here we illustrate Prop.~\ref{prop:ReducedSystemStability} for a two-dimensional subsystem of the cart-pendulum system.

Recall that the system state decomposition was already shown in \eqref{eq:Decomposedx1x2}. The insertion map used here is inspired by backstepping as in Remark~\ref{rmk:gamma}. Linearizing the $ x_1 $-subsystem with $u=0$ and selecting $ x_2$ as a stabilizing linear feedback yields
\begin{equation}
\label{eq:BackSteppingInsertionFunction}
\gamma(x_1) =\left[ \begin{array}{cc} 0.03 & 0.1\\ 0 & 0 \end{array} \right]  \left[ \begin{array}{c}p\\ \dot{p} \end{array} \right].
\end{equation}
We do not further explore the choice of $ \gamma $ because we will primarily use the orbit library $ \gamma_{\cal L} $ of Definition \ref{def:libraryGamma} in the remainder of the paper, including the bipedal robot section.

\begin{figure}[t!]
	\centering
	\includegraphics[width=0.9\columnwidth]{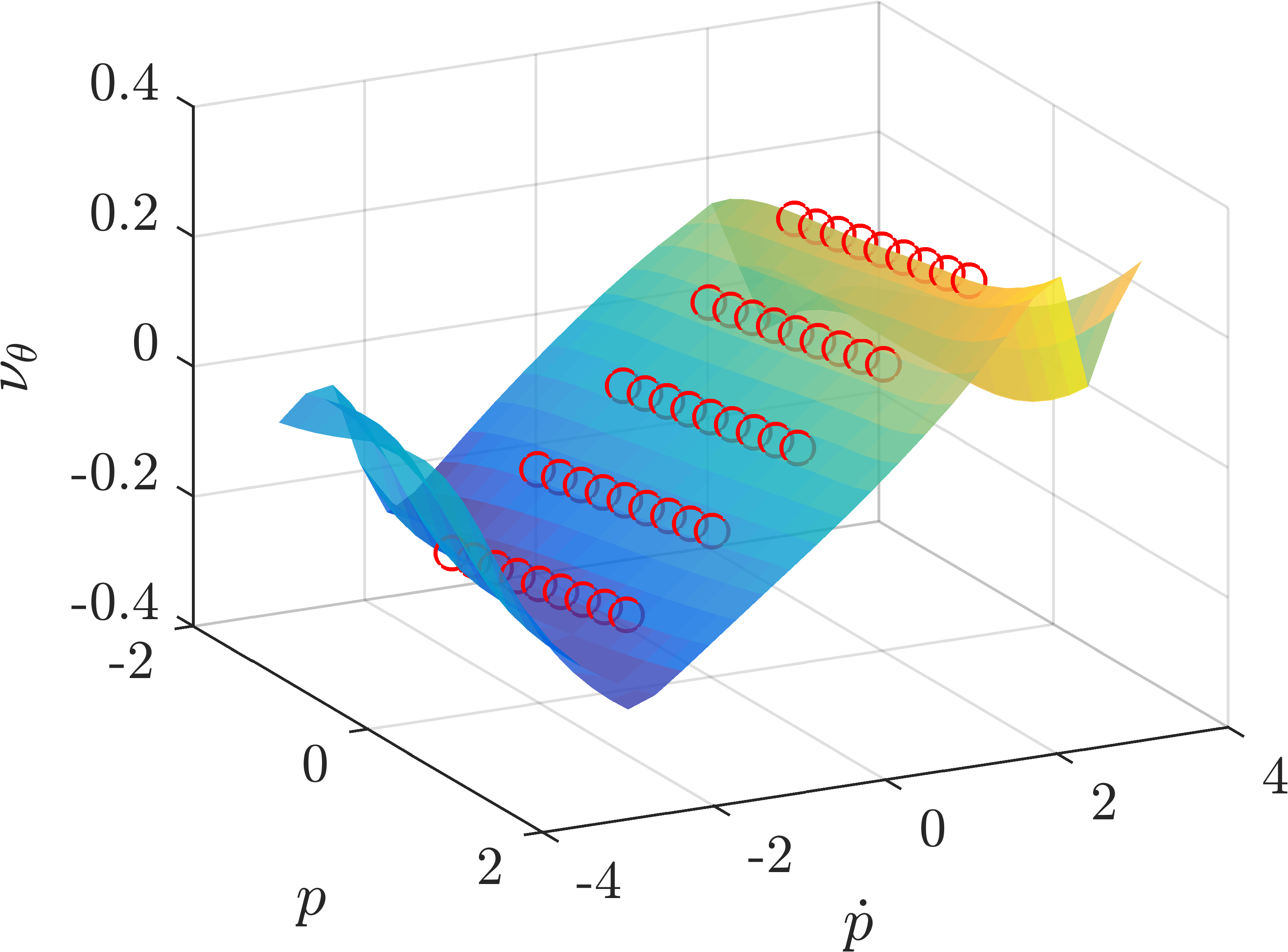}
	\caption{A slice of the function $\nu(t,p, \dot{p})$, with $t=0$. This is the $ \nu $ associated with $ x_2 $ coordinate $ \theta $. The circles are training and validation data. The interpolation is smooth while the extrapolation may be not. }
	\label{fig:fittingReduced}
\end{figure}

With this insertion map, the trajectories required by A-\ref{itm:A5} are determined via optimization with
 \begin{equation}
\label{eq:sampleReducedState}
B_1= \{-1 \le p \le 1, -2 \le \dot{p} \le 2 \}.
\end{equation}
In anticipation of using the results here in Corollary~\ref{cor:OverallSystemStabilityExtraU}, the boundary condition of A-\ref{itm:A6} is also imposed.

The set of initial conditions $ \xi_1^j,~j\in J $  now has 25 points instead of 625 points. The mapping \eqref{eq:PsiReduced} is checked to be injective by evaluating the numerical rank of the $x_1$-features in Table~\ref{tab:How2DoRegressionReduced} via SVD. Just as before, the function $\nu(t,x_1)$ is obtained from the data via the MATLAB Neural Network Fitting Toolbox, again with the parameters indicated in Table~\ref{tab:NNparameters}. The same holds for the function $\mu(t,x_1)$. An example of the fitting of $ \nu $ is shown in Figure~\ref{fig:fittingReduced}.

The evolution of the target model is shown in the next subsection when a disturbance is applied after $y=x_2-\nu(t,x_1)$ has nearly converged to zero.

\subsubsection{Embedding the Target Dynamics in the Original System}
\label{sec:InnerOuterPendulum}

The learned functions from the reduce-order optimization are now used to stabilize the full-order system based on Theorem \ref{them:OverallSystemStability} and Corollary \ref{cor:OverallSystemStabilityExtraU}. To place the system in the form \eqref{eq:TwoBlockSystemTV}, a pre-feedback is applied
\begin{equation}
	\label{eq:general_u}
	\bar{u} := \frac{3\cos(\theta)\sin(\theta)\dot{\theta}^2 - 12g\sin(\theta)- 6\cos(\theta)u}{3\cos(\theta)^2 - 8}
\end{equation}
resulting in
$$\ddot{\theta}=\bar{u}.$$
The original input $ u $ can be computed from $ \bar{u} $ because \eqref{eq:general_u} is invertible in the operational range of interest, namely $-\rfrac{\pi}{2} < \theta < \rfrac{\pi}{2} $. While the function $ \bar{\mu}(t, x) $ of Corollary \ref{cor:OverallSystemStabilityExtraU} can be recovered from $\mu(t,x_1)$ and \eqref{eq:general_u}, it is just as easy to learn it with the features $ (t_j, x_{1}^{j, i}) $ and label $ \bar{u}^{j, i} $. In the full model,
\begin{equation}
\label{eq:barUfeedback}
\bar{u} = \bar{\mu}(t, x_1) -  [K_p~~K_d]\big(x_2 - {\nu}(t, x_1) \big),
\end{equation}
with $K_p=50$ and $K_d=15$.

\begin{figure}[t!]
	\centering
	\includegraphics[width=0.9\columnwidth]{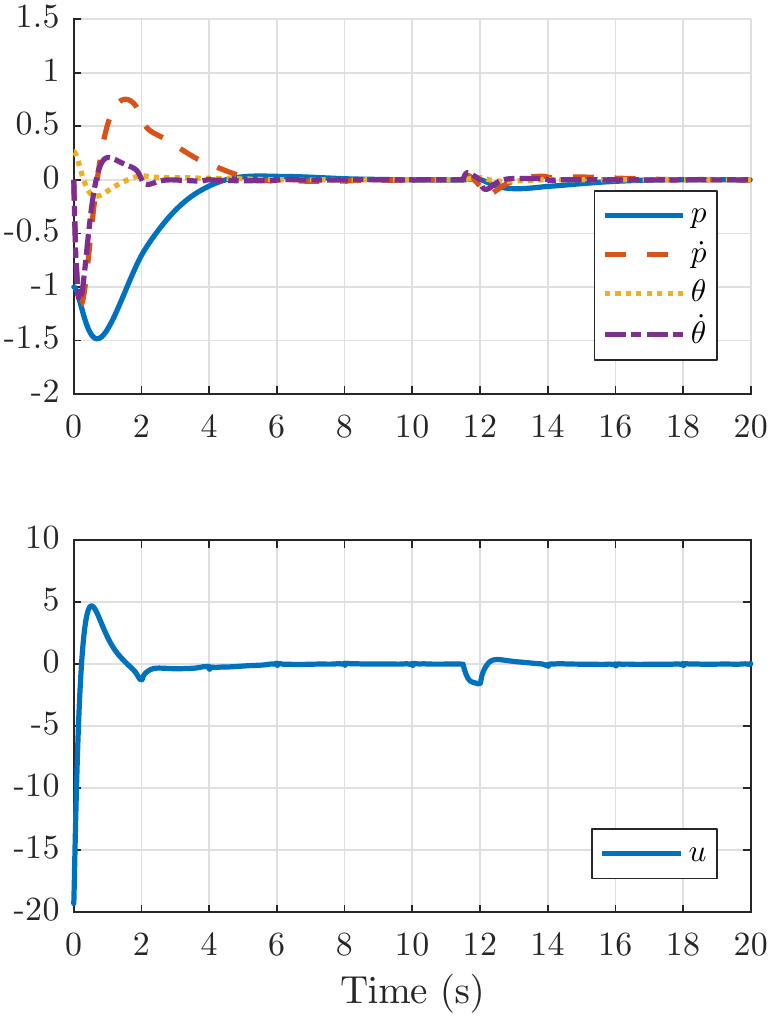}
	\caption{Response of the reduced-order model \eqref{eq:ReducedDynamicsExtraUterm}. The states of the model are $p$ and $\dot{p}$, while $\theta$ and $\dot{\theta}$ are outputs. The initial condition and disturbance are as in Figure~\ref{fig:CH_CUcontrollerDiscussion}.}
	\label{fig:CUreduced}
\end{figure}

\begin{figure}[t!]
	\centering
	\includegraphics[width=0.9\columnwidth]{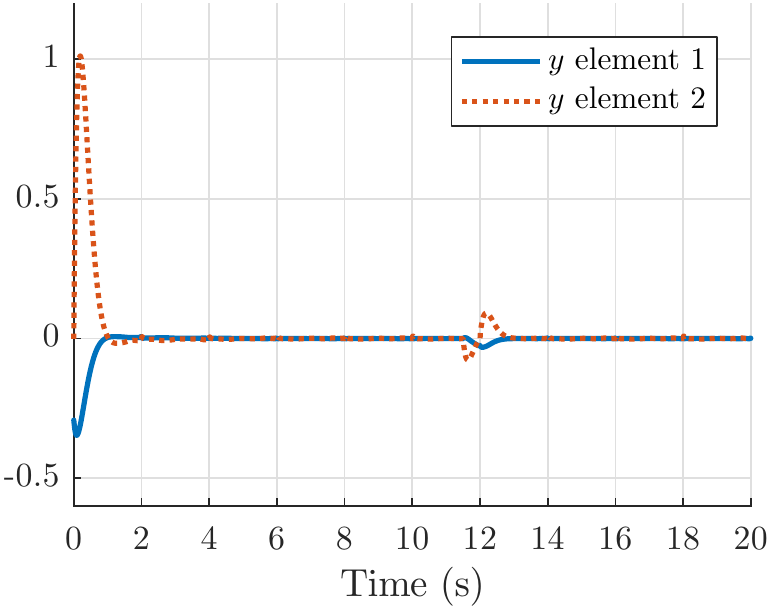}
	\caption{Showing the convergence of the output error in \eqref{eq:NewCoordinates}. Comparing this figure with Figure~\ref{fig:CUreduced} shows that the system converges to the zero dynamics surface more quickly than it converges to the periodic orbit. The disturbance initially drives the system away from the surface.} 
	\label{fig:attractiveness}
\end{figure}

Figure.~\ref{fig:CUreduced} shows the response of the closed-loop system with the same initial condition and perturbation of Figure.~\ref{fig:CH_CUcontrollerDiscussion}. The settling time and disturbance rejection performance is similar to the full state learned feedback. Figure~\ref{fig:attractiveness} illustrates the attractiveness of the surface $x_2 = \nu(t, x_1) $ by showing that the output error in \eqref{eq:NewCoordinates}
of the full-order system decays exponentially to zero. When the disturbance is applied for $11.5 \le t < 12$, the output is driven away from zero and then decays back quickly when the disturbance is removed.

\subsection{Orbit Library and Transitioning Among Periodic Orbits}

The last subsection has gone through the control design process for a trivial periodic orbit where the pendulum is upright, and the cart is at the origin. This subsection designs a controller for a set of periodic orbits, illustrate an insertion map $ \gamma_{\cal L} $ arising from an orbit library, and shows the possibility of the mapping \eqref{eq:PsiReduced} not being injective. To simplify matters, we work directly with the cart-pendulum system \textit{after} the pre-feedback \eqref{eq:general_u} has been applied.

\subsubsection{Orbit Library}
For $T_p=2$ seconds, define a set of periodic motions of the cart by
$$ {p(t)=p_0 + \frac{\dot{p}_0}{\pi} \sin(\pi t) },$$
for $(p_0, \dot{p}_0) \in B_1 $ in \eqref{eq:sampleReducedState}.
 The trajectory for $p(t)$ fixes the acceleration of the cart, which in turn gives trajectories for $\theta(t)$, $\dot{\theta}(t)$, and $u(t)$. Moreover, imposing $-\rfrac{\pi}{2} < \theta < \rfrac{\pi}{2} $ selects among the two possible solutions for the model. These considerations define an orbit library ${\cal L}$, with solutions indexed by $(p_0, \dot{p}_0)$. Denote the set of initial conditions of the orbit library as $ (\xi^{\cal L}_1, \xi^{\cal L}_2) $. Recalling Definition \ref{def:libraryGamma}, an insertion map associated to the orbit library is
 $$ \gamma_{\cal L}(\xi^{\cal L}_1)=\xi^{\cal L}_2. $$
To make it more explicit, for this example, we use linear regression to find $\gamma_{\cal L}:\real^2 \to \real^2$ as
 \begin{equation}
 \label{eq:linearRegressionOrbitGamma}
 \left[ \begin{array}{c}\theta_0\\ \dot{\theta}_0 \end{array} \right] = \left[ \begin{array}{cc} 0 & 0\\ 0 & 0.5911 \end{array} \right]  \left[ \begin{array}{c}p_0\\ \dot{p}_0 \end{array} \right].
 \end{equation}

 \begin{remark} The reader may be wondering why we bring up the orbit library as a means of computing a new insertion function, especially when the `backstepping-inspired' insertion map worked so well? The point is that for a robot, where the dimension of $x_2$ may be twenty, one has no idea how to design a `backstepping-inspired' insertion map, whereas the concept of an orbit library extends naturally as will be seen later in the paper.
 \end{remark}

\subsubsection{Loss and Recovery of Injectivity}

\begin{figure}[t!]
	\centering
	\includegraphics[width=0.9\columnwidth]{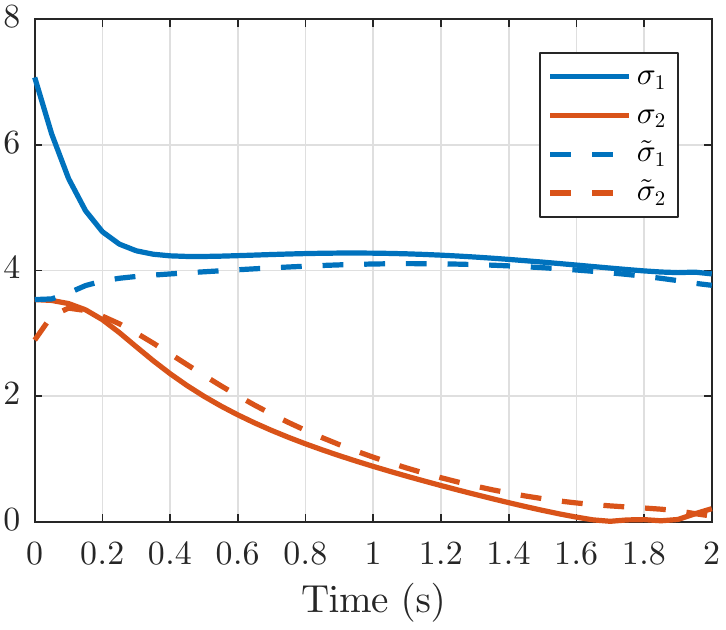}
	\caption{The singular values of the matrix formed by the sampled trajectories versus time the insertion function is given by \eqref{eq:linearRegressionOrbitGamma}. The solid lines are for $(p(t), \dot{p}(t))$, with $\Delta t=0.05$, while the dashed lines are for $(p(t)-\theta(t), \dot{p}(t)-\dot{\theta}(t))$ for the same time samples. For the choice of the insertion function arising from backstepping in \eqref{eq:BackSteppingInsertionFunction}, the minimum value of $\sigma_2$ was $ 0.58 $.}
	\label{fig:svdCheck}
\end{figure}

\begin{figure}[t!]
	\centering
	\includegraphics[width=1\columnwidth]{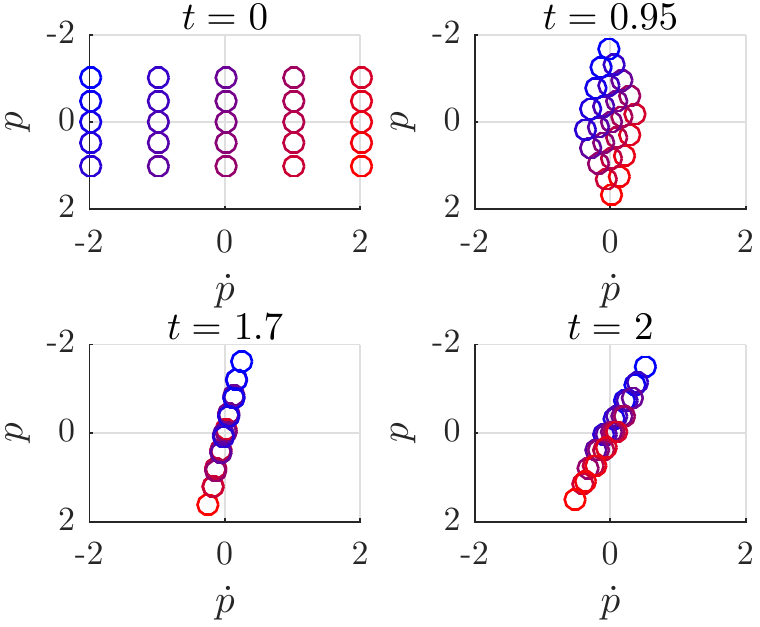}
	\caption{Another perspective on the information in Figure~\ref{fig:svdCheck}. The initial conditions are taken from a grid, as can be seen at $t=0$. At subsequent times, the grid is transformed into a parallelogram at $t=0.95$ and a line at $t=1.7$, where the mapping $\Psi_{1t}:B_1 \to \real^{n_1}$ in \eqref{eq:PsiReduced} loses rank. We have not yet observed this problem in the case of bipedal robots.}
	\label{fig:featureSinglarity}
\end{figure}

Using the new insertion function \eqref{eq:linearRegressionOrbitGamma}, trajectory generation is performed exactly as in Sect.~\ref{sec:StablizeEquilibrium}, with $x_1$ and $x_2$ as given in \eqref{eq:Decomposedx1x2}. The mapping \eqref{eq:PsiReduced} is checked not to be injective. Indeed, for $t\simeq 1.8$, the cart trajectories pass through a one-dimensional surface as shown in Figures~\ref{fig:svdCheck} and \ref{fig:featureSinglarity}, while the mapping \eqref{eq:PsiFull} remains injective. Hence, one expects the existence of a set of coordinates in which the design can proceed. It can be checked that the new coordinates\endnote{Almost any linear combination of $x_1$ and $x_2$ works because it takes rows from the bottom of \eqref{eq:PsiFull} and adds them to the top rows, making \eqref{eq:PsiReduced} full rank, and hence locally injective. There is nothing magic about our choice.}
$$ \tilde{x}_1 = (p - \theta, \dot{p} - \dot{\theta})$$
are ``full rank'' as shown in Figure~\ref{fig:svdCheck}.
It is not necessary to redo the optimization because the new coordinates correspond to a new $\tilde{B}_1$ and a new insertion map, $\tilde{\gamma}_{\cal L}$, and these are not explicitly required in the computation of $\nu$ and $\mu$. The important thing is the feature set for the Supervised Machine Learning is now indexed by $(t_j, \tilde{x}_1^{j,i})$ rather than $(t_j, x_1^{j,i})$.

\subsubsection{Transitioning Between Periodic Orbits}

Next, we use the library insertion map $ \gamma_{\cal L} $ and the new state $ \tilde{x}_1 $ to design a controller for transitioning between periodic orbits; this is analogous to transitioning between walking gaits of various speeds or direction for a bipedal robot. The cost function used in A-\ref{itm:A5} is modified to include the target orbit $A:=(p_0, \dot{p}_0) \in B_1$, per
\eqref{eq:CostFunctionPendulumOrbit} to
\begin{equation}
\label{eq:CostFunctionPendulumOrbit}
J(\xi_1, A)= \min_{u} \int_{0}^{T_h} \big( ||x-\varphi^A||^2_Q + ||u - u^A||^2_R \big)dt
\end{equation}
subject to $ x(T_h) = \varphi^A(T_p)$ and $x(0)=(\xi_1, \gamma_{\cal L}(\xi_1))$. The boundary condition A-\ref{itm:A6} is also applied as $ \gamma_{\cal L}(x_1(T_p))=x_2(T_p). $ Here, the target trajectory and its corresponding input are denoted as $ \varphi^A(t) $ and $ u^A(t) $. To simplify the problem, the cost function excludes the barrier penalty $ L(p, p_b) $ in \eqref{eq:CostFunctionPendulum}. The set $B_1$ is still given by \eqref{eq:sampleReducedState}.

\begin{figure}[t!]
	\centering
	\includegraphics[width=0.9\columnwidth]{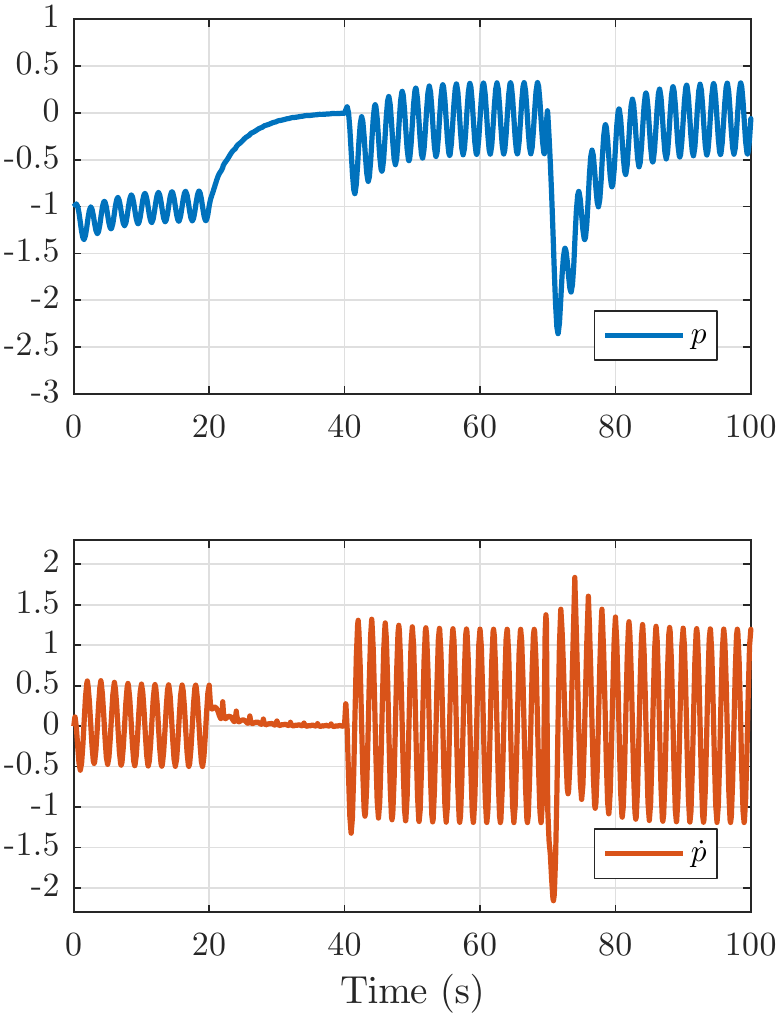}
	\caption{Plots of $p(t)$ (top) and $\dot{p}(t)$ (bottom) as the closed-loop system transitions from one periodic orbit to another, as given in \eqref{eq:transitionTarget}, with a disturbance applied for $t\in [69.5, 70]$. }
	\label{fig:transitionOrbits}
\end{figure}

With this choice of the insertion map, the boundary condition A-\ref{itm:A6} means that each cycle in the transition moves the cart-pendulum from one periodic orbit to the next; this is  because $(\varphi_{1\xi_1}(T_p), \gamma_{\cal L}(\varphi_{1\xi_1}(T_p))$ is an initial condition for a periodic solution of the model.
Denote the family of solutions to the optimization problem by
\begin{equation}
\label{eq:denoteOrbitLearning}
\begin{aligned}
	\tilde{x}_1^{i,j,k} &:= \varphi^{A_k}_{1\xi_{1}^i}(t_j) - \varphi^{A_k}_{2\xi_{1}^i}(t_j) \\
	\nu^{i,j,k} &:=\varphi^{A_k}_{2\xi_{1}^i}(t_j) \\
	\bar{\mu}^{i,j,k} &:=\bar{u}^{A_k}_{\xi^i}(t_j).
\end{aligned}
\end{equation}
The feature set for the Supervised Machine Learning is taken as $(t_j, \tilde{x}_1^{i,j,k}, A_k)$ and the labels are $(\nu^{i,j,k}, \bar{\mu}^{i,j,k})$. Figure~\ref{fig:transitionOrbits} shows a orbit transition
\begin{equation}
\label{eq:transitionTarget}
(-1, 0.5) \rightarrow (0, 0) \rightarrow (0, 1.2)
\end{equation}
of the target orbit $ (p_0, \dot{p}_0) $ at $ t = 20 $ and $ t = 40 $. A constant external force $ d = 20 N$ is applied to the cart for $ t \in [69.5, 70] $.

\begin{figure}[t!]
	\centering
\includegraphics[width=0.9\columnwidth]{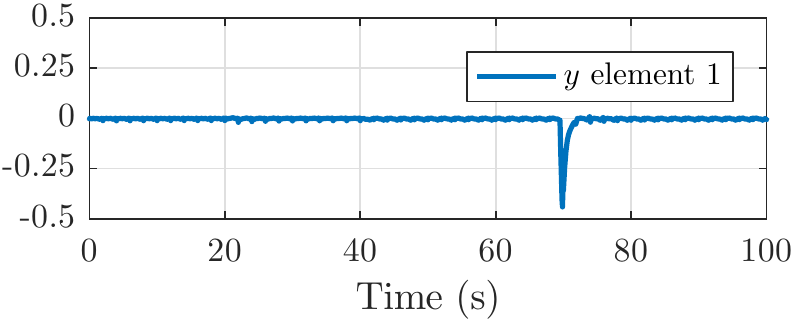}
	\caption{This shows that there is no jump in the output when the transition point takes place at a multiple of $T_p$. The jump corresponds to the disturbance in Figure~\ref{fig:transitionOrbits}. Only the first component of \eqref{eq:NewCoordinates} is shown as the other component is the derivative of this one.}
	\label{fig:transitionOutput}
\end{figure}

\begin{remark}
\
\begin{enumerate}

\item[(i)] Orbit transition from set $ B_1 $ to a target orbit $ A $ can also be reviewed as rejecting state disturbances in $ B_1 $. The distance from a state in $ B_1 $ to $ A $ is not necessarily ``small'', indicating the region of attraction for this controller could be ``large''.

\item[(ii)] There may exist two orbits $ A_m $ and $ A_n $ in $ B_1 $ for which a transition cannot be achieved over $[0, T_h] $. However, one may think of transitions in $ B_1 $ as a graph so that if there exists an orbit $ A_k $ such that
    $$ A_m  \rightarrow A_k \rightarrow A_n $$
    is possible, then the orbits are connected.

\item[(iii)] If the target orbit is modified at multiples of $T_p$, there are no jumps in $\nu$; this is because the orbit-library insertion map transitions the system from one periodic orbit to another as shown in Figure~\ref{fig:transitionOutput}.
\end{enumerate}
\end{remark}

\section{Hybrid Model and Control}
\label{sec:HybridSystemControl}

This section describes an extension of the control policy developed in Sect.~\ref{sec:MainIdeas} to systems with impulse effects \cite{grizzle2014models, WGCCM07, BASI89}, a special class of hybrid models that arises in bipedal robots. The control goals for the hybrid system corresponds to stabilizing periodic walking gaits for various speeds, and to transitioning among these gaits. The robot should also be able to reject a range of force perturbations.

\subsection{Hybrid Model}
\label{sec:HybridModel}

We consider a hybrid system with one continuous-time phase as follows
\begin{equation}\label{eq:open_loop_hybrid_model_with_one_phase}
\begin{array}{l}
\Sigma:\left\{\begin{array}{ll}
~\dot{x}~=f(x,u)\qquad\quad\ & x^{-}\notin\mathcal{S}
\\x^{+}=\Delta(x^{-}) & x^{-}\in\mathcal{S},
\end{array}\right.
\end{array}
\end{equation}
in which $x\in\mathcal{X}$ and $\mathcal{X}\subset\real^{n}$ denote the \textit{vector of state variables} and $n$-dimensional \textit{state manifold}, respectively. The continuous-time control input is represented by $u\in\mathcal{U}$, where $\mathcal{U}\subset\real^{m}$ is an open \textit{set of admissible control values}. In addition, $f:\mathcal{X} \times \mathcal{U} \rightarrow\textrm{T}\mathcal{X}$ is assumed to be continuously differentiable ( $\mathcal{C}^{1}$ ) so that a Poincar\'e map can be computed later when checking stability.  For each $u \in \mathcal{U}$,  $f(\cdot, u)$ is a vector field in $\textrm{T}\mathcal{X}$, the \textit{tangent bundle} of the state manifold $\mathcal{X}$.

The \textit{switching hypersurface} $\mathcal{S}$ is an $(n-1)$-dimensional manifold
\begin{equation}\label{switching_manifold}
\mathcal{S}:=\left\{x\in\mathcal{X}\,|\,p(x)=0\right\},
\end{equation}
on which the state solutions are allowed to undergo a sudden jump according to the \textit{re-initialization rule} $x^{+}=\Delta(x^{-})$. Here, $p:\mathcal{X}\rightarrow  \real$ is a $\mathcal{C}^{1}$-\textit{switching function} which satisfies $\frac{\partial p}{\partial x}(x)\neq0$ for all $x\in\mathcal{S}$.
 Moreover, $\Delta:\mathcal{X}\rightarrow\mathcal{X}$ denotes the $\mathcal{C}^{1}$ \textit{reset map}. $x^{-}(t):=\lim_{\tau\nearrow t}x(\tau)$ and $x^{+}(t):=\lim_{\tau\searrow t}x(\tau)$ represent the left and right limits of the state trajectory $x(t)$, respectively. As in \cite{WEGRKO03}, the solution of the hybrid system \eqref{eq:open_loop_hybrid_model_with_one_phase} is assumed to be right continuous.
 In particular, it is constructed by piecing together the flow of $\dot{x}=f(x,u)$ such that the discrete transition takes place when this flow intersects the switching hypersurface $\mathcal{S}$. The new initial condition for $\dot{x}=f(x,u)$ is then determined by the reset map $x^{+}=\Delta(x^{-})$.

\subsection{Setting up the Optimization Problem}
\label{sec:MPClikel}


We now start the translation of the Assumptions A-\ref{itm:A1} to A-\ref{itm:A6} to the hybrid setting for the purpose of designing a feedback controller to locally exponentially stabilize a periodic solution. For bipedal robots, mid-step is a good time to make adjustments to the gait: (a) impact transients have had a chance to settle out; (b) the swing foot is safely away from the ground; and (c) there is still adequate time to steer the swing leg to a favorable configuration for impact. Hence, we will use mid-step of the periodic orbit understudy for setting the beginning and end of the trajectories that we compute via optimization.

The hybrid model  \eqref{eq:open_loop_hybrid_model_with_one_phase} is assumed to satisfy the following conditions.
\begin{enumerate}
\setcounter{enumi}{0}
\renewcommand\labelenumi{{\bf H-\theenumi}}

\item \label{itm:H1} $f:\mathcal{X} \times \mathcal{U} \to T\mathcal{X}$ and the reset map $\Delta {\cal X} \to {\cal S}$ are Lipschitz continuous. This will allow the stability analysis tools of \cite{AmGaGrSr2014} to be applied later on.


\item \label{itm:H2}  There exists $T_p>0$, $x_{\rm m}^* \in \real^n$, a piecewise continuous input $u^*:[0, T_p] \to \mathcal{U}$ and a solution $\varphi^*(t)$ of \eqref{eq:open_loop_hybrid_model_with_one_phase} satisfying:
\begin{itemize}
\item[\bf a)] $\varphi^{*+}(0)=x_{\rm m}^*$;
\item[\bf b)] $\varphi^{*-}(\rfrac{T_p}{2})\in {\cal S}$, (swing foot touches the ground);
\item[\bf c)] $\forall~t\neq \rfrac{T_p}{2}, \varphi^*(t) \not \in {\cal S}$, (does so only once); and
\item[\bf d)] $\varphi^{*-}(T_p)=x_{\rm m}^*$ (periodicity).
\end{itemize}
It is noted that by the definition of ${\cal S}$, the periodic solution is transversal to ${\cal S}$, namely $\frac{d}{dt} p(\varphi^{*-}(\rfrac{T_p}{2}))<0$. And yes, the motion is being ``clocked'' with the middle of the step.

\mbox{ } \hfill \(\blacksquare\)
\end{enumerate}

The point $x_{\rm m}^*$ is the midpoint of the periodic trajectory, as measured by time. The controller we build will start from mid-stance, follow the Lagrangian model, undergo impact, and then once again evolve according to the Lagrangian model. To formulate the trajectory designs and the closed-loop system, we need to split the continuous phase of the model \eqref{eq:open_loop_hybrid_model_with_one_phase} into part-$({\rm i})$, after mid-stance, and part-$({\rm ii})$, the first half of the ensuing stance phase.

\begin{equation}\label{eq:hybridSplitPhaseModel}
\begin{array}{l}
\Sigma_{\rm i}:\left\{
\begin{aligned}
\dot{\tau}&=1, \\
\dot{x}&=f(x,u), &x^{-}\notin\mathcal{S}\\
\tau^{+}&=\tau^- \\
x^{+}&=\Delta(x^{-}), &x^{-} \in \mathcal{S}
\end{aligned}\right.
\\
\\\Sigma_{\rm ii}:\left\{
\begin{aligned}
\dot{\tau}&=1, & \tau^- < T_p\\
\dot{x}&=f(x,u),  \\
\tau^{+}&=0 & \tau^- = T_p\\
x^{+}&=x^{-}.
\end{aligned}\right.

\end{array}
\end{equation}
The guard condition on the phase-$i$ depends only on the state $x$, whereas the guard condition on the phase-$ii$ depends only on ``time'' as measured by $\tau$.

\begin{enumerate}
\setcounter{enumi}{2}
\renewcommand\labelenumi{{\bf H-\theenumi}}

 \item \label{itm:H3}  The user has selected an open ball $B\subset \real^n$ about $x_{\rm m}^*$, a positive-definite, locally Lipschitz-continuous function $V:B \to \real$, and constants $0<\alpha_1 \le \alpha_2$ such that, $V(x_{\rm m}^*)=0$ and $\forall~x\in B$,
 $$
 \begin{array}{l}
 \alpha_1( x-x_{\rm m}^*)^\top ( x-x_{\rm m}^*) \le \\
 V(x) \le \alpha_2 ( x-x_{\rm m}^*)^\top ( x-x_{\rm m}^*).
 \end{array}
 $$
  \item \label{itm:H4} There is a constant $0 \le c < 1$,  such that, for each initial condition $\xi \in B$, there exists a piecewise continuous input $u_\xi:[0, T_p] \to \real^m$  and a corresponding solution  $\varphi_\xi: [0, T_p]  \to \real^n$ of the hybrid model \eqref{eq:hybridSplitPhaseModel} satisfying
      \begin{itemize}
\item[\bf a)] $\varphi_{\xi}^+(0)=\xi$,
 \item[\bf b)] $\varphi_\xi^-(\rfrac{T_p}{2}) \in {\cal S}$,
\item[\bf c)] $\forall~t\neq \rfrac{T_p}{2}, \varphi_\xi(t) \not \in {\cal S}$,
\item[\bf d)] $\varphi_\xi^+(T_P) \in B$ and there is exponential convergence toward the periodic orbit, namely,
\begin{equation}
\label{eq:LyapunovConstraintH4}
V(\varphi_\xi^+(T_P)) \le c V(\xi),
\end{equation}
and
\item[\bf e)] $\xi = x_{\rm m}^* \Rightarrow u_\xi = u^*$.
\end{itemize}

\mbox{ } \hfill \(\blacksquare\)
\end{enumerate}

\begin{prop}
\label{prop:FullStateHybridModel} Assume the open-loop hybrid model \eqref{eq:open_loop_hybrid_model_with_one_phase} satisfies Assumptions H-\ref{itm:H1} to H-\ref{itm:H4}.
 Assume in addition there exist open sets $B^{e}_{\rm i}$ and $B^{e}_{\rm ii}$ that contain $B$, a $\delta>0$, and two feedbacks
 \begin{align*}
 \mu_{\rm i}&:[0,\rfrac{T_p}{2}+\delta] \times B^{e}_{\rm i} \to \real^m \\
 \mu_{\rm ii}&:[\rfrac{T_p}{2}-\delta, T_p] \times B^{e}_{\rm ii} \to \real^m
 \end{align*}
 that are piecewise continuous in $t$, locally Lipschitz continuous in $x$, and, such that, for $0\le t < T_p$ and $\xi \in B$,
\begin{equation}
\label{eq:LearningConditionAHybrid}
\begin{aligned}
\mu_{\rm i}(t,\varphi_{\xi}(t))&=u_\xi(t),~\ 0 \le t < \rfrac{T_p}{2} \\
\mu_{\rm ii}(t,\varphi_{\xi}(t))&=u_\xi(t),~\ \rfrac{T_p}{2} \le t < T_p.
\end{aligned}
\end{equation}
Then $\varphi^*:[0, T_p]\to \real^n$ is a locally exponentially stable periodic solution of the closed-loop system
\begin{equation}\label{eq:hybridSplitPhaseModelClosedLoopFullState}
\begin{array}{l}
\Sigma_{\rm i}^{cu}:\left\{
\begin{aligned}
\dot{\tau}&=1, \\
\dot{x}&=f(x,\mu_{\rm i}(\tau,x)), &x^{-}\notin\mathcal{S}\\
\tau^{+}&=\tau^- \\
x^{+}&=\Delta(x^{-}), &x^{-} \in \mathcal{S}
\end{aligned}\right.
\\
\\\Sigma_{\rm ii}^{cu}:\left\{
\begin{aligned}
\dot{\tau}&=1, & \tau^- < T_p\\
\dot{x}&=f(x,\mu_{\rm ii}(\tau,x)),  \\
\tau^{+}&=0 & \tau^- =  T_p  \\
x^{+}&=x^{-}.
\end{aligned}\right.
\end{array}
\end{equation}
\mbox{ } \hfill \(\blacksquare\)
\end{prop}

\subsection{Generalized Hybrid Zero Dynamics}
\label{sec:GHZD}

Following Appendix-\ref{sec:MechanicalModelDecomposition}, assume now that the continuous phase of the hybrid model has been decomposed as
\begin{equation}
\label{eq:TwoBlockSystemHybrid}
\begin{aligned}
\dot{x}_1&= f_1(x_1,x_2,u)\\
\dot{x}_2&= f_2(x_1,x_2,u),
\end{aligned}
\end{equation}
with
$$x_2=\left[ \begin{array}{c} x_{2a} \\ x_{2b}  \end{array} \right]~\text{and}~f_2= \left[ \begin{array}{c} x_{2b} \\ u \end{array} \right].$$
Let $\gamma:\real^{n_1} \to \real^{n_2}$ be a locally Lipschitz  continuous insertion function that preserves the periodic orbit, namely, writing $x_m^*=:(x_{1m}^*; x_{2m}^*)$, it follows that $\gamma(x_{1m}^*)=x_{2m}^*$.

\begin{enumerate}
\setlength{\itemsep}{.15cm}
\setcounter{enumi}{4}
\renewcommand\labelenumi{{\bf H-\theenumi}}
\item \label{itm:H5}  There is a constant $0 \le c < 1$,  such that, for each initial condition $\xi=(\xi_1, \gamma(\xi_1)) \in B$, there exists a continuous input $u_{\xi_1}:[0, T_p] \to \real^m$ and a corresponding solution $\varphi_{\xi_1}: [0, T_p]  \to \real^n$ of the hybrid model \eqref{eq:hybridSplitPhaseModel} satisfying,
          \begin{itemize}
\item[\bf a)] $\varphi_{\xi_1}^+(0)=(\xi_1;\gamma(\xi_1))$,
 \item[\bf b)] $\varphi_{\xi_1}^-(\rfrac{T_p}{2}) \in {\cal S}$,
\item[\bf c)] $\forall~t\neq \rfrac{T_p}{2}, \varphi_{\xi_1}(t) \not \in {\cal S}$,
\item[\bf d)] $\varphi_\xi^+(T_P) \in B$ and there is exponential convergence toward the periodic orbit, namely,
\begin{equation}
\label{eq:LyapunovConstraintReducedHybrid}
V((\varphi^+_{1\xi_1}(T_p), \gamma(\varphi^+_{1\xi_1}(T_p)) ) \le c V(\xi_1, \gamma(\xi_1)),
\end{equation}
and
\item[\bf e)] $(\xi_1,\gamma(\xi_1)) = x_{\rm m}^* \Rightarrow u_{\xi_1} = u^*$.
\end{itemize}
where  a solution of the $(n_1 + n_2)$-dimensional model \eqref{eq:TwoBlockSystemHybrid} has been decomposed as $\varphi_{\xi_1}(t)=:(\varphi_{1\xi_1}(t), \varphi_{2\xi_1}(t)).$

\mbox{ } \hfill \(\blacksquare\)
\end{enumerate}

The following result generalizes the hybrid zero dynamics defined in \cite{WEGRKO03,MOGR08,WGCCM07}. Even in the case of one degree of underactuation, one is able to achieve exponential stability with this method for gaits that could not be rendered stable with the previous formulation of virtual constraints. See Appendix~\ref{app:zeroDynamics} for an example. More important that this fact, however, the new formulation allows a systematic approach to robot models with more than one degree of underactuation. This is illustrated in Sect.~\ref{sec:biped}.

\begin{prop}
\label{prop:HybridZeroDynamicsStability} Assume the open-loop hybrid system \eqref{eq:open_loop_hybrid_model_with_one_phase} with $f$ given by \eqref{eq:TwoBlockSystemHybrid} satisfies Hypotheses H-\ref{itm:H1} to H-\ref{itm:H3} and H-\ref{itm:H5}, and define $B_1:=\{ \xi_1 \in \real^{n_1}~|~ (\xi_1,\gamma(\xi_1))\in B \}$.
Assume in addition there exist open sets $B^{e}_{\rm 1.i}$ and $B^{e}_{\rm 1.ii}$ that contain $B_1$, a $\delta>0$, and two feedbacks
 \begin{align*}
 \nu_{\rm i}&:[0,\rfrac{T_p}{2}+\delta] \times B^{e}_{\rm 1.i} \to \real^{n_1} \\
 \nu_{\rm ii}&:[\rfrac{T_p}{2}-\delta, T_p] \times B^{e}_{\rm 1.ii} \to \real^{n_1}
 \end{align*}
 and
 \begin{align*}
 \mu_{\rm i}&:[0,\rfrac{T_p}{2}+\delta] \times B^{e}_{\rm 1.i} \to \real^m \\
 \mu_{\rm ii}&:[\rfrac{T_p}{2}-\delta, T_p] \times B^{e}_{\rm 1.ii} \to \real^m
 \end{align*}
that are piecewise continuous in $t$, locally Lipschitz continuous in $x$, and, such that, for $0\le t < T_p$ and $\xi \in B$,
\begin{equation}
\label{eq:LearningConditionHZDnu}
\begin{aligned}
\nu_{\rm i}(t,\varphi_{1\xi_1}(t))&=\varphi_{2\xi_1}(t),~\ 0 \le t < \rfrac{T_p}{2} \\
\nu_{\rm ii}(t,\varphi_{1\xi_1}(t))&=\varphi_{2\xi_1}(t),~\ \rfrac{T_p}{2} \le t < T_p
\end{aligned}
\end{equation}
and
\begin{equation}
\label{eq:LearningConditionHZDmu}
\begin{aligned}
\mu_{\rm i}(t,\varphi_{1\xi_1}(t))&=u_{\xi_1}(t),~\ 0 \le t < \rfrac{T_p}{2} \\
\mu_{\rm ii}(t,\varphi_{1\xi_1}(t))&=u_{\xi_1}(t),~\ \rfrac{T_p}{2} \le t < T_p.
\end{aligned}
\end{equation}
Then $x_1^*:[0, T_p]\to \real^{n_1}$ is a locally exponentially stable periodic solution of the reduced-order hybrid system
%
\begin{equation}\label{eq:GHZD}
\begin{array}{l}
\Sigma_{\rm i}:\left\{
\begin{aligned}
\dot{\tau}&=1, \\
\dot{x}_1&=f_1(x_1,\nu_{\rm i}(\tau,x_1),\mu_{\rm i}(\tau,x_1)),  \\
&\text{when}\begin{bmatrix} x_1^{-} \\ \nu_{\rm i}(\tau^-,x_1^-)\end{bmatrix} \notin \mathcal{S}
\\
\tau^{+}&=\tau^-, \\
x_1^{+}&=\Delta_1(x_1^{-},\nu_{\rm i}(\tau^-,x_1^-)),
\\
&\text{when}\begin{bmatrix} x_1^{-} \\ \nu_{\rm i}(\tau^-,x_1^-)\end{bmatrix} \in \mathcal{S}
\end{aligned}\right.
\\
\\
\Sigma_{\rm ii}:\left\{
\begin{aligned}
\dot{\tau}&=1, \hspace*{.85cm} \tau^- <   T_p\\
\dot{x}_1&=f(x,\nu_{\rm ii}(\tau,x_1),\mu_{\rm ii}(\tau,x_1)),  \\
\tau^{+}&=0, \hspace*{1.0cm} \tau^- = T_p\\
x^{+}&=x^{-}.&
\end{aligned}\right.
\end{array}
\end{equation}
\mbox{ } \hfill \(\blacksquare\)
\end{prop}

\begin{remark}
\
\begin{enumerate}
\setlength{\itemsep}{.15in}
\item[(i)] In principle, $\tau^*:[0, T_p] \to \real $  needs to be defined to complete the periodic orbit, but clearly, the trivial solution, $\tau^*(t)=t,$ is the only possibility.

\item[(ii)] As in the non-hybrid case, using the trajectories in H-\ref{itm:H5}, define
\begin{equation}
\label{eq:PsiFullHybrid}
\Psi_{t}:B_1 \to \real^{n},~ \Psi_{t}(\xi_1):=\left[ \begin{array}{c} \varphi_{1\xi_1}(t) \\ \varphi_{2\xi_1}(t) \end{array} \right]
\end{equation}
and $\Psi_{e}: [0, T_p) \times B_1 \to \real^{n+1}$ by
\begin{equation}
\label{eq:PsiFullHybridExtended}
\Psi_{e}(\tau,\xi_1):= \left[ \begin{array}{c} \tau \\ \Psi_{\tau}(\xi_1) \end{array} \right].
\end{equation}
By H-\ref{itm:H5}-b), $\forall~\xi_1 \in B_1$, $\Psi_{e}^-(\rfrac{T_p}{2},\xi_1) \in {\cal S}$. Hence, the loss of dimension is in the $\tau$-component, and therefore
$$ \dim \left( \Psi_{e}^-(\rfrac{T_p}{2},B_1) \cap {\cal S} \right) = \dim (B_1).$$

\item[(iii)] The Generalized Hybrid Zero Dynamics Manifold (G-HZD) is therefore\endnote{Modulo $T_p$ is not required here because $\tau$ is reset at $T_p$, whereas in the non-hybrid case, it was required in \eqref{eq:ZeroDynamicsManifold}.}
\begin{equation}
\label{eq:GeneralZeroDynamicsManifold}
Z_e:=\Psi_{e}([0, T_p),B_1),
\end{equation}
which has two components,
$$ Z_{e,i}:=\Psi_{e}([0, \rfrac{T_p}{2}),B_1)$$
and
$$Z_{e,ii}:=\Psi_{e}([\rfrac{T_p}{2}, T_p),\Psi_{e}^-(\rfrac{T_p}{2},B_1) \cap {\cal S}).$$

\item[(iv)] The corresponding restriction dynamics is given by \eqref{eq:GHZD}, which is then the G-HZD.

\end{enumerate}
\end{remark}

%

\subsection{Stabilizing the Original Model}
\label{sec:InnerOuterLoopHybrid}

We can now obtain and explain the controller we use on bipeds. Similar to Sect.~\ref{sec:ExtendedModels}, assume the continuous phase of the hybrid model has the form
\begin{equation}
\label{eq:TwoBlockSystemHybridReduced}
\begin{aligned}
\dot{x}_1&= f_1(x_1,x_2,u_1)\\
\dot{x}_2&= f_2(x_1,x_2,u_2),
\end{aligned}
\end{equation}
with
$$x_2=\left[ \begin{array}{c} x_{2a} \\ x_{2b}  \end{array} \right]~\text{and}~f_2= \left[ \begin{array}{c} x_{2b} \\ u_2 \end{array} \right],$$
and $u=(u_1, u_2)$. The reason to split the input and not allow the $u_2$-component to enter the $x_1$-dynamics will be clear shortly. We allow the $u_1$-component to be empty.

\begin{enumerate}
\setlength{\itemsep}{.15cm}
\setcounter{enumi}{5}
\renewcommand\labelenumi{{\bf H-\theenumi}}
\item \label{itm:H6} The solutions in H-\ref{itm:A5} also satisfy
\begin{equation}
\label{eq:boundaryConditionHybrid}
\gamma(\varphi_{1\xi_1}(T_p))=\varphi_{2 \xi_1}(T_p).
\end{equation}
\mbox{ } \hfill \(\blacksquare\)
\end{enumerate}

\begin{theorem}
\label{them:StabilityThanksToHybridZeroDynamics} Assume the open-loop hybrid system \eqref{eq:open_loop_hybrid_model_with_one_phase} with $f$ given by \eqref{eq:TwoBlockSystemHybridReduced} satisfies Hypotheses H-\ref{itm:H1} to H-\ref{itm:H3}, H-\ref{itm:H5} and H-\ref{itm:H6}. Define $B_1:=\{ \xi_1 \in \real^{n_1}~|~ (\xi_1,\gamma(\xi_1))\in B \}$.
Assume in addition there exist open sets $B^{e}_{\rm 1.i}$ and $B^{e}_{\rm 1.ii}$ that contain $B_1$, a $\delta>0$, and feedbacks
 \begin{align*}
 \nu_{\rm i}&:[0,\rfrac{T_p}{2}+\delta] \times B^{e}_{\rm 1.i} \to \real^{n_1} \\
 \nu_{\rm ii}&:[\rfrac{T_p}{2}-\delta, T_p] \times B^{e}_{\rm 1.ii} \to \real^{n_1}
 \end{align*}
 and
 \begin{align*}
 \mu_{\rm i}&:[0,\rfrac{T_p}{2}+\delta] \times B^{e}_{\rm 1.i} \to \real^m \\
 \mu_{\rm ii}&:[\rfrac{T_p}{2}-\delta, T_p] \times B^{e}_{\rm 1.ii} \to \real^m
 \end{align*}
that are piecewise continuous in $t$, locally Lipschitz continuous in $x$, and, such that, for $0\le t < T_p$ and $\xi_1 \in B_1$,
\begin{equation}
\label{eq:LearningConditionHybridTheorem_nu}
\begin{aligned}
\nu_{\rm i}(t,\varphi_{1\xi_1}(t))&=\varphi_{2\xi_1}(t),~\ 0 \le t < \rfrac{T_p}{2} \\
\nu_{\rm ii}(t,\varphi_{1\xi_1}(t))&=\varphi_{2\xi_1}(t),~\ \rfrac{T_p}{2} \le t < T_p
\end{aligned}
\end{equation}
and
\begin{equation}
\label{eq:LearningConditionHyrbidTheorem_mu}
\begin{aligned}
\mu_{\rm i}(t,\varphi_{1\xi_1}(t))&=u_{\xi_1}(t),~\ 0 \le t < \rfrac{T_p}{2} \\
\mu_{\rm ii}(t,\varphi_{1\xi_1}(t))&=u_{\xi_1}(t),~\ \rfrac{T_p}{2} \le t < T_p.
\end{aligned}
\end{equation}
Then for all $\frac{n_2}{2} \times \frac{n_2}{2} $ positive definite matrices $K_p$ and $K_d$, $\exists~\epsilon^\ast >0$, such that $\forall~ 0 < \epsilon \le \epsilon^\ast$, $x^*:[0, T_p]\to \real^{n_1+n_2}$ is a locally exponentially stable periodic solution of the closed-loop hybrid system
\begin{equation}
\label{eq:hybridSplitPhaseModelClosedLoopInnerOuter}
\begin{array}{l}
\Sigma_{\rm i}:\left\{
\begin{aligned}
\dot{\tau}&=1 \\
\dot{x}&=f(x_1,x_2,u_1,u_2) \hspace*{2.4 cm} x^-\notin\mathcal{S}\\
u_1&= \mu_{1i}(\tau,x_1) \\
u_2&= \mu_{2i}(\tau,x_1) - \left[\frac{K_p}{\epsilon^2}\frac{K_d}{\epsilon}\right]\big(x_2 - \nu_i(\tau, x_1) \big) \\
\tau^{+}&=\tau^- \\
x^{+}&=\Delta(x^-), \hspace*{3.5 cm}  x^-\in \mathcal{S}
\end{aligned}\right.
\\
\\
\Sigma_{\rm ii}:\left\{
\begin{aligned}
\dot{\tau}&=1 \hspace*{4.45 cm}  \tau^- < T_p\\
\dot{x}&=f(x_1,x_2,u_1,u_2) \\
u_1&= \mu_{1ii}(\tau,x_1) \\
u_2&= \mu_{2ii}(\tau,x_1) - \left[\frac{K_p}{\epsilon^2}\frac{K_d}{\epsilon}\right]\big(x_2 - \nu_{ii}(\tau, x_1) \big) \\
\tau^{+}&=0   \hspace*{4.45 cm}  \tau^- =  T_p\\
x^{+}&=x^-.
\end{aligned}\right.
\end{array}
\end{equation}
Moreover, the closed-loop system possesses a G-HZD and it is given by \eqref{eq:GHZD}.
\end{theorem}

\begin{remark}
The control was split so that the high-gain part of the feedback does not directly enter the states of the zero dynamics, namely $x_1$. This allows the system to conform with existing theorems \cite{AmGaGrSr2014} for establishing the exponential stability of the periodic orbit --- in the full-order hybrid model --- on the basis of its stability in the zero dynamics, \eqref{eq:GHZD}. Isolating the action of the high-gain controller to the $x_2$-dynamics was not necessary in the case of ODEs.
\end{remark}

\section{Bipedal Walking Gaits}
\label{sec:biped}

\begin{figure}[t!]
	\centering
	\includegraphics[width=1\columnwidth]{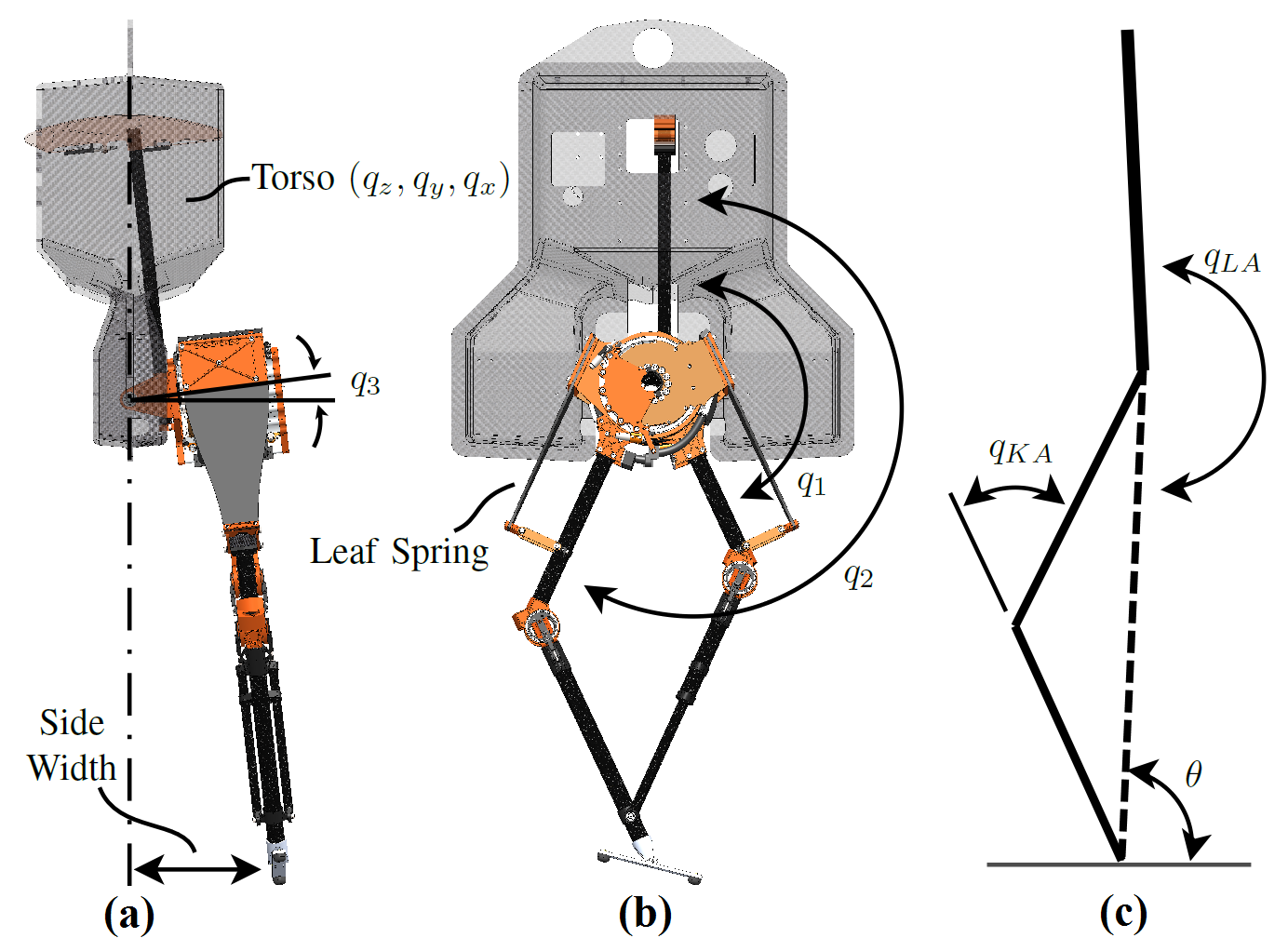}
	\caption{Based on \cite[Figure~2]{DaHartleyGrizzle2017}. Coordinate representation for MARLO, the Michigan copy of an Atrias-series robot.}
	\label{fig:MARLOv2}
\end{figure}

This section applies the Generalized Hybrid Zero Dynamics (G-HZD) developed in Section~\ref{sec:HybridSystemControl} to a bipedal robot, namely, the University of Michigan copy of an ATRIAS-series 3D robot that we call MARLO \cite{RAHUAKGR14}; see Figure~\ref{fig:MARLOv2}. As shown in \cite{BuHaGrGr16,GriffinIJRR2016, DaHartleyGrizzle2017, HartleyGrizzleCCTA2017}, the robot is capable of walking forward and backward at various speeds and over challenging terrain, both indoors and outdoors. The robot can step in place, but it cannot stand in place because its feet are passive. In the experiments, shoes are placed over the passive feet to prevent excessive yawing about the point of contact. The robot's hips have 2 DoF (pitch and roll). Because the hips lack yaw motion, turning is not one of the robot's strengths!

\begin{figure}[t!]
	\centering
	\includegraphics[width=1\columnwidth]{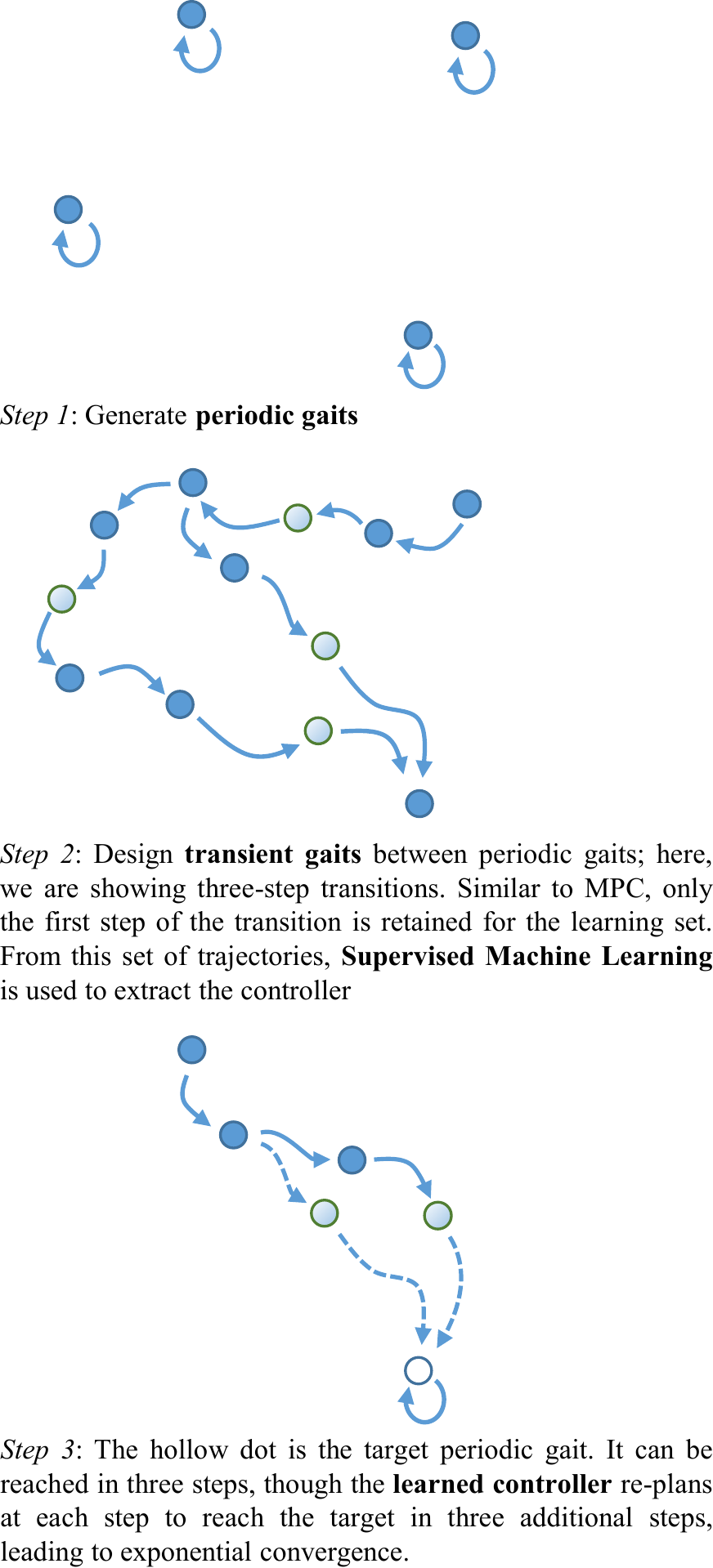}
	\caption{This can be thought of as an alternative representation of Figure~\ref{fig:OurApproach} when the surface $Z_0$ is built from periodic solutions of the full-order model. The light dots represent transient trajectories while the other dots (solid or hollow) are periodic orbits.}
	\label{fig:bipedControlDesignProcess}
\end{figure}

The control laws developed here will illustrate stepping in place, walking forward and backward, and transitions among such gaits. The work illustrates a theoretically sound method for gait design that unifies and significantly extends many of our previous results. The first control designs will rely on the optimization package in \cite{Jo2014}, which can only handle a planar model of the robot. For these designs, lateral stability is achieved via a heuristic foot placement policy. Since March of 2017, we have had access to the trajectory optimization package \cite{HeCoHuAm16}, which easily handles the full 3D model of the robot. Figure~\ref{fig:bipedControlDesignProcess} summarizes our control design process.

%

\subsection{MARLO}

The robot is described in detail in \cite{RAHUAKGR14}. For the planar model, the configuration variables are joint angles and one absolute coordinate. The angle $ \theta $, the absolute stance leg angle, is unactuated because the feet are passive. Hence, we define
$$x_1 = \begin{bmatrix} \theta \\ \dot{\theta} \end{bmatrix}.$$
For the purposes of controller design, the regulated quantities are $ q_a = (q_{x}, q_{sw, LA}, q_{st, KA}, q_{sw, KA}) $, that is, the torso, stance knee, swing leg, and swing knee angle, and hence
$$x_2 = \begin{bmatrix} q_a\\ \dot{q}_a \end{bmatrix}.$$
In the simulation and control design, we constrain the stance foot to remain in contact with the ground with no foot slip. In the experiments, we estimate the ground reaction forces through the deflection of the leaf springs to decide whether to control the torso angle or the stance leg angle $ q_{st, LA} $, as in \cite{ReHuJoPeVaJoAbHu15}. The model decomposition is done as in Appendix~\ref{sec:MechanicalModelDecomposition}.


\subsection{Design of Planar Periodic Gaits and Transition Trajectories via Optimization}

For robots, an orbit library is called a gait library. We first design a gait library
\begin{equation}
\label{eq:gaitLibraryPeriod}
	{\cal L}:=\{  \bar{v} ~|~-0.8 \le \bar{v} \le 0.8 \}
\end{equation}

consisting of periodic gaits for various average walking speeds
satisfying H-\ref{itm:H2}. In this example, we reuse the gaits described in \cite{DaHartleyGrizzle2017}, where each gait has period $ T_p = 0.4$,
and the cost function is
\begin{equation}
\label{eq:CostOneStep}
 J= \int_{-\rfrac{T_p}{2}}^{\rfrac{T_p}{2}} ||u(t)||^2 dt.
 \end{equation}
Denote the trajectory of the periodic gait and the corresponding input as $ \varphi^{\bar{v}}(t) $ and $ u^{\bar{v}}(t) $, respectively, and the midpoint of the periodic trajectory as $ x^{\bar{v}}_m $. The insertion function is built from the gait library and is denoted $\gamma_{\cal L}$.



For a given periodic orbit in $\cal L$, we define
\begin{equation}
\label{eq:domainOfxi1}
\begin{aligned}
B_1^{\bar{v}} = \{\xi_1:= (\theta, \dot{\theta}) ~|~ -\frac{\pi}{12} &\le |\theta - \theta_m^{\bar{v}}| \le \frac{\pi}{12},\\ -0.8 &\le |\dot{\theta} - \dot{\theta}_m^{\bar{v}}| \le 0.8 \},
\end{aligned}
\end{equation}
a sliding window\endnote{Because the legs are 1 m long, average walking speed and angular rate at the middle of the step are nearly the same.} about the target speed $ \bar{v}$. For $\xi_1 \in B_1^{\bar{v}} $, trajectories are generated as in H-\ref{itm:H5} using optimization with cost function
\begin{equation}
	\label{eq:costThreeSteps}
	J(\xi_1, \bar{v}) = \sum_{k = 1}^{6}\int_{(k-1)\frac{T_p}{2}}^{k\frac{T_p}{2}} 	\big( ||x-\varphi^{\bar{v}}||^2_Q + ||u - u^{\bar{v}}||^2_R \big)dt
\end{equation}
for a horizon of length $T_h= 3T_p $.  The optimization is performed subject to the hybrid dynamics describing MARLO, the physical constraints shown in Table.~\ref{tab:physicalConstraint}, and the terminal constraint $ x(3T_p) = \varphi^{\bar{v}}(T_p) $. So that the  trajectories can be used to stabilize the full-order model, the boundary constraints in H-\ref{itm:H6} is also imposed. The resulting trajectory and input for the given $\bar{v}$ and selected $\xi_1$ are denoted $ \varphi^{\bar{v}}_{\xi_1}(t) $ and $ u^{\bar{v}}_{\xi_1}(t) $, respectively.

\begin{table}
	\centering
	\caption {Optimization constraints}
	\label{tab:physicalConstraint}
	{\renewcommand{\arraystretch}{1.5}%
		\begin{tabular}{c l}
			Motor Toque $|u|$        & $<5$~Nm   \\
			Step Duration $T$			& $=0.4$~s\\
			Friction Cone $\mu$       & $<0.6$    \\
			Impact Impulse $Fe$       & $<15$~Ns  \\
			Vertical Ground Reaction Force & $>300$~N  \\
			Mid-step Swing Foot Clearance  & $>0.15$~m
	\end{tabular} }
	\vspace*{-.5cm}
\end{table}

\begin{remark}
		The initial condition set $ B_1^{\bar{v}} $ is related to the notion of 3-step capture region defined in \cite{KODBREGOPR12,ZaHaRu15}. In our experience, three-steps is a reasonable balance between planning horizon and computational burden.
\end{remark}


\subsection{Controller Design via Machine Learning}

The gait library \eqref{eq:gaitLibraryPeriod} is assumed to be discretized by 5 evenly spaced average speeds $\bar{v}_k$, each $\xi_1^i=: (\theta^i, \dot{\theta}^i)$ is drawn from  a uniformly spaced gird of 25 points, and time interval $[0, T_p]$ is evenly sampled into 21 points, $t_j$. The combined training and validation data set is therefore denoted by
\begin{equation}
\label{eq:denoteOrbitLearningBiped}
\begin{aligned}
x_1^{i,j,k} &:= \varphi^{\bar{v}_k}_{1\xi_{1}^i}(t_j) \\
\nu^{i,j,k} &:=\varphi^{\bar{v}_k}_{2\xi_{1}^i}(t_j) \\
\mu^{i,j,k} &:=u^{\bar{v}_k}_{\xi^i}(t_j).
\end{aligned}
\end{equation}
Any infeasible optimization problems\endnote{ For example, due to torque limits, there is no solution for $ \bar{v} = 0.8 $ and $ \xi_1 = (\theta_m^{\bar{v}} + \rfrac{\pi}{12}, 0.8 + \dot{\theta}_m^{\bar{v}}) $.} are removed from the data set before processing it by Supervised Machine Learning.

\begin{figure}[t!]
	\centering
	\includegraphics[width=0.9\columnwidth]{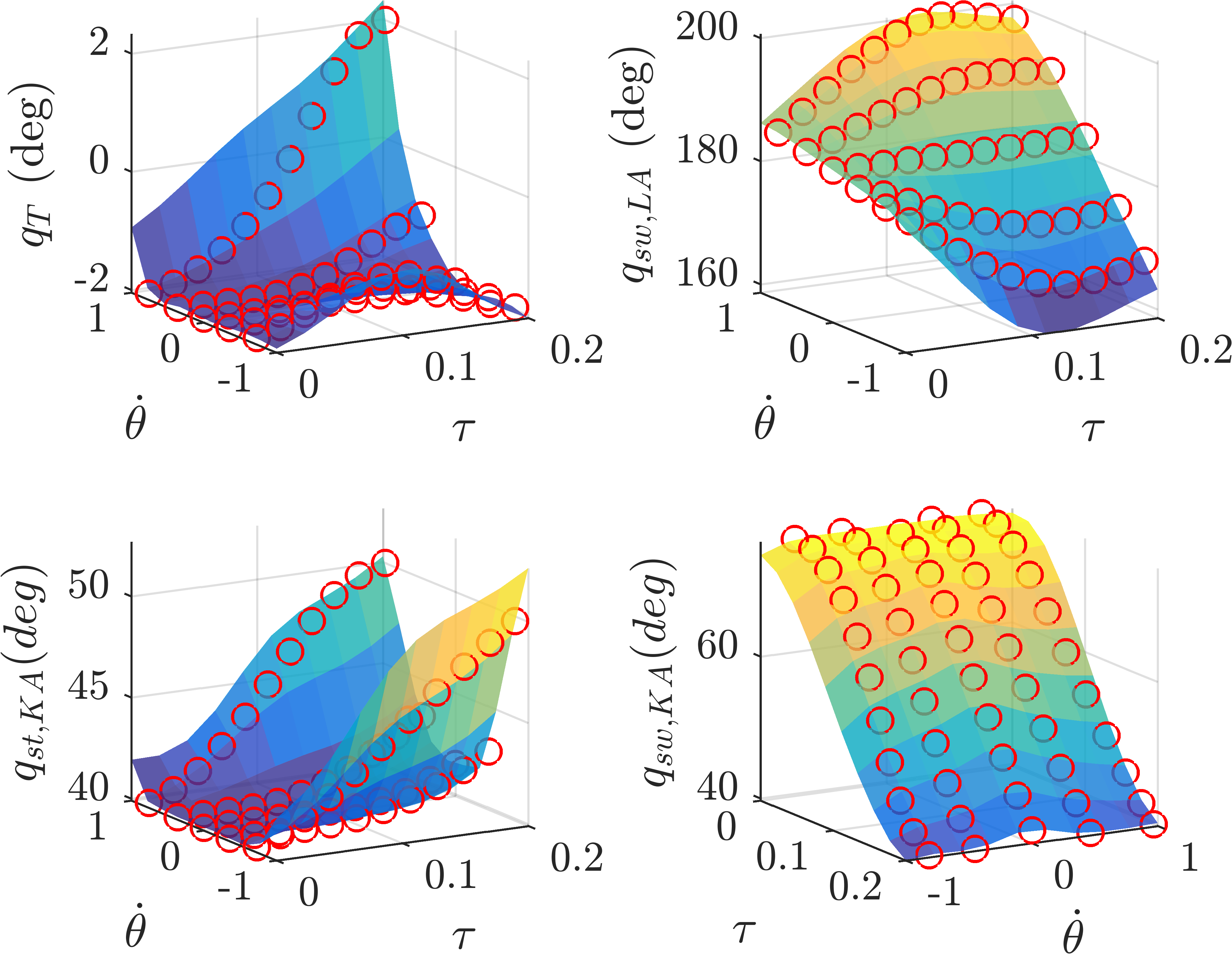}
	\caption{The fitting results of the Supervised Machine Learning. The features $ \theta = \theta_m^{\bar{v}} $ and $ \bar{v} = 0$ are fixed at constant values, while $ \tau \in [0, 0.2] $ is from mid-step to ground contact. The plots show four of the components in $ \nu_{\rm i}(t,x_1,{\bar{v}}) $.}
	\label{fig:BipedFitting}
\end{figure}

\begin{figure}[t!]
	\centering
	\includegraphics[width=0.9\columnwidth]{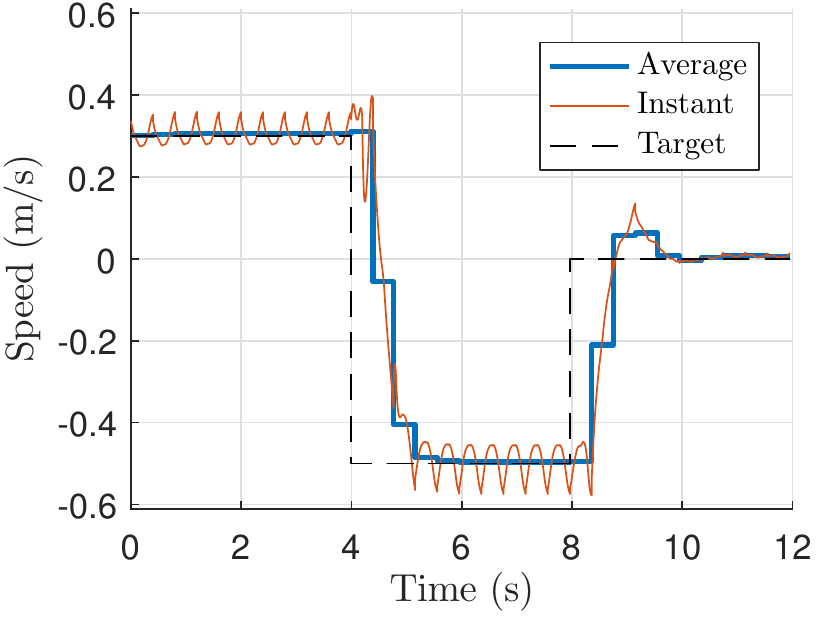}
	\caption{The target speed $ \bar{v} $ changes from $ 0.3 $~m/s to $ -0.5 $~m/s and to $ 0 $~m/s. The gait transition takes less than five steps to reach the target speed. The error between target and average speed is small.}
	\label{fig:transitionBiped}
\end{figure}

We next learn the functions in Theorem \ref{them:StabilityThanksToHybridZeroDynamics}. The features are $ (t_j, x_1^{i,j,k}, {\bar{v}}_k) $ and the labels are $ (\nu^{i,j,k}, \mu^{i,j,k}) $ and the data base is split at $ t = \rfrac{T_p}{2} $ so that the functions
$$
\begin{aligned}
\nu_{\rm i}(t,x_1,{\bar{v}}) \\
\nu_{\rm ii}(t,x_1,{\bar{v}}) \\
\mu_{\rm i}(t,x_1,{\bar{v}}) \\
\mu_{\rm ii}(t,x_1,{\bar{v}})
\end{aligned}
$$
are learned individually. Part of the fitting is shown in Figure~\ref{fig:BipedFitting}. These functions are enough to construct the G-HZD in \eqref{eq:GHZD}. To complete the control design as in \eqref{eq:hybridSplitPhaseModelClosedLoopInnerOuter}, the feedback gains $ K_p $, $ K_d $ and $ \epsilon $ must be selected. While in principle these last gains may have to vary with $\bar{v}$, for MARLO, we have found that a single set of gains\endnote{They essentially correspond the low-level PD gains which are straightforward to tune on MARLO. In simulation, $ K_p = 800$, $K_d = 40$ and $ \epsilon = 1$} suffices. The transition among different target speeds are shown in Figure~\ref{fig:transitionBiped}.

\subsection{Example Performance Analysis}

\subsubsection{Stability Analysis}

\begin{figure}[t!]
	\centering
	\includegraphics[width=0.9\columnwidth]{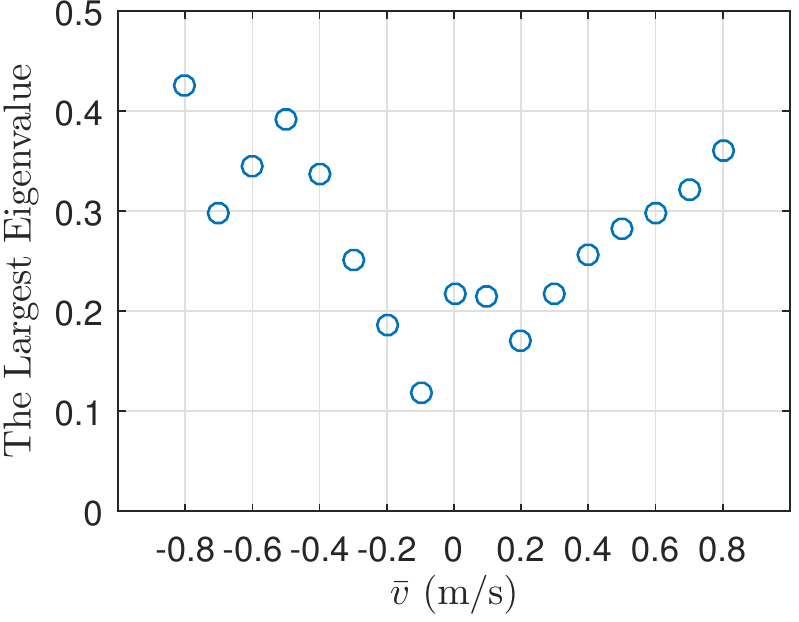}
	\caption{The largest eigenvalue of the Poincar\'e map is given for some target speeds. This indicates the target speed in the gait library \eqref{eq:gaitLibraryPeriod} is exponentially stable.}
	\label{fig:PoincareMapCollection}
\end{figure}

We know that the periodic orbits in the full-order model should be locally exponentially stable by the results in Sect.~\ref{sec:HybridSystemControl}. We formally verify this by numerically evaluating the Jacobian of the Poincar\'e map for twenty evenly-spaced points in the interval $-0.8 \le \bar{v} \le 0.8$. The magnitude of the largest eigenvalue is shown in Figure~\ref{fig:PoincareMapCollection}, which proves local exponential stability for each fixed target speed $\bar{v}$. Note that only five of these points were in the training data. The learned feedback functions have provided stable gaits for a continuum of target speeds.

\begin{figure}[t!]
	\centering
	\includegraphics[width=0.9\columnwidth]{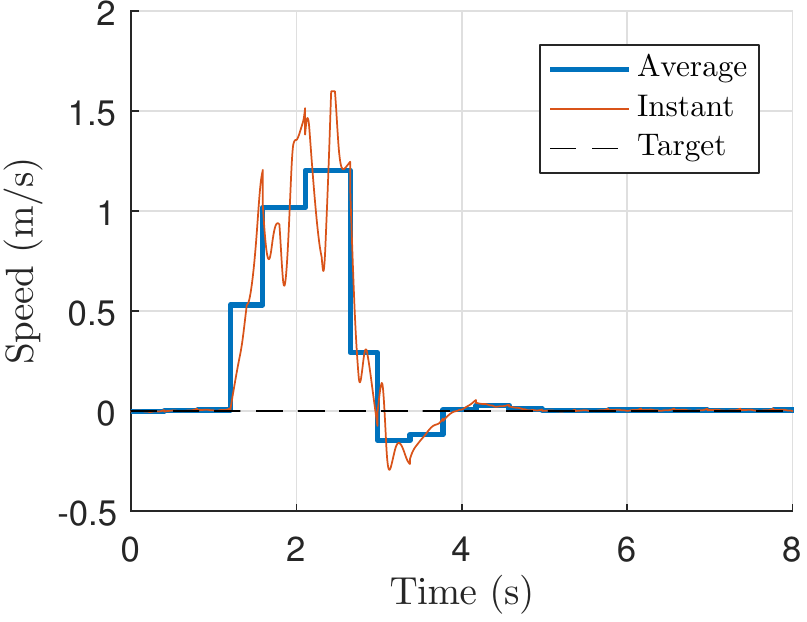}
	\caption{A perturbation is applied from 1.2 seconds to 2 seconds in a magnitude of $ 150 $~N. The maximum speed is larger than the maximum speed ($ 0.8 $~m/s) in the training set. The extrapolation may be the cause of the speed oscillation. Even though, the speed convergence near the target speed in less than five steps.}
	\label{fig:perturbationBiped}
\end{figure}

The stability of the overall closed-loop system is further illustrated by applying force perturbations, which is a more ``realistic'' test. We apply a longitudinal force on the hip at 1.2 seconds for 0.8 second (i.e., $ 2T_p $) and examine the time to recover the nominal gait. For the stepping-in-place gait $ \bar{v} = 0 $, the largest force from which the robot can recover without violating the physical constraints is $150$~N. Figure~\ref{fig:perturbationBiped} shows the resulting longitudinal velocity of the robot. The peak speed is approximately $1.5$~m/s, which is beyond the maximum training speed of $0.8$~m/s. The speed is once again less than $0.05$ m/s within five steps. The convergence rate is relatively fast given that the optimizer uses a horizon of three steps.

\subsubsection{Interpretation of the posture changes employed by the controller}

\begin{figure}[t!]
	\centering
	\includegraphics[width=0.9\columnwidth]{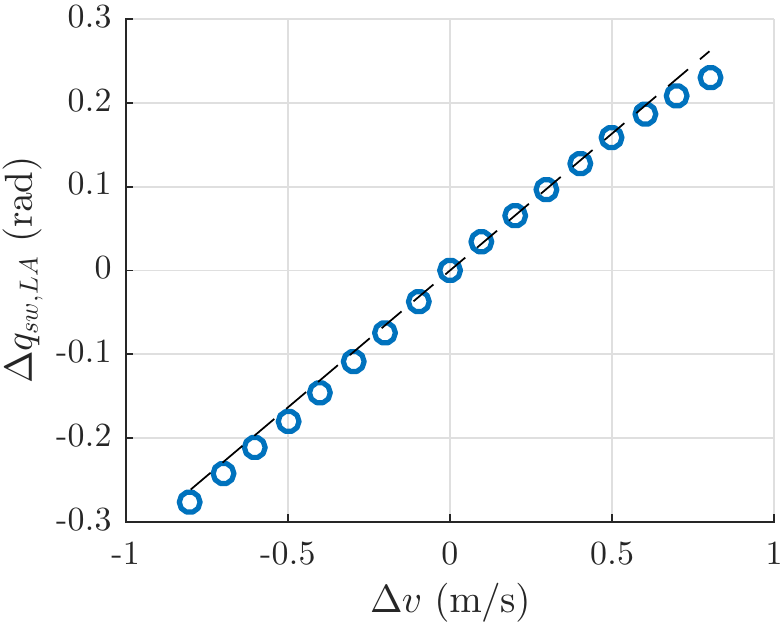}
	\caption{Change in swing leg angle vs change in velocity. One part of the learned-optimal strategy is a standard linear leg-angle adjustment policy as in \cite{PRTE06, dunn1996foot}.}
	\label{fig:optimal_swLAGain}
\end{figure}

\begin{figure}[t!]
	\centering
	\includegraphics[width=0.9\columnwidth]{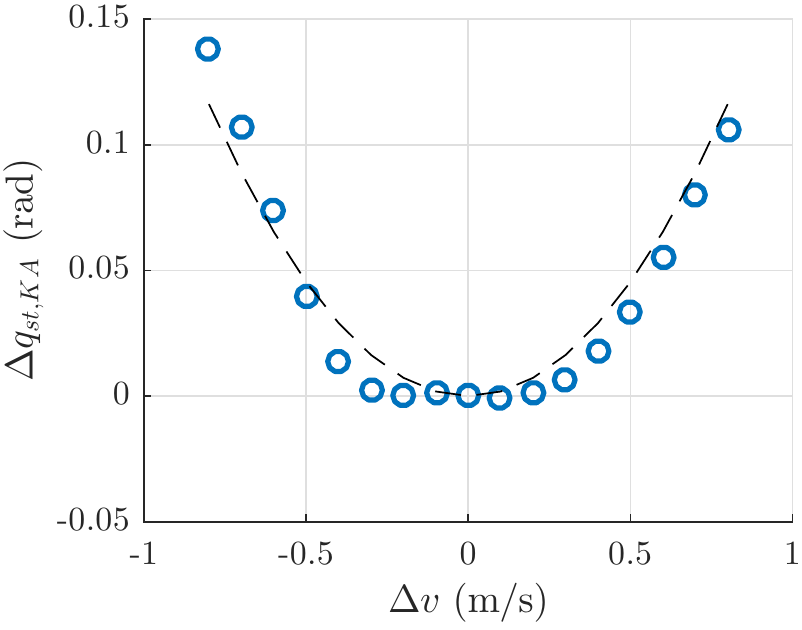}
	\caption{Change in stance knee angle vs change in velocity. This is not part of the standard recommendations in \cite{PRTE06, dunn1996foot}.}
	\label{fig:optimal_stKAGain}
\end{figure}

\begin{figure}[t!]
	\centering
	\includegraphics[width=0.9\columnwidth]{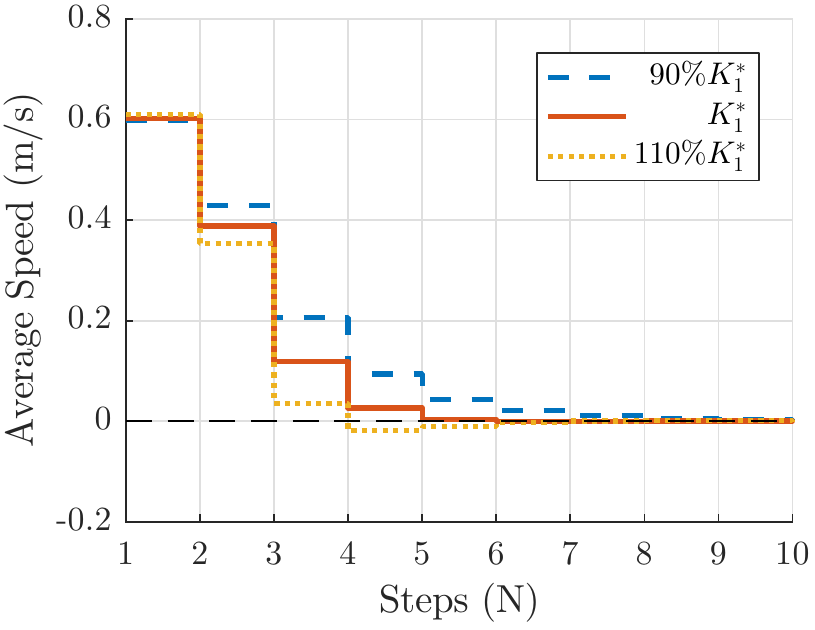}
	\caption{When the effective linear policy from the learned-optimal strategy is modified by $\pm 10\%$, the convergence to the nominal speed of zero either slows down or overshoots.}
	\label{fig:optimalPDgain}
\end{figure}

We now provide some physical intuition for how the controller coordinates the links to achieve stability. In fact, we evaluate $\nu(\rfrac{T_p}{2},x_1)$ with $x_1=(\rfrac{\pi}{12}, \dot{\theta})$, $-0.8 \le \dot{\theta} \le 0.8$. Figures~\ref{fig:optimal_swLAGain} and \ref{fig:optimal_stKAGain} show the changes in the swing leg angle and the stance knee angle at touchdown. The swing leg is seen to obey an approximately linear relationship with respect to velocity, just as in the foot placement controllers in \cite{PRTE06, dunn1996foot} designed on the basis of an inverted pendulum model or a linear inverted pendulum model, viz
$$\Delta q_{sw, LA} = K_1 \Delta v,$$
and the scalar $K_1$ is constant. Denote the regressed linear fit in Figure~\ref{fig:optimal_swLAGain} by $K^*_1$. We add $ \pm 10\% $ to $K^*_1$ and compare the resulting foot placement strategies in Figure~\ref{fig:optimalPDgain}. It is seen that with the smaller gain the velocity takes longer to settle whereas with the larger gain, there is overshoot.


\begin{figure}[t!]
	\centering
	\includegraphics[width=0.9\columnwidth]{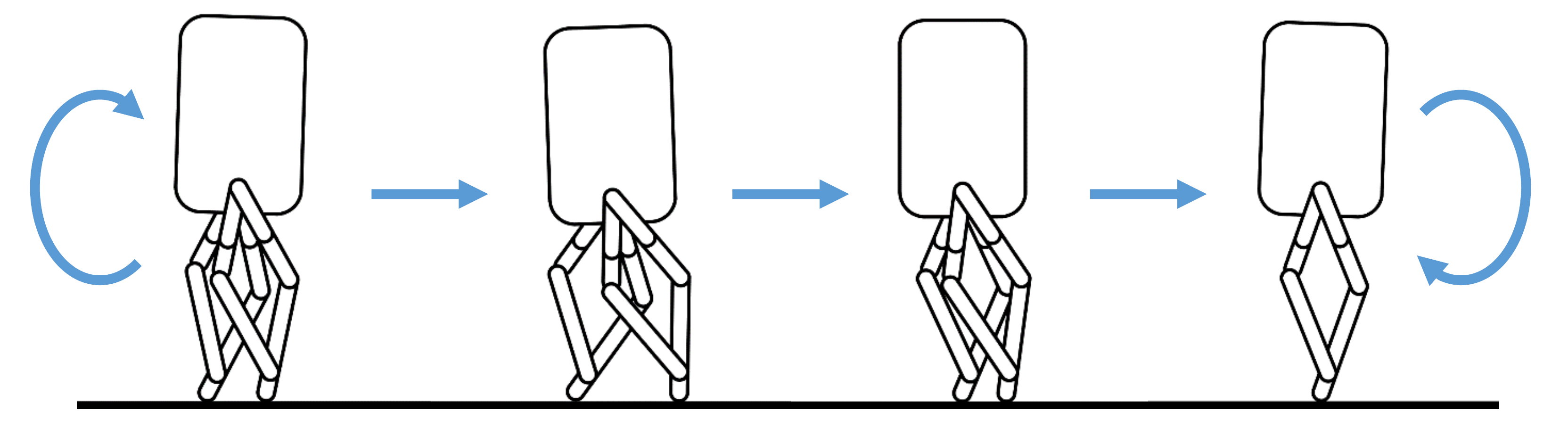}
	\caption{Stick figure showing the coordinated action of torso angle, knee bend, and leg angle provided by the learned-three-step optimization.}
	\label{fig:threeStepTransition}
\end{figure}

The learned controller is more than just providing foot placement. Figure~\ref{fig:optimal_stKAGain} also shows a quadratic relationship for knee angle versus velocity just before touchdown, viz
$$\Delta q_{st, KA} = K_2 (\Delta v)^2.$$
As the velocity moves away from zero in either direction, the stance knee angle increases. Perhaps this is to lower the center of mass and to make it easier for the swing foot to touch the ground. Furthermore, in addition to changing the swing leg angle, the learned controller also straightens the swing knee angle, thereby extending the foot further out. Finally, it also leans the torso backward, keeping the center of mass over the stance toe, as shown in Figure~\ref{fig:threeStepTransition}. These adjustments are all coordinated by the optimization and automatically extracted from the transition trajectories by the supervised learning. Unlike a classical foot placement controller that adjusts only swing leg angle, the learned controller uses the many degrees of freedom of the robot to achieve better performance.

%
%
%

\section{Experiments on a 3D Bipedal Robot}

This section extends the learning controller of the last section to the full-3D model of MARLO. Hence, both the sagittal and lateral planes of the robot are included, while yaw rotations are assumed to be small due to the foot. This control design will allow the physical 3D-robot to walk and step in place. The 3D-controller design will mostly follow the process of the planar example that was illustrated through simulations, though some modifications have been made during the experimental implementation to deal with model uncertainty, impact model uncertainty, and signal noise; these will be clearly explained and justified.

\subsection{3D-MARLO Configuration}

We still use Figure~\ref{fig:MARLOv2} to describe the generalized coordinates on the 3D robot. Let $ (p_x, p_y) $ denote the sagittal and lateral position of the center of mass, and be $ (v_x, v_y) $ be the corresponding velocities. We define
$$x_1 = \begin{bmatrix}
 p_x \\ p_y \\ v_x \\ v_y
\end{bmatrix}$$ as the ``weakly actuated'' states. The regulated angles are $ q_a = (q_{y}, q_{x}, q_{sw, 3}, q_{sw, LA}, q_{st, KA}, q_{sw, KA}) $, that is, the torso roll and pitch, swing hip, swing leg, stance knee and swing knee angles, and hence the ``strongly actuated'' states are
$$x_2 = \begin{bmatrix} q_a\\ \dot{q}_a \end{bmatrix}.$$
We note that the coordinates $ (x_1, x_2) $ describe the robot dynamics in \eqref{eq:TwoBlockSystemHybrid}, or \eqref{eq:SpongStateForm} after a coordinate change.

\subsection{Optimization}

We first use optimization, with a cost function as in \eqref{eq:CostOneStep}, to build a periodic gait library
\begin{equation}
\label{eq:gaitLibraryPeriod3D}
\begin{aligned}
{\cal L}:= \{(\bar{v}_x, \bar{v}_y) ~|~ -0.6 &\le \bar{v}_x \le 0.6,\\ -0.4 &\le \bar{v}_y \le 0.4 \},
\end{aligned}
\end{equation}
as a grid in two-dimensional Cartesian space; i.e., longitudinal and lateral speed of the robot. The gaits are designed, without loss of generality, such that the associated 2D-Cartesian positions $(p_x, p_y)$ are equal to the origin at a nominal point in each of the gaits. The insertion function associated to the gait library is constructed using linear regression as in the pendulum model \eqref{eq:linearRegressionOrbitGamma} and in the planar biped examples; specifically,
\begin{equation}
\label{eq:linearRegressionOrbitGamma2}
\begin{aligned}
x_2 = \gamma_{\cal L}(x_1) = a_0 + a_1 x_1,
\end{aligned}
\end{equation}
where $ a_0 $ and $ a_1 $ are constant vectors. A linear fitting is good enough for 3D MARLO. While one could do a more sophisticated fit, the maximum root-mean-square-error (RMSE) for all joints is less than 1 degree, and for all joint velocities it is less than 4 degrees per second, even with the linear fit. In part, this is a benefit of using the middle of the gait for building the controller.

The next step is to design transition trajectories from periodic orbits in the library $ {\cal L} $ to a target periodic orbit. In this paper, we will only illustrate the target orbit as the stepping-in-place gait, that is, $ (\bar{v}^*_x, \bar{v}^*_y) = 0 $. Further details are not given because they follow the planar example of the last section. One can also refer to work in \cite{DaHartleyGrizzle2017} for the design of transition gaits for different ground slopes.

We perform three-step trajectory optimizations as in \eqref{eq:costThreeSteps} and denote the collection of transition trajectories to the target speed of zero as $ \varphi_{{\cal L} \to 0} $. The orbit library is evenly sampled per $ \bar{v}_x = \{-0.6, -0.4, -0.2, \dots, 0.6\} $ and $ \bar{v}_x = \{-0.4, -0.3, -0.2, \dots, 0.4\} $, so that the total number of transition trajectories $ \varphi{^i}_{{\cal L} \to 0} $ is 63. The time interval $ [0, 0.4] $ is evenly sampled into 21 points, $ t_j $.

\subsection{Machine Learning}

For the 3D robot, we illustrate a different philosophy in building the feature set. Recall that in the planar examples, we included time and all of $x_1$-coordinates in the feature vector. Here, we ill select the feature vector as simply time and the velocity coordinates, namely,
$$ (t, v_x, v_y), $$
and leave out the Cartesian positions $(p_x, p_y)$. There are several reasons for this:
\begin{enumerate}
\item The impact map in the hybrid model \eqref{eq:hybridSplitPhaseModel} resets the Cartesian position variables to a constant, assumed to the origin. Hence, these variables do not need to be stabilized.
\item On the robot MARLO, we are placing the torso sagittal and roll angles in the $x_2$-coordinates, and hence if these are kept upright, the position of the center of mass does not provide significant additional information.
\item It keeps the dimension of the feature set as low as possible, which allows a smaller training data set.
\end{enumerate}

On 3D MARLO, the labels are taken as
$$ \nu = \varphi_{2, {\cal L} \to 0}, $$
which represents the $ x_2$-coordinates of the optimized trajectories. The control input $\mu$ is not learned because, as in previous experiments \cite{GriffinIJRR2016, DaHartleyGrizzle2017} on this robot, we use PD control
\begin{equation}
	\label{eq:PDcontrol3D}
	u_{PD} = \left[K_p, K_d\right]\big(x_2 - \nu(t, v_x, v_y)\big)
\end{equation}
to regulate the joint angles, without a feedforward term. The feedforward torque is not applied because of uncertainty in the model. Specifically, the model does not include the motor drive friction, which consumes about 20\% of the torque in nominal operation (stepping in place and walking), nor does the model include backlash or compliance in the harmonic drives. Finally, the leaf springs in Figure~\ref{fig:MARLOv2} are excluded from the model; they deflect about 5 degrees when supporting the robot's weight.

Since the impact happens at $ \rfrac{T_p}{2} = 0.2 $, the functions
$$
\begin{aligned}
\nu_{\rm i}(t, v_x, v_y)  \\
\nu_{\rm ii}(t, v_x, v_y) \\
\end{aligned}
$$
are learned individually.

\subsection{Experimental Implementation}

Another difference between the model and the physical robot occurs in the velocity signal. Due to spring deflection, impacts, and joint compliance, the estimated Cartesian velocities $ (v_x, v_y) $ are noisy even if each of the individual joint angular velocity signals is relatively ``clean''. We thus use a strong\endnote{The time constant is $ 0.2 $ second.} first-order filter to clean up the Cartesian velocity signals $ (v_x, v_y) $ appearing in $ \nu(t,v_x,v_y), $ the reference for low-level PD controllers \eqref{eq:PDcontrol3D}.

\begin{figure}[t!]
	\centering
	\includegraphics[width=0.9\columnwidth]{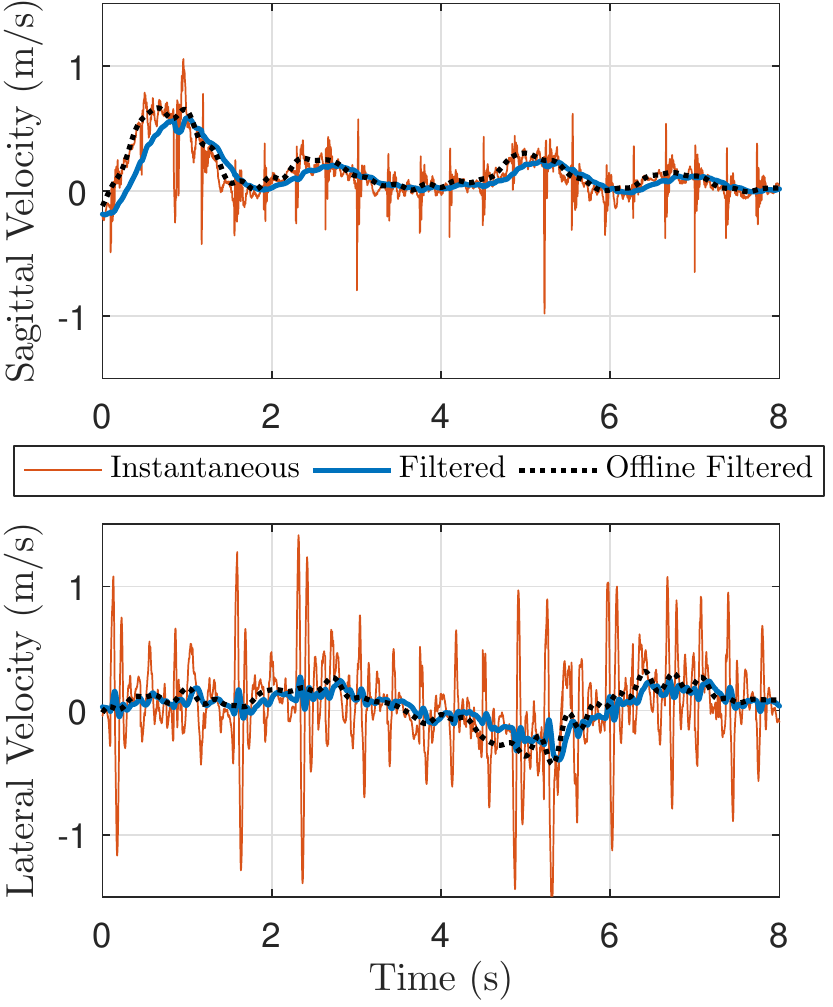}
	\caption{A perturbation is applied from 1.2 seconds to 2 seconds in a magnitude of $ 150 $~N. The maximum speed is larger than the maximum speed ($ 0.8 $~m/s) in the training set. The extrapolation may be the cause of the speed oscillation. Even though, the speed convergence near the target speed in less than five steps.}
	\label{fig:velocity_raw}
\end{figure}

The filtered signal, shown in Figure~\ref{fig:velocity_raw}, is relatively clean, but causes phase lag. Moreover, the energy loss at impact is less on the robot than in the design model because the springs store energy at impact release it throughout a step. This two factors lead the learned controller to generate overshoot in the Cartesian velocities. We mitigate the overshoot by introducing a speed-damping term on the torque of the swing leg and the swing hip,
\begin{equation}
\label{eq:DerivativeTerm}
\begin{aligned}
u^{sw, LA}_{d} &= N_{x,d} (v_x - v_x^{k})  \\
u^{sw, 3}_{d} &= N_{x,d} (v_y - v_y^{k}),
\end{aligned}
\end{equation}
which is the same term used in \cite{ReHuJoPeVaJoAbHu15}. The overall torque is
\begin{equation}
\label{eq:torque3D}
u = u_{PD} + u_{d}.
\end{equation}

\subsection{Results}

The learned controller for stabilizing the stepping-in-place gait was implemented on MARLO. The nominal Cartesian velocities are thus zero. Forces were applied in the sagittal and lateral directions, or a mix of the two, by the experimenter applying a push or a kick to the robot. The amount of force has not been estimated, but the reader can judge of his-or-herself based on the experiment video (see also Extension~\hyperlink{ex:video}{2}).

\begin{figure}[t!]
	\centering
	\includegraphics[width=0.9\columnwidth]{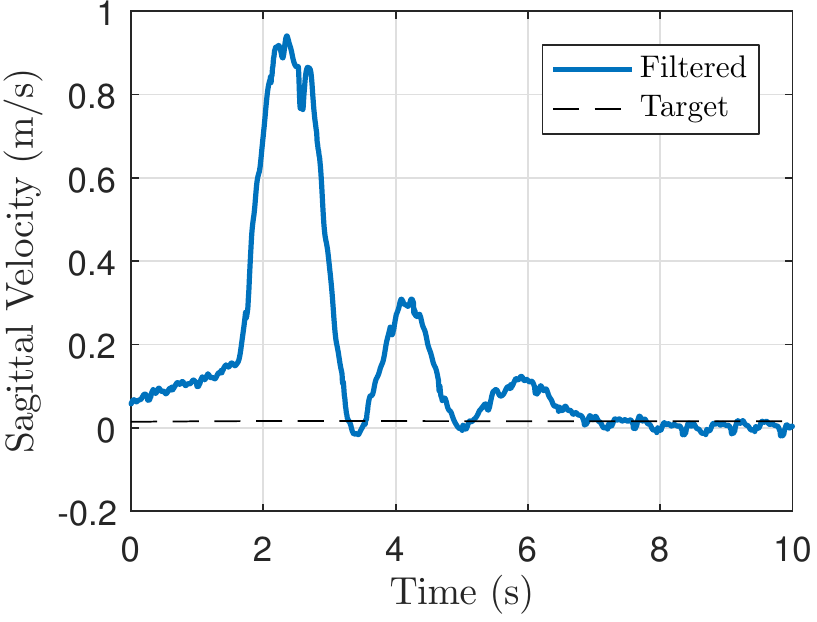}
	\caption{An example of the velocity response for a kick in the forward (sagittal) direction. The perturbation is applied at around 1 second driving the robot to peak speed of 0.9 m/s.}
	\label{fig:velocity_x}
\end{figure}

In the first experiment, five successive kicks were applied to MARLO in the forward (sagittal) direction. MARLO consistently recovered from the disturbances. The peak speed varied from 0.8 to 1 m/s, as shown in Figure~\ref{fig:velocity_x}. A harder kick was not applied since the training set only includes speed up to 0.6 m/s. After reaching the peak speed, MARLO slowed down to 0.1 m/s in less than 5 seconds. The robot acted as an underdamped spring-load system. This may be caused by the leg springs in the physical which compressed and unloaded on the second step, while the model did not include this effect.  We have added the damping term in \eqref{eq:DerivativeTerm} to mitigate the overshoot effect. A larger derivative gain will further reduce the overshoot, but will so increase the settling time. Still, the speed slowed down to 0.1 m/s in less than 5 seconds. The leg motor torque is shown in Figure~\ref{fig:velocity_x}. The torque bound (5 Nm) was reached when robot moved around the peak speed. This could explain why the optimization can only find the solution of transition gait from 0.6 m/s to zero.

\begin{figure}[t!]
	\centering
	\includegraphics[width=0.9\columnwidth]{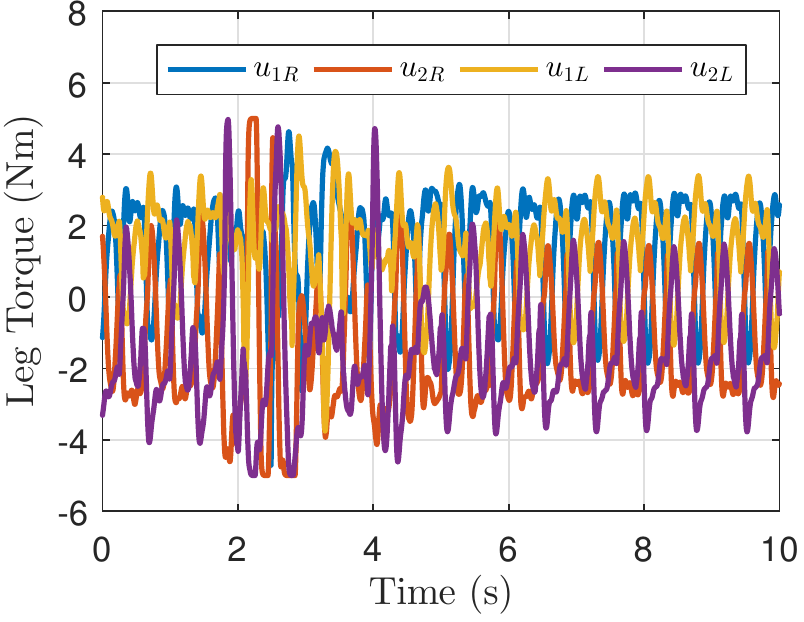}
	\caption{The leg torques are applied before a 1:50 gear transmission. The torques $ (u_{1R}, u_{2R}, u_{2L}, u_{2L}) $ are associate with the leg joints $ (q_{1R}, q_{2R}, q_{1L}, q_{2L}) $ in the robot configuration, Figure~\ref{fig:MARLOv2}. }
	\label{fig:torque_x}
\end{figure}

\begin{figure}[t!]
	\centering
	\includegraphics[width=0.9\columnwidth]{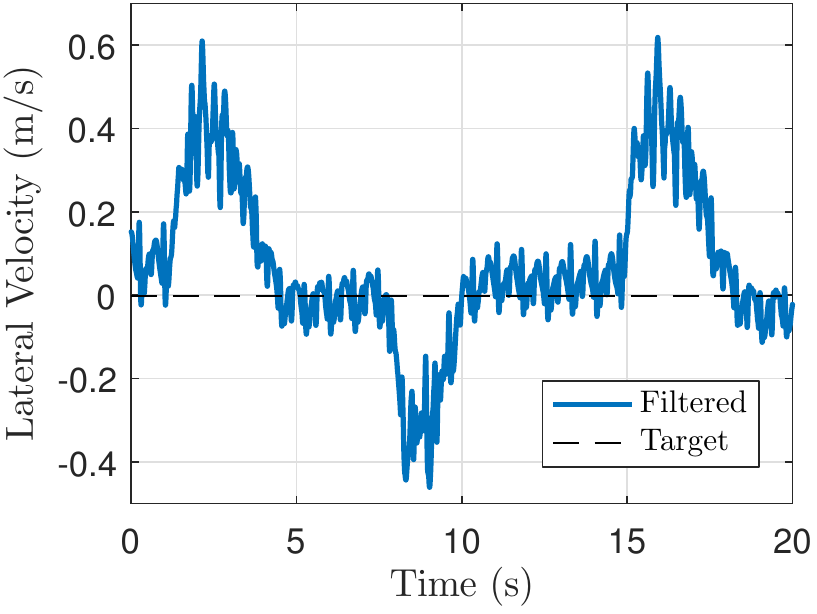}
	\caption{An example of the velocity response for multiple pushes in the lateral direction. The positive sign is the right direction whereas the negative is the left.}
	\label{fig:velocity_y}
\end{figure}

\begin{figure}[t!]
	\centering
	\includegraphics[width=0.9\columnwidth]{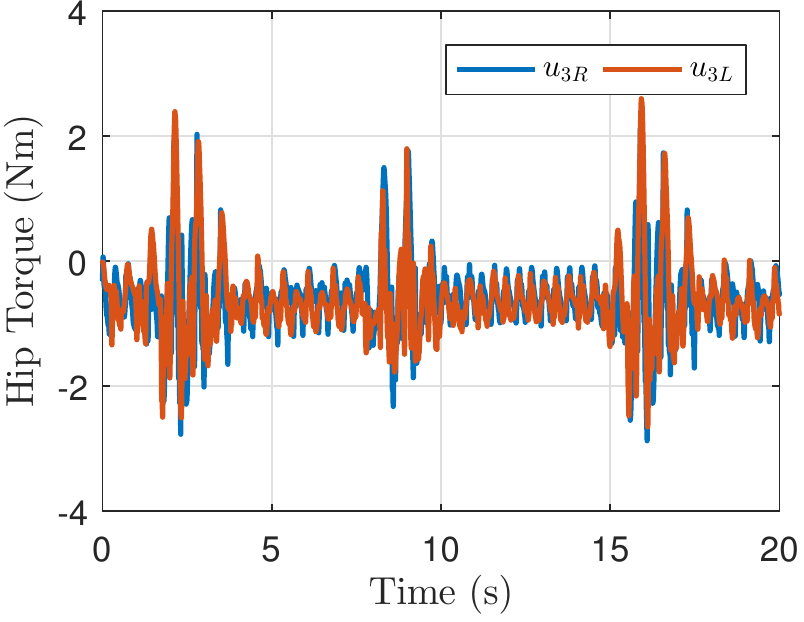}
	\caption{The hip torques are applied before a 1:27.5 gear-belt drive. The positive torque is to move leg inward whereas the positive torque is to move leg outward.
	}
	\label{fig:torque_y}
\end{figure}

The second experiment was to push MARLO in the lateral direction. Because of the parallelogram shape of the legs, one foot cannot place across the other, which limits the available range of foot placement. Plus the weaker motor on the hip, the lateral stability of MARLO is weaker than the sagittal direction. The push drove the lateral speed to 0.6 m/s, which is higher than the training speed 0.4 m/s. The push was applied on both left and right sides, shown in Figure~\ref{fig:velocity_y}. The hip motor torque is shown in Figure~\ref{fig:torque_y}. The torque bound is 3 Nm.

Random direction pushes and kicks were included in the last part of the experiment video. We applied force to move MARLO backward and to turn around.

\section{Discussion and Conclusions}
\label{sec:Conclusions}

\subsection{Strategy}
This paper is building on the recent revolution in open-loop trajectory optimization. It is now possible to compute in minutes gaits that used to take us hours or more. Armed with a set of open-loop trajectories, the question we posed was, how to turn them into a feedback controller? Our strategy was to attempt to build a surface from the trajectories and to induce a vector field on that surface that had desirable properties, such as (1) it contained a periodic solution of the model that met important physical constraints; (2) trajectories on the surface, by design, converged to the periodic solution; and (3), we could find a feedback controller for the complete model of the system that would render this surface exponentially attractive.

We have used Supervised Machine Learning as a ``universal'' function approximation, in other words, glorified regression. The functions we sought were implicitly contained in the data, and to our knowledge, closed-form solutions seem unlikely to exist. Hence, they had to be computed numerically in one fashion or another, and we believe an important contribution of the paper is to show how the functions needed to build a feedback controller can be extracted from a set of trajectories over a fixed time interval.

\subsection{Curse of Dimensionality} The method we use to build a vector field from open-loop trajectories works, at least in principle, in large dimensions. Even with the large strides made in optimization, high-dimensional state spaces are still the bottleneck, at least with our approach. Hence, it was important to find a structural property in our robot models that would allow us to focus the optimization effort on a low-dimensional portion of the system. We chose to exploit the local input-output linearizability of the actuators associated with knees and hips for example and put into the ``weakly-actuated category'' things like the global orientation of the body and possibly ankle joints. This allowed us to build trajectories of the full-model parameterized by initial conditions of a small subsystem, without making any approximations. In the end, we do build the control law for the full model around a low-dimensional model, just as advocates of pendulum models do, however, and this is important, all of the solutions of our low-dimensional model are feasible solutions of the robot itself and they meet whatever constraints were included in the design of the trajectories.

\subsection{Original HZD vs G-HZD} Once one understands how the G-HZD tool works, it's hard to believe how much the previous method could accomplish with a single optimization. The work presented in \cite{WGCCM07} uses a single optimization to design the periodic orbit. And that is it. For robots with one degree of underactuation and for which a ``mechanical phase variable'' can be found, that is a strictly monotonically increasing generalized coordinate along the periodic orbit, the basic HZD result in \cite{WGCCM07} shows how to build an invariant surface, relate stability of the periodic orbit in the surface to a physical property of the periodic orbit\endnote{The velocity of the robots center of mass should point downward at the end of the step.}, and how to render the surface exponentially attractive. G-HZD can handle more than one degree of underactuation. As in \cite{ReCoHeHuAm16,Reher2016Algorithmic,kolathaya2017parameter}, G-HZD includes ``time'' as a monotonically increasing generalized coordinate. What makes it quite distinct from these references is that the equivalent of ``G-HZD virtual constraints'', the function $\nu(t,x_1)$, depends on time and the full state of the zero dynamics, thereby enriching the set of possible closed-loop behaviors.

With the previous work on HZD and its extensions in \cite{BuHaGrGr16,Kaveh2016BMI,GriffinIJRR2016}, we were unable to handle challenging terrain such as the Michigan Wave Field . This motivated the introduction of a family of periodic orbits in \cite{DaHaHaGrGr16} and a first attempt at including transitions among the periodic orbits in \cite{DaHartleyGrizzle2017} as a means of building in stability. While this latter paper also used Supervised Machine Learning to design a feedback function, it also made some false steps relying on analogies with model predictive control that were not supported by deeper analysis. The present paper is our attempt to provide a consistent design framework. An attractive feature of the original HZD approach is that it has an easily verifiable set of sufficient conditions for its many results. We hope in some future paper, a similarly clean analytical framework will be developed for G-HZD.

%

\subsection{Stability Mechanism} Not only do the new results handle higher degrees of underactuation but even in the case of one degree of under actuation, the way stability is achieved is quite different in G-HZD vs. HZD. As discussed in the Introduction, with G-HZD, one does not have to count on the impacts to create stability. More general stability mechanisms, such as foot placement, or as shown in Figure~\ref{fig:optimal_stKAGain}, lowering the center of mass, naturally arise.

\subsection{Future Work}

We see this paper as a first cut in developing a happy marriage among trajectory optimization, machine learning, and geometric nonlinear control. We hope the results in the paper can be reinforced with easy-to-check sufficient conditions for our many assumptions. Beyond these technical considerations, we also see several other directions. The recent work in \cite{ConsoliniCoordinateFree17} may provide a geometric formulation for building the invariant surface, which would also clarify the choice of coordinates for making the mapping $\Psi$ in \eqref{eq:PsiFull} and its projection to be full rank. We believe the feedback linearization assumption on the ``strongly actuated'' part of the dynamics can be weakened considerably. Replacing the terminal condition in the optimization with a terminal penalty is another direction that needs to be investigated.

\begin{acks}
Professor Jonathan Hurst (Oregon State University), Mikhail Jones (Agility Robotics) and the whole team in the Dynamic Robotics Laboratory (Oregon State University) are sincerely thanked for sharing their copy of ATRIAS. The experiments reported here have also benefited greatly from the participation of Ph.D. students Omar Harib, Ross Hartley, and Yukai Gong (University of Michigan). The optimization and gait generation for the 3D walking and experiments were contributed by a postdoctoral researcher Ayonga Hereid (University of Michigan).
\end{acks}

\begin{dci}
The authors declare that there is no conflict of interest.
\end{dci}

\begin{funding}
	This work was supported in part by NSF Grants CPS-1239037 and NRI-1525006, and in part by Toyota Research Institute under award TRI-N021515;  however this article solely reflects the opinions and conclusions of its authors and not any of the funding entities.
\end{funding}

\theendnotes

\bibliographystyle{SageH}

\begin{thebibliography}{86}
\providecommand{\natexlab}[1]{#1}
\providecommand{\url}[1]{\texttt{#1}}
\providecommand{\urlprefix}{URL }
\expandafter\ifx\csname urlstyle\endcsname\relax
  \providecommand{\doi}[1]{DOI:\discretionary{}{}{}#1}\else
  \providecommand{\doi}{DOI:\discretionary{}{}{}\begingroup
  \urlstyle{rm}\Url}\fi

\bibitem[{Aghasadeghi et~al.(2013)Aghasadeghi, Zhao, Hargrove, Ames, Perreault
  and Bretl}]{aghasadeghi2013learning}
Aghasadeghi N, Zhao H, Hargrove LJ, Ames AD, Perreault EJ and Bretl T (2013)
  Learning impedance controller parameters for lower-limb prostheses.
\newblock In: \emph{Intelligent robots and systems (IROS), 2013 IEEE/RSJ
  international conference on}. IEEE, pp. 4268--4274.

\bibitem[{Agrawal et~al.(2017)Agrawal, Harib, Hereid, Finet, Masselin, Praly,
  Ames, Sreenath and Grizzle}]{Agrawal2017First}
Agrawal A, Harib O, Hereid A, Finet S, Masselin M, Praly L, Ames AD, Sreenath K
  and Grizzle JW (2017) First steps towards translating {HZD} control of
  bipedal robots to decentralized control of exoskeletons.
\newblock \emph{{IEEE} Access} 5: 9919--9934.
\newblock \doi{10.1109/ACCESS.2017.2690407}.

\bibitem[{Ames et~al.(2014)Ames, Galloway, Grizzle and
  Sreeenath}]{AmGaGrSr2014}
Ames A, Galloway K, Grizzle JW and Sreeenath K (2014) Rapidly exponentially
  stabilizing control {L}yapunov functions and hybrid zero dynamics.
\newblock \emph{IEEE Transactions on Automatic Control} 59(4): 876--891.

\bibitem[{Ames(2014)}]{Ames2014Human}
Ames AD (2014) Human-inspired control of bipedal walking robots.
\newblock \emph{IEEE Transactions on Automatic Control} 59(5): 1115--1130.
\newblock \doi{10.1109/TAC.2014.2299342}.

\bibitem[{Ames et~al.(2012)Ames, Cousineau and Powell}]{Ames2012Dynamically}
Ames AD, Cousineau EA and Powell MJ (2012) Dynamically stable bipedal robotic
  walking with {NAO} via human-inspired hybrid zero dynamics.
\newblock In: \emph{Proceedings of the 15th ACM international conference on
  Hybrid Systems: Computation and Control}. ACM, pp. 135--144.

\bibitem[{Apostolopoulos et~al.(2015)Apostolopoulos, Leibold and
  Buss}]{ApLeBu15}
Apostolopoulos S, Leibold M and Buss M (2015) Settling time reduction for
  underactuated walking robots.
\newblock In: \emph{Intelligent Robots and Systems (IROS), 2015 IEEE/RSJ
  International Conference on}. IEEE, pp. 6402--6408.

\bibitem[{Astolfi and Ortega(2003)}]{ASOR03}
Astolfi A and Ortega R (2003) Immersion and invariance: a new tool for
  stabilization and adaptive control of nonlinear systems.
\newblock \emph{IEEE Transactions on Automatic Control} 48(4): 590 -- 606.

\bibitem[{Azevedo et~al.(2002)Azevedo, Poignet and Espiau}]{AzPoEs02}
Azevedo C, Poignet P and Espiau B (2002) Moving horizon control for biped
  robots without reference trajectory.
\newblock In: \emph{Robotics and Automation, 2002. Proceedings. ICRA '02. IEEE
  International Conference on}, volume~3. pp. 2762--2767.
\newblock \doi{10.1109/ROBOT.2002.1013650}.

\bibitem[{Azevedo et~al.(2004)Azevedo, Poignet and Espiau}]{AzPoEs04}
Azevedo C, Poignet P and Espiau B (2004) Artificial locomotion control: from
  human to robots.
\newblock \emph{Robotics and Autonomous Systems} 47(4): 203--223.

\bibitem[{Bainov and Simeonov(1989)}]{BASI89}
Bainov D and Simeonov P (1989) \emph{Systems with Impulse Effects : Stability,
  Theory and Applications}.
\newblock Chichester: Ellis Horwood Limited.

\bibitem[{Buss et~al.(2016)Buss, Hamed, Griffin and Grizzle}]{BuHaGrGr16}
Buss BG, Hamed KA, Griffin BA and Grizzle JW (2016) Experimental results for
  {3D} bipedal robot walking based on systematic optimization of virtual
  constraints.
\newblock In: \emph{American control conference}.

\bibitem[{Buss et~al.(2014)Buss, Ramezani, {Akbari Hamed}, {Griffin, B. A.,
  Galloway} and Grizzle}]{Buss2014Preliminary}
Buss BG, Ramezani A, {Akbari Hamed} K, {Griffin, B A, Galloway} KS and Grizzle
  JW (2014) Preliminary walking experiments with underactuated {3D} bipedal
  robot {MARLO}.
\newblock In: \emph{Intelligent Robots and Systems (IROS)}. IEEE, pp.
  2529--2536.

\bibitem[{Chevallereau et~al.(2003)Chevallereau, Abba, Aoustin, Plestan,
  Westervelt, Canudas-De-Wit and Grizzle}]{Chevallereau2003RABBIT}
Chevallereau C, Abba G, Aoustin Y, Plestan F, Westervelt ER, Canudas-De-Wit C
  and Grizzle JW (2003) Rabbit: a testbed for advanced control theory.
\newblock \emph{IEEE Control Systems} 23(5): 57--79.
\newblock \doi{10.1109/MCS.2003.1234651}.

\bibitem[{Clarke et~al.(1997)Clarke, Ledyaev, Sontag and
  Subbotin}]{clarke1997asymptotic}
Clarke FH, Ledyaev YS, Sontag ED and Subbotin AI (1997) Asymptotic
  controllability implies feedback stabilization.
\newblock \emph{IEEE Transactions on Automatic Control} 42(10): 1394--1407.

\bibitem[{Consolini et~al.(2016)Consolini, Costalunga and
  Maggiore}]{ConsoliniCoordinateFree17}
Consolini L, Costalunga A and Maggiore M (2016) A coordinate-free theory of
  virtual holonomic constraints.
\newblock \emph{arXiv preprint arXiv:1709.07726 [math.OC]} .

\bibitem[{Coron(1994)}]{co94}
Coron JM (1994) On the stabilization of controllable and observable systems by
  an output feedback law.
\newblock \emph{Mathematics of Control, Signals, and Systems (MCSS)} 7(3):
  187--216.

\bibitem[{Coron et~al.(1995)Coron, Praly and Teel}]{coron1995feedback}
Coron JM, Praly L and Teel A (1995) Feedback stabilization of nonlinear
  systems: Sufficient conditions and lyapunov and input-output techniques.
\newblock In: \emph{Trends in control}. Springer, pp. 293--348.

\bibitem[{Da et~al.(2016)Da, Harib, Hartley, Griffin and
  Grizzle}]{DaHaHaGrGr16}
Da X, Harib O, Hartley R, Griffin B and Grizzle JW (2016) From {2D} design of
  underactuated bipedal gaits to {3D} implementation: Walking with speed
  tracking.
\newblock \emph{IEEE Access} 4: 3469--3478.
\newblock \doi{10.1109/ACCESS.2016.2582731}.

\bibitem[{Da et~al.(2017)Da, Hartley and Grizzle}]{DaHartleyGrizzle2017}
Da X, Hartley R and Grizzle JW (2017) Supervised learning for stabilizing
  underactuated bipedal robot locomotion, with outdoor experiments on the wave
  field.
\newblock In: \emph{2017 IEEE International Conference on Robotics and
  Automation (ICRA)}. pp. 3476--3483.
\newblock \doi{10.1109/ICRA.2017.7989397}.

\bibitem[{Dunn and Howe(1996)}]{dunn1996foot}
Dunn ER and Howe RD (1996) Foot placement and velocity control in smooth
  bipedal walking.
\newblock In: \emph{Robotics and Automation, 1996. Proceedings., 1996 IEEE
  International Conference on}, volume~1. IEEE, pp. 578--583.

\bibitem[{Embry et~al.(2016)Embry, Villarreal and Gregg}]{EmViGr16}
Embry KR, Villarreal DJ and Gregg RD (2016) A unified parameterization of human
  gait across ambulation modes.
\newblock In: \emph{Submit to: International Conference of the IEEE Engineering
  in Medicine and Biology Society (EMBC)}. IEEE.

\bibitem[{Erez et~al.(2013)Erez, Lowrey, Tassa, Kumar, Kolev and
  Todorov}]{ErLoTaKuKoTo13}
Erez T, Lowrey K, Tassa Y, Kumar V, Kolev S and Todorov E (2013) An integrated
  system for real-time model predictive control of humanoid robots.
\newblock In: \emph{2013 13th IEEE-RAS International Conference on Humanoid
  Robots (Humanoids)}. pp. 292--299.
\newblock \doi{10.1109/HUMANOIDS.2013.7029990}.

\bibitem[{Faraji et~al.(2014)Faraji, Pouya, Atkeson and Ijspeert}]{FaPoAtIj14}
Faraji S, Pouya S, Atkeson CG and Ijspeert AJ (2014) Versatile and robust {3D}
  walking with a simulated humanoid robot ({Atlas}): a model predictive control
  approach.
\newblock In: \emph{2014 IEEE International Conference on Robotics and
  Automation (ICRA)}. IEEE, pp. 1943--1950.

\bibitem[{Full and Koditschek(1999)}]{FUKO99}
Full R and Koditschek D (1999) Templates and anchors: Neuromechanical
  hypotheses of legged locomotion on land.
\newblock \emph{Journal of Experimental Biology} 202: 3325--3332.

\bibitem[{Gregg et~al.(2014)Gregg, Lenzi, Hargrove and Sensinger}]{GrLeHaSe14}
Gregg RD, Lenzi T, Hargrove LJ and Sensinger JW (2014) Virtual constraint
  control of a powered prosthetic leg: From simulation to experiments with
  transfemoral amputees.
\newblock \emph{IEEE Transactions on Robotics} .

\bibitem[{Griffin and Grizzle(2016)}]{GriffinIJRR2016}
Griffin B and Grizzle J (2016) Nonholonomic virtual constraints and gait
  optimization for robust walking control.
\newblock \emph{The International Journal of Robotics Research} :
  0278364917708249.

\bibitem[{Grizzle et~al.(2014)Grizzle, Chevallereau, Sinnet and
  Ames}]{grizzle2014models}
Grizzle JW, Chevallereau C, Sinnet RW and Ames AD (2014) Models, feedback
  control, and open problems of {3D} bipedal robotic walking.
\newblock \emph{Automatica} 50(8): 1955--1988.

\bibitem[{Hamed et~al.(2016)Hamed, Buss and Grizzle}]{Kaveh2016BMI}
Hamed KA, Buss BG and Grizzle JW (2016) Exponentially stabilizing
  continuous-time controllers for periodic orbits of hybrid systems:
  Application to bipedal locomotion with ground height variations.
\newblock \emph{The International Journal of Robotics Research} 35(8):
  977--999.
\newblock \doi{10.1177/0278364915593400}.

\bibitem[{Han and Tedrake(2017)}]{HanTedrake2017}
Han W and Tedrake R (2017) Feedback design for multi-contact push recovery via
  {LMI} approximation of the piecewise-affine quadratic regulator.
\newblock In: \emph{Under Review}.

\bibitem[{Hartley et~al.(2017)Hartley, Da and Grizzle}]{HartleyGrizzleCCTA2017}
Hartley R, Da X and Grizzle JW (2017) Stabilization of {3D} underactuated biped
  robots: Using posture adjustment and gait libraries to reject velocity
  disturbances.
\newblock In: \emph{{IEEE} Conference on Control Technology and Applications
  (CCTA)}.

\bibitem[{Hereid et~al.(2016{\natexlab{a}})Hereid, Cousineau, Hubicki and
  Ames}]{HeCoHuAm16}
Hereid A, Cousineau EA, Hubicki CM and Ames AD (2016{\natexlab{a}}) {3D}
  dynamic walking with underactuated humanoid robots: A direct collocation
  framework for optimizing hybrid zero dynamics.
\newblock In: \emph{IEEE International Conference on Robotics and Automation
  (ICRA)}.

\bibitem[{Hereid et~al.(2016{\natexlab{b}})Hereid, Kolathaya and
  Ames}]{Hereid2016Online}
Hereid A, Kolathaya S and Ames AD (2016{\natexlab{b}}) Online hybrid zero
  dynamics optimal gait generation using legendre pseudospectral optimization.
\newblock In: \emph{To appear in: IEEE Conference on Decision and Control
  (CDC)}. IEEE.

\bibitem[{Hereid et~al.(2014)Hereid, Kolathaya, Jones, Van~Why, Hurst and
  Ames}]{Hereid2014Dynamic}
Hereid A, Kolathaya S, Jones MS, Van~Why J, Hurst JW and Ames AD (2014) Dynamic
  multi-domain bipedal walking with {ATRIAS} through {SLIP} based
  human-inspired control.
\newblock In: \emph{Proceedings of the 17th International Conference on Hybrid
  Systems: Computation and Control}, HSCC '14. New York, NY, USA: ACM.
\newblock ISBN 978-1-4503-2732-9, pp. 263--272.
\newblock \doi{10.1145/2562059.2562143}.

\bibitem[{Isidori(1995)}]{ISI95}
Isidori A (1995) \emph{Nonlinear Control Systems}.
\newblock Third edition. Berlin: Springer-Verlag.

\bibitem[{Jones(2014)}]{Jo2014}
Jones MS (2014) \emph{Optimal control of an underactuated bipedal robot}.
\newblock Master's Thesis, Oregon State University, ScholarsArchive@OSU.

\bibitem[{Kajita and Tani(1991)}]{Kajita1991}
Kajita S and Tani K (1991) Study of dynamic biped locomotion on rugged
  terrain-theory and basic experiment.
\newblock Pisa, pp. 741--746.

\bibitem[{Kajita et~al.(1992)Kajita, Yamaura and Kobayashi}]{KAYAKO92}
Kajita S, Yamaura T and Kobayashi A (1992) Dynamic walking control of biped
  robot along a potential energy conserving orbit.
\newblock \emph{IEEE Transactions on Robotics and Automation} 8(4): 431--37.

\bibitem[{Kaneko et~al.(2011)Kaneko, Kanehiro, Morisawa, Akachi, Miyamori,
  Hayashi and Kanehira}]{robots_HRP4}
Kaneko K, Kanehiro F, Morisawa M, Akachi K, Miyamori G, Hayashi A and Kanehira
  N (2011) Humanoid robot {HRP-4} -- humanoid robotics platform with
  lightweight and slim body.
\newblock pp. 4400--4407.

\bibitem[{Karagiannis et~al.(2005)Karagiannis, Astolfi and
  Ortega}]{karagiannis2005nonlinear}
Karagiannis D, Astolfi A and Ortega R (2005) Nonlinear stabilization via system
  immersion and manifold invariance: survey and new results.
\newblock \emph{Multiscale Modeling \& Simulation} 3(4): 801--817.

\bibitem[{Khalil(2002)}]{KHA02}
Khalil H (2002) \emph{Nonlinear Systems - 3rd Edition}.
\newblock Upper Saddle River, NJ: Prentice Hall.

\bibitem[{Kolathaya and Ames(2017)}]{kolathaya2017parameter}
Kolathaya S and Ames AD (2017) Parameter to state stability of control lyapunov
  functions for hybrid system models of robots.
\newblock \emph{Nonlinear Analysis: Hybrid Systems} 25: 174--191.

\bibitem[{Koolen et~al.(2012)Koolen, de~Boer, Rebula, Goswami and
  Pratt}]{KODBREGOPR12}
Koolen T, de~Boer T, Rebula J, Goswami A and Pratt J (2012)
  {Capturability-based analysis and control of legged locomotion, Part 1:
  Theory and application to three simple gait models}.
\newblock \emph{The International Journal of Robotics Research} 31(9):
  1094--1113.
\newblock \doi{10.1177/0278364912452673}.

\bibitem[{Krause et~al.(2012)Krause, Englsberger, Wieber and Ott}]{KrEnWiOt12}
Krause M, Englsberger J, Wieber PB and Ott C (2012) Stabilization of the
  capture point dynamics for bipedal walking based on model predictive control.
\newblock \emph{IFAC Proceedings Volumes} 45(22): 165--171.

\bibitem[{Kuindersma et~al.(2016)Kuindersma, Deits, Fallon, Valenzuela, Dai,
  Permenter, Koolen, Marion and Tedrake}]{KuDeFa16}
Kuindersma S, Deits R, Fallon M, Valenzuela A, Dai H, Permenter F, Koolen T,
  Marion P and Tedrake R (2016) Optimization-based locomotion planning,
  estimation, and control design for the atlas humanoid robot.
\newblock \emph{Autonomous Robots} 40(3): 429--455.

\bibitem[{Liu et~al.(2013)Liu, Atkeson and Su}]{LiAtSu13}
Liu C, Atkeson CG and Su J (2013) Biped walking control using a trajectory
  library.
\newblock \emph{Robotica} 31(02): 311--322.

\bibitem[{Manchester and Umenberger(2014)}]{MaUm14}
Manchester IR and Umenberger J (2014) Real-time planning with primitives for
  dynamic walking over uneven terrain.
\newblock In: \emph{Robotics and Automation (ICRA), 2014 IEEE International
  Conference on}. IEEE, pp. 4639--4646.

\bibitem[{Marcucci et~al.(2017)Marcucci, Deits, Gabiccini, Bicchi and
  Tedrake}]{MarcucciTedrake2017}
Marcucci T, Deits R, Gabiccini M, Bicchi A and Tedrake R (2017) Approximate
  hybrid model predictive control for multi-contact push recovery in complex
  environments.
\newblock In: \emph{Under Review}.

\bibitem[{Martin et~al.(2014)Martin, Post and Schmiedeler}]{Martin2014Design}
Martin AE, Post DC and Schmiedeler JP (2014) Design and experimental
  implementation of a hybrid zero dynamics-based controller for planar bipeds
  with curved feet.
\newblock \emph{The International Journal of Robotics Research} 33(7):
  988--1005.
\newblock \doi{10.1177/0278364914522141}.

\bibitem[{Mayne et~al.(2000)Mayne, Rawlings, Rao and
  Scokaert}]{mayne2000constrained}
Mayne DQ, Rawlings JB, Rao CV and Scokaert PO (2000) Constrained model
  predictive control: Stability and optimality.
\newblock \emph{Automatica} 36(6): 789--814.

\bibitem[{Morris and Grizzle(2009)}]{MOGR08}
Morris B and Grizzle JW (2009) Hybrid invariant manifolds in systems with
  impulse effects with application to periodic locomotion in bipedal robots.
\newblock \emph{IEEE Transactions on Automatic Control} 54(8): 1751--1764.

\bibitem[{Nguyen et~al.(2016)Nguyen, Agrawal, Da, Martin, Geyer, Grizzle and
  Sreenath}]{NguyenWAFR2017}
Nguyen Q, Agrawal A, Da X, Martin WC, Geyer H, Grizzle JW and Sreenath K (2016)
  Dynamic walking on stepping stones with gait library and control barrier.
\newblock In: \emph{Workshop on Algorithimic Foundations of Robotics (WAFR)}.

\bibitem[{Nguyen et~al.(2017)Nguyen, Agrawal, Da, Martin, Geyer, Grizzle and
  Sreenath}]{NguyenRSS2017}
Nguyen Q, Agrawal A, Da X, Martin WC, Geyer H, Grizzle JW and Sreenath K (2017)
  Dynamic walking on randomly-varying discrete terrain with one-step preview.
\newblock In: \emph{Robotics: Science and Systems (RSS)}.

\bibitem[{Ogura et~al.(2006)Ogura, Aikawa, Shimomura, Morishima, Lim and
  Takanishi}]{Takanishi06_WABIAN2}
Ogura Y, Aikawa H, Shimomura K, Morishima A, Lim H and Takanishi A (2006)
  Development of a new humanoid robot {WABIAN-2}.
\newblock In: \emph{Proc. of the 2006 IEEE International Conference on Robotics
  and Automation, Orlando, FL}. pp. 76--81.

\bibitem[{Parisini and Zoppoli(1995)}]{parisini1995receding}
Parisini T and Zoppoli R (1995) A receding-horizon regulator for nonlinear
  systems and a neural approximation.
\newblock \emph{Automatica} 31(10): 1443--1451.

\bibitem[{Park et~al.(2013)Park, Ramezani and Grizzle}]{PaRaGr13}
Park H, Ramezani A and Grizzle JW (2013) A finite-state machine for
  accommodating unexpected large ground height variations in bipedal robot
  walking.
\newblock \emph{IEEE Transactions on Robotics} 29(29): 331--345.

\bibitem[{Park et~al.(2005)Park, Kim, Lee and Oh}]{robots_KHR3}
Park I, Kim J, Lee J and Oh J (2005) Mechanical design of humanoid robot
  platform {KHR-3} {(KAIST Humanoid Robot 3: HUBO)}.
\newblock In: \emph{Proc. 5th IEEE---RAS Int. Conf. Humanoid Robots}. pp.
  321--326.

\bibitem[{Parker and Chua(1989)}]{PACH89}
Parker TS and Chua LO (1989) \emph{Practical Numerical Algorithms for Chaotic
  Systems}.
\newblock New York, NY: Springer-Verlag.

\bibitem[{Pfeiffer et~al.(2002)Pfeiffer, Loffler and Gienger}]{PFLOGI02}
Pfeiffer F, Loffler K and Gienger M (2002) The concept of jogging {J}ohnnie.
\newblock In: \emph{IEEE International Conference on Robotics and Automation}.
  Washington, DC.

\bibitem[{Powell et~al.(2013)Powell, Hereid and Ames}]{PoHeAm13}
Powell MJ, Hereid A and Ames AD (2013) Speed regulation in {3D} robotic walking
  through motion transitions between human-inspired partial hybrid zero
  dynamics.
\newblock In: \emph{Robotics and Automation (ICRA), 2013 IEEE International
  Conference on}. IEEE, pp. 4803--4810.

\bibitem[{Pratt et~al.(2012)Pratt, Koolen, de~Boer, Rebula, Cotton, Carff,
  Johnson and Neuhaus}]{PRKODBRECOJONE12}
Pratt J, Koolen T, de~Boer T, Rebula J, Cotton S, Carff J, Johnson M and
  Neuhaus P (2012) {Capturability-based analysis and control of legged
  locomotion, Part 2: Application to M2V2, a lower-body humanoid}.
\newblock \emph{The International Journal of Robotics Research} 31(10):
  1117--1133.
\newblock \doi{10.1177/0278364912452762}.

\bibitem[{Pratt and Tedrake(2006)}]{PRTE06}
Pratt J and Tedrake R (2006) Velocity-based stability margins for fast bipedal
  walking.
\newblock In: Diehl M and Mombaur K (eds.) \emph{Fast Motions in Biomechanics
  and Robotics}, \emph{Lecture Notes in Control and Information Sciences},
  volume 340. Springer Berlin Heidelberg.
\newblock ISBN 978-3-540-36118-3, pp. 299--324.
\newblock \doi{10.1007/978-3-540-36119-0_14}.

\bibitem[{Raibert(1986{\natexlab{a}})}]{RA86b}
Raibert MH (1986{\natexlab{a}}) Legged robots.
\newblock \emph{Communications of the ACM} 29(6): 499--514.

\bibitem[{Raibert(1986{\natexlab{b}})}]{RA86a}
Raibert MH (1986{\natexlab{b}}) \emph{Legged robots that balance}.
\newblock Mass.: MIT Press.

\bibitem[{Ramezani et~al.(2014)Ramezani, Hurst, {Akbari Hamed} and
  Grizzle}]{RAHUAKGR14}
Ramezani A, Hurst JW, {Akbari Hamed} K and Grizzle JW (2014) {Performance
  Analysis and Feedback Control of ATRIAS, A Three-Dimensional Bipedal Robot}.
\newblock \emph{Journal of Dynamic Systems, Measurement, and Control} 136(2).

\bibitem[{Reher et~al.(2016{\natexlab{a}})Reher, Cousineau, Hereid, Hubicki and
  Ames}]{ReCoHeHuAm16}
Reher J, Cousineau EA, Hereid A, Hubicki CM and Ames AD (2016{\natexlab{a}})
  Realizing dynamic and efficient bipedal locomotion on the humanoid robot
  {DURUS}.
\newblock In: \emph{2016 IEEE International Conference on Robotics and
  Automation (ICRA)}. pp. 1794--1801.
\newblock \doi{10.1109/ICRA.2016.7487325}.

\bibitem[{Reher et~al.(2016{\natexlab{b}})Reher, Hereid, Kolathaya, Hubicki and
  Ames}]{Reher2016Algorithmic}
Reher JP, Hereid A, Kolathaya S, Hubicki CM and Ames AD (2016{\natexlab{b}})
  Algorithmic foundations of realizing multi-contact locomotion on the humanoid
  robot {DURUS}.
\newblock In: \emph{The International Workshop on the Algorithmic Foundations
  of Robotics (WAFR)}.

\bibitem[{Reyhanoglu et~al.(1999)Reyhanoglu, van~der Schaft, McClamroch and
  Kolmanovsky}]{REVAMCKO99}
Reyhanoglu M, van~der Schaft A, McClamroch N and Kolmanovsky I (1999) Dynamics
  and control of a class of underactuated mechanical systems.
\newblock \emph{IEEE Transactions on Automatic Control} 44(9): 1663--1671.

\bibitem[{Rezazadeh et~al.(2015)Rezazadeh, Hubicki, Jones, Peekema, Van~Why,
  Abate and Hurst}]{ReHuJoPeVaJoAbHu15}
Rezazadeh S, Hubicki C, Jones M, Peekema A, Van~Why J, Abate A and Hurst JW
  (2015) Spring-mass walking with {ATRIAS} in {3D}: Robust gait control
  spanning zero to 4.3 kph on a heavily underactuated bipedal robot.
\newblock \emph{ASME Dynamic Systems and Control Conference (DSCC)} : 23.

\bibitem[{Saglam and Byl(2015)}]{SaBy15}
Saglam CO and Byl K (2015) Meshing hybrid zero dynamics for rough terrain
  walking.
\newblock In: \emph{2015 IEEE International Conference on Robotics and
  Automation (ICRA)}. IEEE, pp. 5718--5725.

\bibitem[{Sakagami et~al.(2002)Sakagami, Watanabe, Aoyama, Matsunaga, Higaki
  and Fujimura}]{SAWAAO02}
Sakagami Y, Watanabe R, Aoyama C, Matsunaga S, Higaki N and Fujimura K (2002)
  The intelligent {ASIMO}: system overview and integration.
\newblock In: \emph{IEEE/RSJ International Conference on Intelligent Robots and
  Systems (IROS)}. Lausanne, Switzerland, pp. 2478--83.

\bibitem[{Sch{\"u}rmann and Althoff(2017)}]{schurmann2017convex}
Sch{\"u}rmann B and Althoff M (2017) Convex interpolation control with formal
  guarantees for disturbed and constrained nonlinear systems.
\newblock In: \emph{Proceedings of the 20th International Conference on Hybrid
  Systems: Computation and Control}. ACM, pp. 121--130.

\bibitem[{Shim and Teel(2003)}]{shim2003asymptotic}
Shim H and Teel AR (2003) Asymptotic controllability and observability imply
  semiglobal practical asymptotic stabilizability by sampled-data output
  feedback.
\newblock \emph{Automatica} 39(3): 441--454.

\bibitem[{Spong(1994)}]{SP94}
Spong MW (1994) Partial feedback linearization of underactuated mechanical
  systems.
\newblock In: \emph{Proc. of the IEEE/RSJ International Conference on
  Intelligent Robots and Systems, Munich, Germany}. pp. 314--321.

\bibitem[{Spong(1996)}]{SPO96}
Spong MW (1996) Energy based control of a class of underactuated mechanical
  systems.
\newblock In: \emph{Proc. of IFAC World Congress, San Francisco, CA}. pp.
  431--435.

\bibitem[{Sreenath et~al.(2012)Sreenath, Park and Grizzle}]{Sreenath2012Design}
Sreenath K, Park H and Grizzle JW (2012) Design and experimental implementation
  of a compliant hybrid zero dynamics controller with active force control for
  running on {MABEL}.
\newblock In: \emph{{IEEE} International Conference on Robotics and Automation,
  {ICRA} 2012, 14-18 May, 2012, St. Paul, Minnesota, {USA}}. pp. 51--56.
\newblock \doi{10.1109/ICRA.2012.6224944}.

\bibitem[{Sreenath et~al.(2011)Sreenath, Park, Poulakakis and
  Grizzle}]{KOPAPOG210}
Sreenath K, Park H, Poulakakis I and Grizzle J (2011) A compliant hybrid zero
  dynamics controller for stable, efficient and fast bipedal walking on
  {MABEL}.
\newblock \emph{International Journal of Robotics Research} 30(9): 1170--1193.

\bibitem[{Sreenath et~al.(2013)Sreenath, Park, Poulakakis and
  Grizzle}]{SrPaPoGr13}
Sreenath K, Park HW, Poulakakis I and Grizzle J (2013) Embedding active force
  control within the compliant hybrid zero dynamics to achieve stable, fast
  running on {MABEL}.
\newblock \emph{The International Journal of Robotics Research} 32(3):
  324--345.

\bibitem[{Vidyasagar(2002)}]{vidyasagar2002nonlinear}
Vidyasagar M (2002) \emph{Nonlinear systems analysis}.
\newblock SIAM.

\bibitem[{Wang et~al.(2017)Wang, Forni, Ortega, Liu and Su}]{wang2017immersion}
Wang L, Forni F, Ortega R, Liu Z and Su H (2017) Immersion and invariance
  stabilization of nonlinear systems via virtual and horizontal contraction.
\newblock \emph{IEEE Transactions on Automatic Control} 62(8): 4017--4022.

\bibitem[{Wang and Chevallereau(2011)}]{WANGChevallereau2011}
Wang T and Chevallereau C (2011) Stability analysis and time-varying walking
  control for an under-actuated planar biped robot.
\newblock \emph{Robotics and Autonomous Systems} 59(6): 444 -- 456.
\newblock \doi{https://doi.org/10.1016/j.robot.2011.03.002}.

\bibitem[{Westervelt et~al.(2003)Westervelt, Grizzle and Koditschek}]{WEGRKO03}
Westervelt E, Grizzle JW and Koditschek D (2003) Hybrid zero dynamics of planar
  biped walkers.
\newblock \emph{IEEE Transactions on Automatic Control} 48(1): 42--56.

\bibitem[{Westervelt et~al.(2007)Westervelt, Grizzle, Chevallereau, Choi and
  Morris}]{WGCCM07}
Westervelt ER, Grizzle JW, Chevallereau C, Choi J and Morris B (2007)
  \emph{Feedback Control of Dynamic Bipedal Robot Locomotion}.
\newblock Control and Automation. Boca Raton, FL: CRC Press.

\bibitem[{Westervelt et~al.(2002)Westervelt, Grizzle and Koditschek}]{WEGRKO02}
Westervelt ER, Grizzle JW and Koditschek DE (2002) Zero dynamics of planar
  biped walkers with one degree of under actuation.
\newblock In: \emph{IFAC World Congress}. Barcelona, Spain.

\bibitem[{Yamaguchi et~al.(1999)Yamaguchi, Soga, Inoue and
  Takanishi}]{Yamaguchi99}
Yamaguchi J, Soga E, Inoue S and Takanishi A (1999) Development of a bipedal
  humanoid robot-control method of whole body cooperative dynamic biped
  walking.
\newblock In: \emph{Proc. of the 1999 IEEE International Conference on Robotics
  and Automation, Detroit, MI}. pp. 368--374.

\bibitem[{Zaytsev et~al.(2015)Zaytsev, Hasaneini and Ruina}]{ZaHaRu15}
Zaytsev P, Hasaneini SJ and Ruina A (2015) Two steps is enough: No need to plan
  far ahead for walking balance.
\newblock In: \emph{2015 IEEE International Conference on Robotics and
  Automation (ICRA)}. pp. 6295--6300.
\newblock \doi{10.1109/ICRA.2015.7140083}.

\bibitem[{Zhao et~al.(2015)Zhao, Horn, Reher, Paredes and
  Ames}]{Zhao2015hybrid}
Zhao H, Horn J, Reher J, Paredes V and Ames AD (2015) A hybrid systems and
  optimization-based control approach to realizing multi-contact locomotion on
  transfemoral prostheses.
\newblock In: \emph{{IEEE} Conference on Decision and Control (CDC)}. pp.
  1607--1612.
\newblock \doi{10.1109/CDC.2015.7402440}.

\end{thebibliography}

\appendices

\appendix

\section{Index to Multimedia Extensions}
\label{sec:multimedia}

\begin{table}[h]
	\label{tab:multimediaExtenstion}
	{
		\renewcommand{\arraystretch}{1.5}%
		\begin{center}
			\begin{tabular}{c c c}
			Extension & Type & Description \\
				\hline
		 	  \hypertarget{ex:code}{1} & Code & \href{https://www.dropbox.com/s/mr6ctpkzho688fe/Pendulum-on-a-cartMATLABcode.zip?dl=0}{Inverted pendulum example} \\
			  \hypertarget{ex:video}{2} & Video & \href{https://youtu.be/3gGEH9qaXXM}{3D robot experiments}

			\end{tabular}
		\end{center}
	}
\end{table}

\section{Normal Forms for Mechanical Models}
\label{sec:MechanicalModelDecomposition}

Consider a standard mechanical model
$$D(q) \ddot{q} + C(q, \dot{q}) \dot{q} + G(q) = B u$$
and let $$\Omega(q,\dot{q}) := C(q,\dot{q})\dot{q} + G(q).$$
We suppose the system is \textit{underactuated}, that is, there are fewer independent actuators than generalized coordinates. In fact, we suppose there exists a partition of the coordinates in which the model takes the form
\begin{equation}
\label{eqn:MechanicalModelpartitioned}
\begin{aligned}
 D_{11}(q) \ddot{q}_1+ D_{12}(q) \ddot{q}_2 + \Omega_1(q,\dot{q}) & = 0\\
 D_{21}(q) \ddot{q}_1+ D_{22}(q) \ddot{q}_2 + \Omega_2(q,\dot{q}) & = B_2 u,
\end{aligned}
\end{equation}
with $B_2$ square and invertible. Because $D(q)$ is positive definite, by the Shur Complement Lemma,  $D_{11}(q)$, $D_{22}(q)$, and
\begin{equation}
 \bar{D}(q):=D_{22}(q) - D_{21}(q) D_{11}^{-1}(q)D_{12} (q)
 \end{equation}
 are all positive definite as well.

Following \cite{SP94}, define
\begin{equation} \label{eqn:normal_form_terms}
\begin{aligned}
J^{\rm norm}(q) & := D_{11}^{-1}(q)D_{12}(q)\\
 \bar{\Omega}_1(q,\dot{q}) & := -D_{11}^{-1}(q)\Omega_1(q,\dot{q}) \\
   \bar{\Omega}_2(q,\dot{q}) & := \Omega_2(q,\dot{q}) - D_{21}(q) D_{11}^{-1}(q)\Omega_1(q,\dot{q}),
\end{aligned}
\end{equation}
 Then the (regular) feedback
\begin{equation}\label{eqn:preliminary_state_feedback}
u = B_2^{-1}(q)\left(\bar{D}(q) v + \bar{\Omega}_2(q,\dot{q}) \right),
\end{equation}
results in the Spong normal form:
\begin{equation} \label{eqn:mech_model_spong_normal_form}
\begin{aligned}
\ddot{q}_1 & = \bar{\Omega}_1(q,\dot{q})- {J^{\rm norm}}(q)v, \\
\ddot{q}_2 & = v.
\end{aligned}
\end{equation}
Defining $x_1 = (q_1, \dot{q}_1)$, $x_{2a} = q_2$, and $\dot{x}_{2b}= \dot{q}_2$, it follows that the model can be expressed as
\begin{equation}
\label{eq:SpongStateForm}
\begin{aligned}
\dot{x}_1 &= f_1(x_1,x_2,v) \\
\dot{x}_{2a} &= x_{2b} \\
\dot{x}_{2b} & = v.
\end{aligned}
\end{equation}

An alternative form is developed in \cite[pp.~62]{WGCCM07}, which uses the conjugate momenta that arises from Lagrange's equations. It has the advantage that the input only appears in the second row of the model. Define
\begin{align}
\sigma_1 &:= D_{11}(q) \dot{q}_1 + D_{12}(q) \dot{q}_2\\
\dot{D}(q,\dot{q})&:=\frac{d}{dt}D(q)
\end{align}
Then the model can also be written as
\begin{equation} \label{eqn:mech_model_mixed_normal_form}
\begin{aligned}
\dot{q}_1 & = D_{11}^{-1}(q) \left[ \sigma_1 - D_{12}(q) \dot{q}_2 \right]\\
\dot{\sigma}_1& = \kappa_1(q, \sigma_1, \dot{q}_2) \\
\ddot{q}_2 & = v
\end{aligned}
\end{equation}
where
\begin{equation} \label{eqn:mech_model_mixed_normal_form_detail}
\begin{aligned}
\kappa_1(q, \sigma_1, \dot{q}_2)&:=  \Big( \dot{D}_{11}(q, \dot{q}) \dot{q}_1 +  \dot{D}_{12}(q, \dot{q}) \dot{q}_2- \\
  &~~~~~~ \left.  \Omega_1(q, \dot{q}) \Big)\right|~_{\dot{q}_1 = D_{11}^{-1}(q) \left[ \sigma_1 - D_{12}(q) \dot{q}_2 \right]}.
 \end{aligned}
\end{equation}
With $x_2$ defined as above and $x_1:=(q_1,\sigma_1)$, the model takes the form
\begin{equation}
\label{eq:MixedStateForm}
\begin{aligned}
\dot{x}_1 &= f_1(x_1,x_2) \\
\dot{x}_{2a} &= x_{2b} \\
\dot{x}_{2b} & = v.
\end{aligned}
\end{equation} 
Various authors prefer one of \eqref{eq:SpongStateForm} and \eqref{eq:MixedStateForm} to the other; both are useful. 

\section{Proofs}
\label{sec:Proofs}

\subsection{Proof of Proposition~\ref{prop:XDaMPCStability}}

The proof is most easily done using the method of Poincar\'e sections \cite{PACH89}. By A-\ref{itm:A1} and the assumption on $u^{cu}$, the closed-loop system \eqref{eq:ClosedLoopCUFullStateFeedback} has period $T_p>0$ and the origin is an equilibrium point. Due to the time-varying nature of the closed-loop system, we make time a state variable, and because the system is $T_p$-periodic, we add in a time-based reset map
\begin{equation}
\label{eq:hybridVersionPeriodicallyTimeVarying}
\left\{
\begin{aligned}
\dot{\tau}&=1, & \tau^- <T_p\\
\dot{x}&=f^{cu}(\tau,x):=f(\tau,x,u^{cu}(\tau,x)),  \\
\tau^{+}&=0 & \tau^- =T_p\\
x^{+}&=x^{-}.
\end{aligned}\right.
\end{equation}
The notation $\tau^-$, $\tau^+$, $x^-$ and $x^+$ is explained in Section~\ref{sec:HybridSystemControl}. Because the state reset map is trivial, namely $x^{+}=x^{-}$, the solutions of \eqref{eq:ClosedLoopCUFullStateFeedback} and \eqref{eq:hybridVersionPeriodicallyTimeVarying} are identical. Define a Poincar\'e section by
\begin{equation}
\label{eq:PoincareSectionForTVProof}
{\cal S}_n := \{ (\tau, x) \in \real^{n+1}~|~\tau=T_p, x\in B \},
\end{equation}
which is an $n$-dimensional hypersurface in the state space of the model.
Then, by construction of the closed-loop system, for $\xi \in {\cal S}_n$, the Poincar\'e map $P:{\cal S}_n \to {\cal S}_n$ is given by
\begin{equation}
\label{eq:PoincareMapForTVProof}
P(\xi) = \varphi_\xi(T_p).
\end{equation}
Indeed, for $\xi \in B$,
$$\begin{aligned}
\varphi_\xi(t) &= \xi + \int_0^t f(\varphi_\xi(s), u_\xi(s)) ds \\
&= \xi + \int_0^t f(\varphi_\xi(s), u^{cu}(s,\varphi_\xi(s)) ds,
\end{aligned}$$
due to \eqref{eq:LearningConditionA}.
By A-\ref{itm:A4}, $\xi^*=0$ is a fixed point of the Poincar\'e map. Also by A-\ref{itm:A4}, $P$ is a contraction because for each $\xi \in B$, $V\circ P(\xi) \le c V(\xi)$, and hence by induction,
$$V\circ P^k(\xi) \le c^k V(\xi),$$
and by A-\ref{itm:A3},
$$||P^k(\xi)||^2 \le c^k \frac{\alpha_2}{\alpha_1 } ||\xi||^2 \xrightarrow[k \to \infty]{} 0,$$
proving local exponential stability of the fixed point. The uniformity in $t_0$ follows from periodicity.
\mbox{ } \hfill \(\blacksquare\)

%

\subsection{Proof of Proposition~\ref{prop:ReducedSystemStability}}

Without loss generality, we assume that $B_{1cu}$ is bounded so that its closure is compact. Then there exists $L$, a Lipschitz constant, such that
$$||\gamma(x_1)||_2 \le L ||x_1||_2$$
for all $x_1 \in B_{1cu}$. Define $V_1: B_1 \to \real$ by $V_1(x_1):= V(x_1, \gamma(x_1))$. It follows that
$$\alpha_1 x_1^\top x_1 \le V_1(x_1) \le \alpha_2 (1+L^2) x_1^\top x_1, $$
and hence $V_1$ is positive definite, with quadratic lower and upper bounds. From \eqref{eq:LyapunovConstraintReduced},
$$ V_1(\varphi_{1\xi_1}(T_p)) \le c V_1(\xi_1).$$
From here, the proof of Prop.~\ref{prop:XDaMPCStability} can be repeated and the result follows.
\mbox{ } \hfill \(\blacksquare\)

\subsection{Proof of Theorem.~\ref{them:OverallSystemStability}}

From the hypotheses of the Theorem and Prop.~\ref{prop:ReducedSystemStability}, the closed-loop system \eqref{eq:TwoBlockSystemNewCoordinates} is a cascade of two locally exponentially systems, namely, the second row of \eqref{eq:TwoBlockSystemNewCoordinates} and the reduced-order system \eqref{eq:ClosedLoopReducedFeedback}. By standard results, the overall system is locally exponentially stable. By \cite[Thm. 43, Section~5.1, pp.~143]{vidyasagar2002nonlinear}, because the system is periodically time-varying, the stability is uniform in $t_0$.

\mbox{ } \hfill \(\blacksquare\)

\subsection{Proof of Corollary.~\ref{cor:OverallSystemStabilityExtraU}}

Defining $y$ as in \eqref{eq:NewCoordinates} results in the closed-loop system having the form
\begin{equation}
\label{eq:TwoBlockSystemNewCoordinatesExtraUterm}
\begin{aligned}
\dot{x}_1&= f_1(t,x_1,\nu(t,x_1)+y, \mu(t,x_1)-[K_p~~K_d] y)\\
\dot{y}&= Ay,
\end{aligned}
\end{equation}
with $A$ Hurwitz. Hence, the proof of Theorem~\ref{them:OverallSystemStability} can be repeated and we are done.

\mbox{ } \hfill \(\blacksquare\)

\subsection{Proof of Proposition~\ref{prop:FullStateHybridModel} }

The proof is very similar to that of Prop.~\ref{prop:XDaMPCStability}.
Define a Poincar\'e section by
\begin{equation}
\label{eq:PoincareSectionFoHybridFullState}
{\cal S}_n := \{ (\tau, x) \in \real^{n+1}~|~\tau=T_p, x\in B \},
\end{equation}
which is an $n$-dimensional hypersurface in the state space of the model.
Then, by construction of the closed-loop system, for $\xi \in {\cal S}_n$, the Poincar\'e map $P:{\cal S}_n \to {\cal S}_n$ is given by
\begin{equation}
\label{eq:PoincareMapForHybridProof}
P(\xi) = \varphi_\xi(T_p).
\end{equation}
By H-\ref{itm:H4}, $\xi^*:=x_m^*$ is a fixed point of the Poincar\'e map. Also by H-\ref{itm:H4}, $P$ is a contraction because for each $\xi \in B$, $V\circ P(\xi) \le c V(\xi)$, and hence by induction,
$$V\circ P^k(\xi) \le c^k V(\xi),$$
and by H-\ref{itm:H3},
$$||P^k(\xi)||^2 \le c^k \frac{\alpha_2}{\alpha_1 } ||\xi-\xi^*||^2 \xrightarrow[k \to \infty]{} 0,$$
proving local exponential stability of the fixed point. Because the closed-loop system is locally Lipschitz continuous, local exponential stability of the fixed point implies exponential attractivity of the orbit
$${\cal O}:=\{ (\tau, \varphi_{\xi^*}(\tau)~|~0\le \tau < T_p \}. $$
Because $\tau(t)=t$, we have local exponential stability of the periodic solution.
\mbox{ } \hfill \(\blacksquare\)

\subsection{Proof of Proposition~\ref{prop:HybridZeroDynamicsStability} }

Without loss generality, we assume that $B_{1cu}$ is bounded so that its closure is compact. Then there exists $L$, a Lipschitz constant, such that
$$||\gamma(x_1-x_1^*)||_2 \le L ||x_1-x_1^*||_2$$
for all $x_1 \in B_{1cu}$. Define $V_1: B_1 \to \real$ by $V_1(x_1):= V(x_1, \gamma(x_1))$. It follows that

$$ \begin{array}{r}
\alpha_1 (x_1-x_1^*)^\top  (x_1-x_1^*) \le V_1(x_1)  \le \\
\alpha_2 (1+L^2) (x_1-x_1^*)^\top (x_1-x_1^*),
\end{array}
$$
and hence $V_1$ is positive definite, with quadratic lower and upper bounds. From \eqref{eq:LyapunovConstraintReduced},
$$ V_1(\varphi_{1\xi_1}(T_p)) \le c V_1(\xi_1).$$
From here, the proof of Prop.~\ref{prop:FullStateHybridModel} can be repeated and the result follows.
\mbox{ } \hfill \(\blacksquare\)

\subsection{Proof of Theorem.~\ref{them:StabilityThanksToHybridZeroDynamics} }

The Poincar\'e section is defined as in \eqref{eq:PoincareSectionFoHybridFullState}.  References \cite{grizzle2014models} and \cite[Chap.~4]{WGCCM07} show how to reduce the stability analysis of a hybrid model with two continuous phases to that of an equivalent hybrid system with a single continuous phase. We build the equivalent hybrid system with the continuous phase from $\Sigma_{\rm ii}$ and a reset map $\Delta_{eq}$ that captures the flow of $\Sigma_{\rm i}$, viz
\begin{equation}
\label{eq:HybridEquivalent}
\left\{
\begin{aligned}
\dot{\tau}&=1 \hspace*{4.45 cm}  \tau^- < T_p\\
\dot{x}&=f(x_1,x_2,u_1,u_2) \\
u_1&= \mu_{1ii}(\tau,x_1) \\
u_2&= \mu_{2ii}(\tau,x_1) - \left[\frac{K_p}{\epsilon^2}~~\frac{K_d}{\epsilon}\right]\big(x_2 - \nu_{ii}(\tau, x_1) \big) \\
\left[\begin{array}{c} \tau \\ x \end{array} \right]^{+} &= \Delta_{eq}(\tau^-,x^-),
  \hspace*{2.45 cm}  \tau^- =  T_p
\end{aligned}\right.
\end{equation}
With this construction and Prop.~\ref{prop:HybridZeroDynamicsStability}, the zero dynamics manifold is
 $${\cal Z}:=\{(\tau,x_1,x_2)~|~ x_2=\nu_{ii}(\tau, x_1) \}, $$
and the restricted Poincar\'e map $\rho:{\cal S}_e \cap {\cal Z} \to {\cal S}_e \cap {\cal Z}$ has $x_1^*$ as a locally exponentially stable fixed point. The equivalent hybrid system \eqref{eq:HybridEquivalent} therefore satisfies all the hypotheses of  \cite[Thm.~2]{AmGaGrSr2014}, and hence the periodic orbit
$${\cal O}:=\{ (\tau, \varphi_{\xi_1^*}(\tau)~|~0\le \tau < T_p \} $$
is locally exponentially stable.
\mbox{ } \hfill \(\blacksquare\)

\section{Relation to Backstepping, Zero Dynamics, and Immersion and Invariance}
\label{app:Relations}

For definiteness, consider a standard Lagrangian dynamical model where $ q \in \real^n $ is a set of generalized coordinates  and $ u \in \real^m $ is a vector of torques,
\begin{equation}
	\label{eq:Lagrange}
	D(q)\ddot{q}+H(q, \dot{q})= B(q).
\end{equation}
Assume the system is underactuated, that is, $ n > m $, and that the coordinates have been decomposed as
\begin{equation}
\label{eq:qDecomposed}
q := \begin{bmatrix} q_1 \\ q_2 \end{bmatrix},
\end{equation}
in which the model takes the form
\begin{equation}\label{eq:Lagrangemodel_partitioned}
\begin{aligned}
 D_{11}(q) \ddot{q}_1+ D_{12}(q) \ddot{q}_2+ H_1(q,\dot{q})& = 0 \\
 D_{21}(q) \ddot{q}_1+ D_{22}(q) \ddot{q}_2+ H_2(q,\dot{q}) & = B_2(q)u.
\end{aligned}
\end{equation}
with $B_2(q)$ square and full rank.
References \cite{SPO96,REVAMCKO99,WEGRKO02} show that there is a regular feedback that places the system in the form
\begin{equation}
\label{eq:spong_modified_SS}
\begin{array}{ccl} \medskip
 \dot{x}_1 & = & f_1(x_1,q_2, \dot{q}_2) \\ \medskip
 \ddot{q}_2 & = & v,\\
\end{array}
\end{equation}
with
$$ x_1 := \begin{bmatrix} q_1  \\ \sigma_1 \end{bmatrix}~\text{and}~\sigma_1 := D_{11}(q) \dot{q}_1+ D_{12}(q) \dot{q}_2,$$
the generalized momentum conjugate to $q_1$.

\subsection{Backstepping}
To begin the backstepping process in \eqref{eq:spong_modified_SS}, one needs a feedback
$$ \begin{bmatrix} q_2  \\ \dot{q}_2 \end{bmatrix} =  \begin{bmatrix} \nu_a(x_1)  \\ \nu_b(x_1) \end{bmatrix}$$
that renders the origin of the reduced-order system
\begin{equation}
\label{eq:SpongNormalForm}
\dot{x}_1=f_1(x_1,\nu_a(x_1), \nu_b(x_1))
\end{equation}
locally exponentially stable with a known Lyapunov function. However, to pull this feedback through the double integrator, it must be true that $\nu_b(x_1)=\frac{d}{dt}\nu_a(x_1)$, that is
\begin{equation}
\label{eq:IntegrabilityCondition}
\nu_b(x_1)= \left[\frac{\partial}{\partial x_1} \nu_a(x_1) \right]f_1(x_1,\nu_a(x_1), \nu_b(x_1)).
\end{equation}
Backstepping does not provide any systematic means to meet the required integrability condition. Of course, if the system has the form
\begin{equation}
\label{eq:spong_modified_SS2}
\begin{array}{ccl} \medskip
 \dot{x}_1 & = & f_1(x_1,q_2)\\ \medskip
 \ddot{q}_2 & = & v,\\
\end{array}
\end{equation}
then there is no integrability constraint and backstepping can be done, assuming one is clever enough to find a feedback $q_2 = \nu_a(x_1)$ that renders the origin of
 $$\dot{x_1}=f_1(x_1,\nu_a(x_1))$$
locally exponentially stable. The solution we presented in Sect.~\ref{sec:InnerOuterLoop} uses trajectory optimization to automatically build a feedback that satisfies the integrability condition \eqref{eq:IntegrabilityCondition} and provides for local exponential stability. Moreover, bounds on inputs and other constraints can potentially be included in the trajectory optimization process, whereas they are challenging to incorporate into backstepping.

\subsection{Zero Dynamics}
\label{app:zeroDynamics}
The method of Hybrid Zero Dynamics as developed in \cite{WEGRKO02} exploits the structure of $f_1$ in \eqref{eq:SpongNormalForm}, namely
$$\frac{d}{dt}\begin{bmatrix} q_1 \\ \sigma_1 \end{bmatrix} = \left[ \begin{array}{c} \dot{q}_1  \\ f_{1b}(q_1,\dot{q}_1,q_2, \dot{q}_2)
\end{array} \right]$$
and
$$ \dot{q}_1 = D_{11}^{-1}(q) \left[\sigma_1  - D_{12}(q) \dot{q}_2 \right],$$
to solve for a solution of the form
$$ \begin{bmatrix} q_2  \\ \dot{q}_2 \end{bmatrix} = \left[ \begin{array}{c} h_d(q_1) \\ \left( \frac{\partial}{\partial q_1} h_d(q_1)\right) \dot{q}_1
\end{array} \right],$$
so that the integrability condition \eqref{eq:IntegrabilityCondition} is automatically met. When the computational method in \cite{Jo2014, HeCoHuAm16} does produce a solution, it does not come with a Lyapunov function and hence input-output linearization is often used to ``pull'' the virtual constraints back through the double integrators. Moreover, conditions for the \textit{virtual constraints} $q_2=h_d(q_1)$ to stabilize (a hybrid version of) \eqref{eq:SpongNormalForm} are only known when $q_1$ is a scalar. The solution we have given in Sect.~\ref{sec:HybridSystemControl} works for vector valued $q_1$, hence for models with more than one degree of underactuation. Moreover, even for one degree of underactuation, it provides a more general solution to the boundary value problem in that it naturally produces solutions of the form $h_d(t, q_1, \dot{q}_1)$, that is, the controller depends in a non-trivial way on the full state of the $x_1$-subsystem.

\subsection{Immersion and Invariance} The method of immersion and invariance (I\&I) presented in \cite{ASOR03,karagiannis2005nonlinear,wang2017immersion} is more general than backstepping and can provide alternative cascade realizations to the simple one used in \eqref{eq:TwoBlockSystemNewCoordinates}. However, I\&I still requires a target system to be provided, such as, \eqref{eq:ClosedLoopReducedFeedback}, which is what our method is constructing. In other words, once a feedback satisfying Prop.~\ref{prop:ReducedSystemStability} has been constructed, I\&I can be used to build alternatives to the feedback used in Thm.~\ref{them:OverallSystemStability}, but it will not replace the design of the reduced-order model.

%
%

\section{Standard MPC and Relative to This Work}
\label{app:Standard_MPC}

\begin{figure}[t!]
	\centering
	\includegraphics[width=0.9\columnwidth]{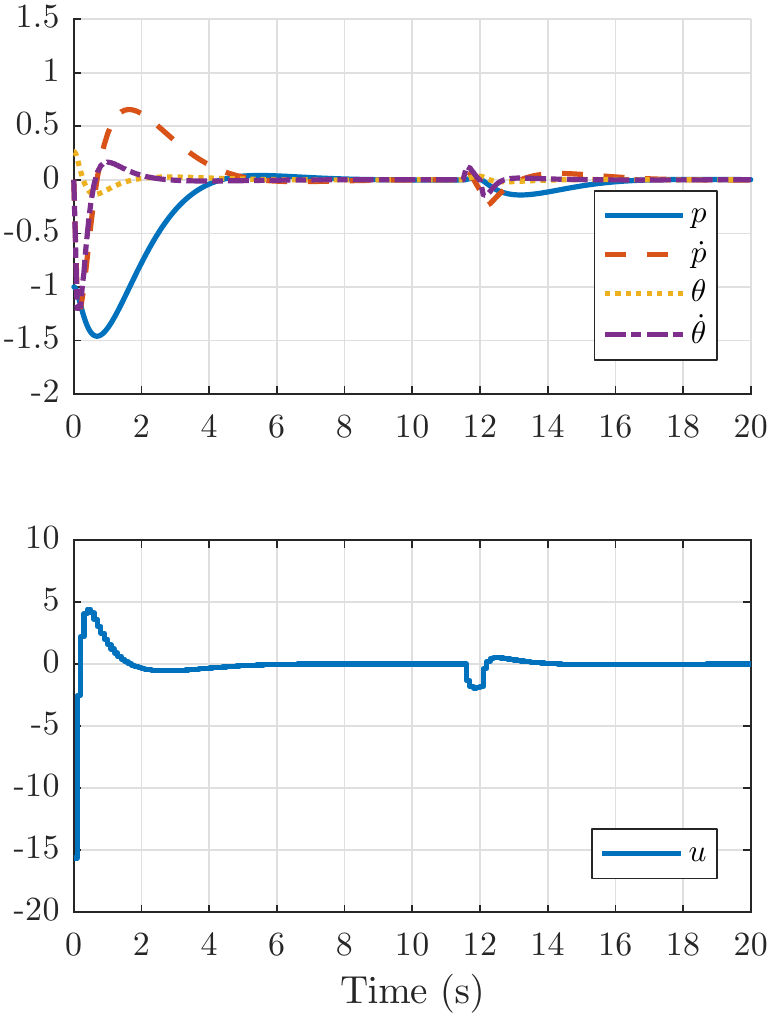}
	\caption{A classic MPC controller is applied to the same system as in Figure~\ref{fig:CH_CUcontrollerDiscussion} 
	}
	\label{fig:Standard_MPC}
\end{figure}

Figure~\ref{fig:Standard_MPC} shows the standard MPC with Zero-Order-Hold (ZOH) condition,
\begin{equation}
	u^{zoh}(t, \xi) = u_{\xi}(0), t \in [0, T_p).
\end{equation}
We would like to implement this controller to the high dimensional system. However, it does not scale well. The learned feedback $ \mu(t, x) $ shows the similar performance as the classic MPC. Sect.~\ref{sec:StabilizingReducedOrderModel} illustrates how to apply it to a reduced-order model while embedding it to the full-order model. In this sense, the curse of dimensionality has been mitigated.

\end{document}